\documentclass{wsbart}

\usepackage[utf8]{inputenc}
\usepackage[T1]{fontenc}

\usepackage[french,english]{babel}

\usepackage{amsmath, amssymb, amsthm, bbm}
\usepackage{mathtools}

\usepackage{enumitem}

\usepackage{orcidlink}

\usepackage[ruled,vlined]{algorithm2e}
\usepackage{caption,subcaption}
\usepackage{booktabs}

\usepackage[capitalise]{cleveref}
\crefname{equation}{}{}
\Crefname{equation}{}{}
\creflabelformat{equation}{(#2#1#3)}
\crefname{algocf}{Algorithm}{Algorithms}
\crefname{figure}{Figure}{Figures}
\crefname{prop}{Proposition}{Propositions}

\numberwithin{equation}{section}

\theoremstyle{plain}
\newtheorem{thm}{Theorem}[section]
\newtheorem{lem}[thm]{Lemma}
\newtheorem{cor}[thm]{Corollary}
\newtheorem{prop}[thm]{Proposition}

\theoremstyle{remark}
\newtheorem{rem}{Remark}[section]
\newtheorem{exm}{Example}

\theoremstyle{definition}
\newtheorem{defn}[thm]{Definition}

\crefname{thm}{Theorem}{Theorems}
\crefname{lem}{Lemma}{Lemmas}
\crefname{exm}{Example}{Examples}

\DeclareMathOperator{\Expect}{\mathbb E}
\DeclareMathOperator{\osc}{osc}
\DeclareMathOperator{\Law}{Law}
\DeclareMathOperator{\Ent}{Ent}
\DeclareMathOperator{\Var}{Var}
\DeclareMathOperator{\Tan}{Tan}

\newcommand{\1}{\mathbbm{1}}
\newcommand{\dd}{\mathop{}\!\mathrm{d}}
\renewcommand{\leq}{\leqslant}
\renewcommand{\geq}{\geqslant}

\title{Uniform-in-time propagation of chaos\\for mean field Langevin dynamics}

\author[1]{Fan Chen\,\orcidlink{0000-0003-0082-7908}}
\author[2]{Zhenjie Ren\,\orcidlink{0000-0003-4656-4074}}
\author[3]{Songbo Wang\,\orcidlink{0009-0009-3214-3587}}
\affil[1]{School of Mathematical Sciences, Shanghai Jiao Tong University,
Shanghai, China}
\affil[2]{CEREMADE, Université Paris-Dauphine, PSL, Paris, France}
\affil[3]{CMAP, École polytechnique, IP Paris, Palaiseau, France}

\subjclass{60J60, 60K35 (primary) 35B40, 35Q83, 35Q84 (secondary)}

\keywords{Langevin diffusion, Fokker–Planck equation, mean field interaction,
convergence to equilibrium, uniform-in-time propagation of chaos,
logarithmic Sobolev inequality, hypercontractivity,
Wasserstein distance, relative entropy}

\begin{document}

\maketitle

\begin{abstract}
We study the mean field Langevin dynamics and the associated particle system.
By assuming the functional convexity of the energy,
we obtain the \(L^p\)-convergence of the marginal distributions toward the
unique invariant measure for the mean field dynamics.
Furthermore, we prove the uniform-in-time propagation of chaos
in both the \(L^2\)-Wasserstein metric and relative entropy.
\end{abstract}

\begin{otherlanguage}{french}
\begin{abstract}
Nous étudions la dynamique de Langevin à champ moyen et le système de
particules correspondant.
En supposant la convexité fonctionnelle de l'énergie,
nous obtenons la convergence dans \(L^p\) des distributions marginales vers
l'unique mesure invariante pour la dynamique à champ moyen.
De plus, nous montrons la propagation du chaos uniforme en temps
à la fois dans la métrique de Wasserstein d'ordre \(2\)
et dans l'entropie relative.
\end{abstract}
\end{otherlanguage}

\tableofcontents

\section{Introduction}

\subsection{Preview of main results}
\label{sec:preview}
Let \(F : \mathcal P_2(\mathbb R^d) \to \mathbb R\) be a \emph{mean field
functional}
and \(D_m F\) be its intrinsic derivative.
In this paper, we study the long-time behavior of the following mean field
Langevin (MFL) dynamics:
\begin{equation}\label{eq:mfl-intro}
    \dd X_t = - D_m F(m_t, X_t) \dd t + \sqrt 2 \dd W_t,
\qquad \text{where $m_t = \Law(X_t)$,}
\end{equation}
as well as the corresponding dynamics of \(N\) particles:
\[
\dd  X^i_t = - D_m F(\mu_{\mathbf X_t}, X^i_t) \dd t + \sqrt 2 \dd W^i_t,\quad
i=1,\ldots,N,\quad\text{where $\mu_{\mathbf X_t} = \frac 1N\sum_{i=1}^N
\delta_{X^i_t}$.}
\]
Here, \(W_t, W^i_t\) are independent \(d\)-dimensional standard Brownian
motions.
We suppose that \(F\) is a functional such that
\begin{itemize}
\item the mapping \(m \mapsto F(m)\) is convex in the functional sense
(as opposed to the optimal transport sense);
\item for every \(x \in \mathbb R^d\), the mapping \(m \mapsto D_m F(m,x)\) is
\(M^F_{mm}\)-Lipschitz continuous with respect to the \(L^1\)-Wasserstein
metric;
\item for every \(m \in \mathcal P_2(\mathbb R^d)\),
the probability measure on \(\mathbb R^d\) that has density proportional to
\(x\mapsto\exp\bigl(-\frac{\delta F}{\delta m} (m,x)\bigr)\)
satisfies the \(\rho\)-logarithmic Sobolev inequality (LSI)
for some \(\rho > 0\).
\end{itemize}

Recently, there has been a growing interest
in modeling the training of neural networks
as a convex mean field optimization problem
(see
\cite{mei2019mean,chizat2018global,sirignano2022mean,rotskoff2018neural,%
HRSS19,kazeykina2020ergodicity, conforti2023game}
and also our \cref{sec:app} for explanations).
With some exceptions (e.g., \cite{chizat2018global, nguyen2020rigorous,
nitanda2021particle, ren2022entropic, CCRW23}),
the majority of the studies \cite{mei2019mean,HRSS19,chizat2022mean,nitanda2022convex} have focused on the entropy-regularized mean field
optimization problem
and the corresponding MFL dynamics in the form of \eqref{eq:mfl-intro}.
It was first proved in \cite{HRSS19} that
under the convexity assumption of $F$,
the marginal distributions of the MFL dynamics converge toward its unique
invariant measure,
which is also the unique minimizer of the mean field optimization problem.
Then it is shown in \cite{nitanda2022convex,chizat2022mean} that,
with the presence of the uniform LSI,
such kind of convergence is exponentially fast.
The main contribution of this paper lies in that, we further explore the fine properties of MFL dynamics
with a particular emphasis on its uniform-in-time propagation of chaos property,
i.e., the time-uniform upper bounds for the distance between the finite-particle
and the mean field dynamics.
Therefore, we provide a theoretical guarantee
for the applicability of the finite-particle approximation
when the dynamics is expected to run for an indefinitely long time.

Recall that we have defined \(m_t = \Law(X_t)\).
Let us also define \(m^N_t = \Law(X^1_t,\ldots,X^N_t)\)
and denote by \(m_\infty\) the unique invariant measure of the mean field
dynamics.
Our main results are summarized as follows:
\begin{itemize}
\item if the Radon--Nikod\'ym derivative
\(\dd m_0 / \! \dd m_\infty\) belongs to \(L^{p_0}(m_\infty)\)
for some \(p_0 > 1\),
then for every \(p \in \mathbb R\),
the norm \(\Vert \dd m_t /\! \dd m_\infty \Vert_{L^p(m_\infty)} \to 1\)
exponentially fast when \(t \to \infty\);
\item the scaled \(L^2\)-Wasserstein distance and the relative entropy
\(\frac 1N W_2^2(m^N_t, m_\infty^{\otimes N})\),
\(\frac 1N H(m^N_t | m_t^{\otimes N})\)
converge to a  \(O(N^{-1})\) neighborhood of zero when \(t \to
\infty\), with an exponential rate that is independent of $N$;
\item if the initial error is zero, i.e., $m^N_0 = m_0^{\otimes N}$,
then \(\sup_{t \in [0,\infty)}\frac 1N W_2^2(m^N_t, m^{\otimes N}_t) \to 0\)
when $N \to \infty$;
further if the assumption of the first claim holds,
then $\sup_{t \in [0,\infty)}\frac 1N H(m^N_t | m^{\otimes N}_t) \to 0$
when $N \to \infty$.
\end{itemize}

We also refer those interested readers to our companion paper \cite{ukploc}, which delves into analogous properties for kinetic MFL dynamics.

\subsection{Related works}
\label{sec:related-works}

\paragraph{Long-time behavior of McKean--Vlasov dynamics.}
Propagation of chaos in finite time for the stochastic McKean--Vlasov dynamics
\[
\dd X_t = b( m_t, X_t) \dd t + \sqrt 2\dd W_t,
\qquad \text{where $m_t = \Law(X_t)$}
\]
is relatively easy to show, using the \emph{synchronous coupling} approach,
given that \(b\) is a jointly Lipschitz function of both measure and space
variables
in the sense of the Wasserstein metric.
The bound obtained by this method, however, generally tends to infinity when
the time interval extends to infinity.
Besides, the dynamics may possess multiple invariant measures,
so uniform-in-time convergence can not be expected without some
additional assumptions or a more general definition of convergence itself
(e.g.\ convergence modulo symmetries).

The research on the long-time behavior of McKean--Vlasov dynamics has been
active in
recent years
and here we introduce a setting that has appeared in many previous works.
Consider functions $U$, $V : \mathbb R^d \to \mathbb R$ and the following
special kind of drift
\[
b(m,x) = - \nabla U(x) - \int \nabla V(x - \tilde x) m(\dd\tilde x).
\]
In this case, \(U\) is referred as the \emph{external potential} and \(V\)
is called
the \emph{interaction potential}.

In this paragraph, we provide a far from exhaustive review
of uniform-in-time propagation of chaos (POC) for McKean--Vlasov dynamics.
First, in the work \cite{malrieu2001logarithmic} of Malrieu in 2001, uniform
POC is established by synchronous coupling
for overdamped dynamics under the assumption that
$U$ is strongly convex and $V$ is convex.
In an alternative way, Carrillo, McCann and Villani set up the mean field
gradient flow framework in their work \cite{carrillo2003kinetic}, which our
paper also relies on. They showed the exponential convergence of the overdamped
mean field system under the assumption that \(U + 2V\) is strongly convex.
In Monmarché's work \cite{monmarche2017long},
uniform POC is extended to the \emph{kinetic} Langevin dynamics,
assuming the same convexity assumption on \(U + 2V\).
This assumption is further relaxed in his follow-up work with Guillin
\cite{guillin2021uniform},
where they incorporate the uniform-in-$N$ log-Sobolev inequality in
\cite{guillin2022uniform}.
In \cite{durmus2020elementary}, Durmus, Eberle, Guillin and Zimmer showed
uniform POC for overdamped Langevin dynamics, under the assumption that the
confining potential \(U\) is only weakly convex and \(V\) is small enough,
utilizing a \emph{reflection coupling} technique.
The reflection coupling technique is then used by Schuh in
\cite{schuh2022global} to show uniform POC for kinetic Langevin dynamics,
albeit in this setting, the form of the confining potential is more restricted
compared to the overdamped case.
The weak uniform-in-time convergence is also demonstrated for the overdamped
dynamics on a
torus in \cite{delarue2021uniform} by Delarue and Tse under various settings.
This research assumes the smallness of interaction without explicitly
specifying its form and employs a \emph{master equation} analysis.
In \cite{lacker2023sharp}, Lacker and Le Flem showed a sharp $O(1/N^2)$ rate
for time-uniform propagation of chaos for the overdamped dynamics, by
studying the relative entropy growth between marginal distributions
with the help of a time-uniform log-Sobolev inequality for the mean field flow.

We now comment on the assumptions and methods of these works.
Apart from the second and third settings of \cite{delarue2021uniform} and that
of \cite{lacker2023sharp},
the aforementioned works all rely on the smallness or the (semi-, weak)
convexity of the interaction potential.
This smallness or convexity is used to control the error between the
coupled processes,
or to deduce a uniform-in-\(N\) log-Sobolev inequality for the \(N\)-particle
system's invariant measure (see \cite{guillin2022uniform}).
Our setting is different from those in other works.
First, our results are built upon the functional convexity
of the mean field energy functional,
which is a different (and even exclusive in some cases) assumption
from the convexity of the interaction potential.
Further details on this alternative assumption of convexity will be provided in
the following paragraph.
Second, our approach does not rely on a uniform-in-\(N\) log-Sobolev inequality
for the invariant measure of the \(N\)-particle system.

Finally, we remark that the translation-invariant models have been studied in
the last setting of \cite{delarue2021uniform} and also in
\cite{durmus2022sticky}. In these cases, there exists a continuum of invariant
measures, and the POC is then obtained modulo the translational symmetry.
Besides, we also mention that in a recent work \cite{guillin2021uniformvortex},
Guillin et al.\ studied the 2D viscous vortex model
where the particles are in singular interactions
and showed the uniform POC estimates.

\paragraph{Linear functional convexity.}

One of our key assumptions is the (linear functional) convexity of the mean
field functional \(F\), formally defined in \cref{eq:convex}.
Except in \cite{suzukiuniform,delarue2021uniform}, this assumption has not been
explicitly exploited to investigate the long-time behavior of the
McKean--Vlasov dynamics.
It is important to distinguish this convexity from the displacement convexity,
which frequently appears in the optimal transport literature and is
defined in, for instance, \cite[Definition 16.1]{villani2009optimal}.
We will clarify in \cref{rem:exclusion-two-convexities} that, for continuous
two-body interaction potentials, Bochner's theorem implies that these two
concepts are even mutually exclusive, except in trivial cases.

This particular form of convexity is implicitly exploited in
\cite{delarue2021uniform} to obtain time-uniform POC estimates.
More precisely, the authors studied McKean--Vlasov drift of form
\(b(m,x) = - \int \nabla V(x - \tilde x) m(\dd\tilde x)\)
on the torus, where all Fourier coefficients
of the interaction potential \(V\) are nonnegative.
Then this property is used to obtain estimates on the master equation
in the long time.
We note that, here, the positivity of the Fourier coefficient
implies that the corresponding energy
\(F(m) = \frac 12\iint V(x - \tilde x) m(\dd x) m(\dd\tilde x)\)
is convex in our functional sense.
Although our results are stated for dynamics in \(\mathbb R^d\), it is
reasonable to expect that our methodology can be extended to the torus and
yield similar results.

The primary motivation for introducing this new setting is to study the
training of two-layer (or one-hidden-layer) neural networks,
which we will explain in \cref{exm:nn,exm:nn-ngd}.

\paragraph{Gradient descent.}

Our dynamics is a special case of McKean--Vlasov with gradient-type drift:
\[
b(m,x) = - D_m F(m,x) = - \nabla \frac{\delta F}{\delta m}(m,x).
\]
This form of drift corresponds to the gradient descent of the \emph{free
energy} \(\mathcal F = F + H\) in \(L^2\)-Wasserstein space,
here, \(H(m) = \int m(x) \log m(x) \dd x\) is the (absolute) entropy
of the measure.
We refer the readers to \cite{jko}
for detailed discussions about the gradient flow
with the linear energy \(F(m) = \int V(x) m(\dd x)\),
and \cite{gf} for a general gradient flow framework in Wasserstein space.
We note that, in a previous work \cite{HRSS19},
this gradient flow structure is exploited to obtain
the ergodicity of the MFL dynamics.
Precisely, the authors established
the following free energy dissipation formula
\[
\frac{\dd\mathcal F(m_t)}{\dd t}
= - \int |D_m F(m_t,x) + \nabla \log m_t(x)|^2
m_t(\dd x),
\]
and then by combining this with LaSalle's invariance principle
and the uniqueness of the invariant measure,
they showed the global convergence of the MFL dynamics.
In this paper, we will prove the same energy descent formula
under weaker assumptions on the regularity of \(x \mapsto D_m F(m,x)\),
thanks to the general framework developed in \cite{gf}.

\subsection{Main contributions}

\paragraph{\texorpdfstring{\(L^p\)}{Lp} convergence and hypercontractivity of
MFL.}

The exponential convergence of relative entropy for the MFL with convex \(F\)
has been proved in \cite{chizat2022mean,nitanda2022convex}
via log-Sobolev inequalities,
extending the classical result \cite{ov} wherein the \(F\) is linear in measure.
In this paper, we introduce a stronger \(L^p\)-convergence in
\cref{thm:lp-convergence}.
To achieve this enhanced convergence result, we require the initial condition
to lie in \(L^{p_0}\) for some \(p_0 > 1\).
This contrasts with the situation of relative entropy,
where elliptic regularization ensures relative entropy to be finite at all
positive times (see \cref{prop:mf-exist-unique-regular}).

Our method of proof is based on the \(L^2\)-convergence and the
hypercontractivity,
which ports the \(L^2\)-convergence to \(L^p\) for all \(p \in \mathbb R\).
Two pivotal observations are the growth of \(L^p\)-norm formula
\cref{eq:lp-growth}
and the hypercontractive inequalities
\cref{eq:hypercontractivity-q>1}-\cref{eq:hypercontractivity-q<1}
for the mean field flow.
Recently the hypercontractivity has also been ultilized in
\cite{decourcel2023sharp}
to show the \(L^p\)-convergence of MFL with Riesz interactions (though on a
torus).
Finally, it is important to mention that the proof of our propagation of
chaos result (\cref{thm:poc}) requires the \(L^p\)-convergence for \(p\) negative.
To address this requirement, we establish the \emph{reverse} hypercontractivity
of the MFL.
This property follows from the analogous formal computations to those employed
in direct hypercontractivity, under the assumption that the invariant measure
satisfies a LSI.

\paragraph{Convergence of particle system.}

Within the mean field setting established in
\cite{chizat2022mean,nitanda2022convex},
we show in \cref{thm:ps-entropy-convergence} that
the particle system's free energy converges to the \(N\)-tensorized invariant
measure of the mean field system exponentially modulo an error of size
\(O(N^{-1})\) per particle.
Our proof approach relies on a decomposition of relative Fisher information
and a componentwise application of the log-Sobolev inequality,
which introduces the \(O(N^{-1})\) error per particle.
Our result differs from that of \cite{guillin2022uniform},
where the precise convergence of the particle system to its invariant measure
is obtained through the use of the uniform-in-\(N\) log-Sobolev inequality.
One notable advantage of our method is that we allow applications involving
potentially significant interactions, including cases such as the training of
neural networks
(as discussed in \cref{exm:nn,exm:nn-ngd}.)

\paragraph{Propagation of chaos.}

By combining the two previous results,
i.e.\ the \(L^p\)-convergence of the MFL
and the entropic convergence of the particle system,
we are able to control the distance between the particle system \(m^N_t\) and
\(N\)-tensorized mean field flow \(m^{\otimes N}_t\),
in terms of Wasserstein distance \emph{and} relative entropy.
The bound on Wasserstein is a direct consequence of Talagrand's \(T_2\)
transport inequality.
To control the relative entropy
we employ a classical duality formula \cref{eq:duality-entropy}
to link \(H(m^N_t | m^{\otimes N}_t)\) to the \(-p\) norm
\(\Vert\dd m_t / \!\dd m_\infty \Vert_{-p}\) for \(p > 0\),
whose exponential convergence is guaranteed by \cref{thm:lp-convergence}.
As a side result, we also obtain the uniform-in-time concentration of measure
of the mean field flow (\cref{thm:mf-uniform-concentration}), based on this
observation.

Let us now compare our method to those of \cite{lacker2023sharp,suzukiuniform}.
In \cite{lacker2023sharp} the authors assumed the mean field flow satisfies a
uniform LSI
and utilized an entropy growth formula similar to our \(L^p\)-growth formula to
estimate the relative entropy bound.
As remarked in \cite{suzukiuniform},
verifying this uniform LSI can be challenging in the mean field setting.
In particular if one wishes to apply the Holley--Stroock perturbation lemma to
the invariant measure \(m_\infty\),
the mean field flow needs
to satisfy \(\log \dd m_t /\!\dd m_\infty \in L^\infty\) uniformly
In \cite{suzukiuniform}, Suzuki, Nitanda and Wu made the assumptions that the
confining potential exhibits a super-quadratic growth, so that this boundedness
follows from the ultracontractivity via super LSI.
However, this confining potential is stronger than the quadratic one in our
setting
and the constants derived from ultracontractivity are dependent on the spatial
dimension.
% An alternative approach, which was discussed in an earlier version of the
% current paper,
% requires the initial condition to satisfy \(\log m_0 / m_\infty \in L^\infty\).
% Using our current method, we remain within the quadratic growth regime, and the
% only additional requirement, as explained earlier,
% is that the initial condition satisfies \(m_0 / m_\infty \in L^{p_0}\) for some
% \(p_0 > 1\).

\subsection{Notations}
\label{sec:notations}
Let \(d\) be a positive integer and \(x\) an element of \(\mathbb R^d\).
We denote the Euclidean norm of \(x \in \mathbb R^d\) by \(|x|\)
and define \(c_d\) as the volume of the \(d\)-dimensional unit ball.
Let \(p \geq 1\),
we define \(\mathcal P_p(\mathbb R^d)\) to be the space of probability measures
on \(\mathbb R^d\) with finite \(p\)-moment, i.e.,
\(\mathcal P_p (\mathbb R^d)
= \{ m \in \mathcal P (\mathbb R^d) : \int |x|^p m (\dd x) < +\infty\}\).
The \(L^p\)-Wasserstein metric is denoted by \(W_p\) and its definition along
with elementary properties, can be found in
\cite[Chapter 7]{gf}.

Consider a mean field functional \(F : \mathcal P_2(\mathbb R^d) \to \mathbb
R\).
We denote by \(\frac{\delta F}{\delta m} : \mathcal P_2(\mathbb R^d) \times
\mathbb R^d \to \mathbb R\) its linear functional derivative
and by \(D_m F = \nabla \frac{\delta F}{\delta m} : \mathcal P_2(\mathbb R^d)
\times \mathbb R^d \to \mathbb R^d\) its intrinsic derivative, provided they
exist.
The definition of linear functional derivative on \(\mathcal P_2(\mathbb R^d)\)
can be found in \cite[Definition 5.43]{probamfg1}.

Let $X$, $Y$ be two random variables.
We denote the distribution of \(X\) as \(\Law(X)\)
and write \(X \sim m\) when \(m = \Law(X)\). Additionally,
we use \(X \overset{d}= Y\) to indicate that \(\Law(X) = \Law(Y)\). The set of
couplings between probability measures $\mu$, $\nu$
is denoted by \(\Pi(\mu, \nu)\). Let \(N \geq 2\) be an integer, we use the bold
letter \(\mathbf x_N = (x^1,\ldots,x^N)\) to represent an \(N\)-tuple of the
elements in \(\mathbb R^{d} \).
We omit the subscript \(N\) when there are no ambiguities.

Let \(I \subset \{1,\ldots,N\}\).
We define \(-I \coloneqq \{1,\ldots,N\} \setminus I\),
i.e., the complementary index set of $I$.
For a probability measure
\(m^N = \Law(\mathbf X) \in \mathcal P(\mathbb R^{dN})\),
we denote its marginal and the (regular) conditional distributions by
\begin{align*}
m^{N,I} &= \Law (X^i)_{i \in I}, \\
m^{N,I|-I} (\mathbf x^{-I}) &= \Law \bigl( (X^i)_{i \in I}
\big| X^j = x^j,~j \in -I \bigr),
\end{align*}
where the latter is defined \(m^{N,-I}\)-almost surely
and \(\mathbf x^{-I}\) denotes the tuple \((x^j)_{j \in -I}\).
We identify \(i\) with the singleton \(\{i\}\) when working with indices.

Given \(\mathbf x_N = (x^1,\ldots,x^N) \in \mathbb R^{dN}\), we denote the
corresponding empirical measure by
\[
\mu_{\mathbf x_N} = \frac 1N \sum_{i=1}^N \delta_{x^i}.
\]
For \(i = 1,\ldots, N\), as introduced in the paragraph above,
the symbol $-i$ denotes the complementary set $\{1,\ldots,N\} \setminus i$.
We denote the empirical measure of the $N-1$ points
\(\mathbf x^{-i}_N=(x_j)_{j \neq i}\) by
\[
\mu_{\mathbf x^{-i}_N} = \frac 1{N-1} \sum_{j=1,\,j\neq i}^N \delta_{x^j}.
\]
For a \(\mathbb R^{dN}\)-valued random variable \(\mathbf X_N = (X^i)_{i=1}^N\),
we can thereby form the random empirical measures $\mu_{\mathbf X_N}$,
$\mu_{\mathbf X_N^{-i}}$.

When a measure \(m \in \mathcal P(\mathbb R^d)\) has a density with respect
to the \(d\)-dimensional Lebesgue measure,
we still denote its density function by \(m : \mathbb R^d \to \mathbb R\).
Let \(\gamma\) be a positive and \(\sigma\)-finite measure on \(\mathbb R^d\).
We define the relative entropy
\begin{align*}
H(m | \gamma ) &= \int \log \frac{\dd m}{\dd\gamma} (x) m(\dd x) \\
\intertext{and the relative Fisher information}
I(m | \gamma) &= \int \biggl| \nabla \log \frac{\dd m}{\dd\gamma} \biggr|^2
m(\dd x)
\end{align*}
provided the corresponding integrals are well defined.
In cases where the integrals are not well defined, we set $H$, $I = +\infty$
respectively.
When \(\gamma = \mathcal L^d\) is the Lebesgue measure on \(\mathbb R^d\), we
omit the dependence on \(\gamma\) and define the \emph{absolute} entropy and
Fisher information as:
\[
H(m) \coloneqq H(m | \mathcal L^d ), \qquad I(m) \coloneqq I(m | \mathcal L^d ),
\]
provided they are well-defined.
For nonnegative functions \(f : \mathbb R^d \to [0,+\infty)\) we also define
its entropy as
\[
\Ent_m f = \Expect_m [ f \log f ] - \Expect_m [ f ] \log \Expect_m [ f ],
\]
which is well defined in \([0,+\infty]\) according to Jensen's inequality.

\paragraph{Organization of paper.}

In \cref{sec:main-results}, we present our assumptions,
introduce the mean field Langevin dynamics and the particle system,
and state our main results.
In section \ref{sec:app}, we offer some examples of MFL, to which our theorems can
be applied, accompanied by numerical experiments of two-layer
neural network training.
The proofs are given in the rest of the paper, and for the most technically
demanding ones, we detailed in Section \ref{sec:mf-technical}.

\section{Main results}
\label{sec:main-results}

\paragraph{Assumptions.}
Let \(F : \mathcal P_2(\mathbb R^d) \to \mathbb R\) be a mean field functional.
We suppose \(F\) is convex in the sense that
for all \(t \in [0,1]\) and all $m$, $m'\in \mathcal P_2(\mathbb R^d)$,
\begin{equation}
\label{eq:convex}
F\bigl((1-t)m + tm'\bigr) \leq (1-t)F(m) + tF(m').
\end{equation}
Suppose also its intrinsic derivative \(D_m F : \mathcal P_2(\mathbb R^d)
\times \mathbb R^d \to \mathbb R^d\)
exists and satisfies
\begin{equation}
\label{eq:lip-in-m}
\forall x \in \mathbb R^d,~\forall m,m' \in \mathcal P_2(\mathbb R^d),~
|D_m F(m,x) - D_m F(m',x)| \leq M^F_{mm} W_1(m,m')
\end{equation}
for some constant \(M^F_{mm} \geq 0\).
For each \(m\in \mathcal{P}_2(\mathbb R^d)\), we define a probability measure
\(\hat m\) by its density
\[\hat m(x) \propto \exp \biggl( - \frac{\delta F}{\delta m}(m,x) \biggr)\]
and suppose \(\hat m\)
satisfies the \(\rho\)-\emph{logarithmic Sobolev inequality} (LSI) uniformly in
\(m\) for some \(\rho > 0\), that is,
for every \(m \in \mathcal P_2(\mathbb R^d)\),
\begin{equation}
\label{eq:lsi}
\forall f \in C^1_b(\mathbb R^d),\qquad
\rho \Ent_{\hat m} (f^2) \leq \Expect_{\hat m} [|\nabla f|^2].
\end{equation}
Here, we implicitly suppose that $\hat m$ is well defined for all
$m \in \mathcal P_2(\mathbb R^d)$, and in particular, we have
$\int \exp \bigl( - \frac{\delta F}{\delta m}(m,x)\bigr)\dd x < \infty$.
We remark that the inequality above can be verified
for mean field functionals $F$ whose linear derivative
$\frac{\delta F}{\delta m}$ is a perturbation of a strongly convex function.
For details, we refer readers to \cref{prop:lsi} in \cref{sec:exm}.
We suppose as well
\begin{equation}
\label{eq:first-der}
\sup_{m \in \mathcal P_2(\mathbb R^d)} \sup_{x \in \mathbb R^d} |\nabla D_m
F(m,x)|
\leq M^F_{mx}
\end{equation}
for some constant $M^F_{mx} \geq 0$.
Finally, for some of the results we additionally suppose that \(x \mapsto D_m
F(m,x)\) belongs to \(C^3\) with the bounds
\begin{equation}
\label{eq:higher-der}
\sup_{m \in \mathcal P_2(\mathbb R^d)} \sup_{x \in \mathbb R^d} |\nabla^k D_m
F(m,x)| < +\infty,\qquad
k = 2,\,3.
\end{equation}

\begin{rem}[Well-definedness of \(\hat m\)]
The definition of $\hat m$ relies on the finiteness
of the normalization constant
\begin{equation}
Z(\hat m) = \int \exp \biggl( - \frac{\delta F}{\delta m}(m,x) \biggr)\dd x.
\label{eq:def-Z-m-hat}
\end{equation}
As mentioned above,
it is assumed implicitly in the condition \cref{eq:lsi} that
$Z (\hat m)$ is finite for every $m \in \mathcal P_2(\mathbb R^d)$.
We will prove in \cref{prop:exist-m-hat} that
the following is sufficient for this finiteness:
\begin{itemize}
\item the condition \cref{eq:lip-in-m} holds, and
\item there exists at least one measure \(m_0\)
such that \(Z(\hat m_0)\) is finite and \(m_0\) satisfies the LSI \cref{eq:lsi}.
\end{itemize}
\end{rem}

\begin{rem}[Functional inequalities]
By approximating the function $f$ by a sequence of functions in
$C^1_\textnormal{b}$, we find that the inequality \cref{eq:lsi}
holds for functions \(f\) whose generalized derivative satisfies
\(\Expect_{\hat m} [|\nabla f|^2] < +\infty\).
It is well known that the LSI \cref{eq:lsi} implies the \emph{Poincaré
inequality}:
\begin{equation}
\label{eq:poincare}
\forall f \in C^1_b(\mathbb R^d),\qquad
2\rho \Var_{\hat m} (f) \leq \Expect_{\hat m} [|\nabla f|^2].
\end{equation}
The restriction \(f \in C^1_\textnormal{b}\) can be analogously removed.
The LSI \cref{eq:lsi} also implies
\emph{Talagrand's \(T_2\)-transport inequality}:
\begin{equation}
\forall \mu \in \mathcal P_2(\mathbb R^d),\qquad
\rho W_2^2(\mu, \hat m) \leq H(\mu | \hat m).
\label{eq:t2}
\end{equation}
See the original work of Otto and Villani \cite{ov} for a proof.
We also sketch their argument in the proof of \cref{lem:lsi-implies-moments}.
All those three inequalities, namely \cref{eq:lsi,eq:poincare,eq:t2},
are stable under tensorization:
if one replaces, for some \(N \geq 2\), the measure \(\hat m\)
by its tensor product \(\hat m^{\otimes N}\),
which is a measure on \(\mathbb R^{dN}\),
and the function \(f : \mathbb R^d \to \mathbb R\)
(resp.\ the probability measure $\mu$ on $\mathbb R^d$)
by function $f^N : \mathbb R^{dN} \to \mathbb R$
having a square-integrable weak derivative $\nabla f^N$
with respect to the measure $\hat m^{\otimes N}$
(resp.\ probability measures $\mu^N$ on $\mathbb R^{dN}$),
then the inequalities hold with the same constant \(\rho\).
\end{rem}

\paragraph{Mean field and particle system.}

We study the \emph{mean field Langevin dynamics}, that is, the following
McKean--Vlasov SDE
\begin{equation}
\dd X_t = - D_m F(m_t, X_t) \dd t + \sqrt 2 \dd W_t,
\qquad \text{where $\Law (X_t) = m_t$.}
\label{eq:mf-sde}
\end{equation}
Let \(N \geq 2\).
The corresponding \(N\)-\emph{particle system} is defined by
\begin{equation}
\label{eq:ps-sde}
\dd X^i_t = - D_m F(\mu_{\mathbf X_t}, X^i_t) \dd t + \sqrt 2 \dd W^i_t,\quad
i=1,\ldots,N,\quad\text{where $\mu_{\mathbf X_t} = \frac 1N \sum_{i=1}^N
\delta_{X^i_t}$.}
\end{equation}
Here, $W$, $W^i$ are standard Brownian motions
in \(\mathbb R^d\), which are independent from each other.
Their marginal distributions \(m_t = \Law(X_t)\),
\(m^N_t = \Law(\mathbf X_t) = \Law(X^1_t,\ldots,X^N_t)\)
then solve the Fokker--Planck equations respectively
\begin{align}
\partial_t m_t &= \Delta m_t + \nabla
\cdot \bigl(D_m F(m_t, \cdot) m_t\bigr),
\label{eq:mf-fp}\\
\partial_t m^N_t &= \sum_{i=1}^N \Bigl( \Delta_i m^N_t
+ \nabla_i \cdot \bigl(D_m F(\mu_{\mathbf x}, x^i) m^N_t\bigr) \Bigr).
\label{eq:ps-fp}
\end{align}
The mean field equation \cref{eq:mf-fp} is non-linear while the \(N\)-particle
system equation \cref{eq:ps-fp} is linear.
We will prove in \cref{prop:mf-exist-unique-regular}
that, if the initial condition $m_0 \in \mathcal P_2(\mathbb R^d)$, the mean field dynamics \cref{eq:mf-fp}
is well posed and enjoys certain regularity.

\begin{rem}
\label{rem:vol-scaling}
We have fixed the volatility (diffusion) constant to be \(\sqrt{2}\) to
simplify our computations.
In order to apply our results to the MFL defined by
\[
\dd X_t = - D_m F(m_t, X_t) \dd t + \sigma \dd W_t,
\qquad \text{where $\Law (X_t) = m_t$,}
\]
with some \(\sigma > 0\),
we apply the rescaling: \(\tilde t = \frac{\sigma^2}{2} t\), \(\tilde F =
\frac{2}{\sigma^2} F\) and \(\tilde X_{\tilde t} = X_t\).
In this way, the new diffusion process \(\tilde t \mapsto \tilde X_{\tilde t}\)
satisfies the SDE \cref{eq:mf-sde}, whose diffusion constant is fixed to
$\sqrt 2$, with the new mean field functional \(\tilde F\).
The same scaling transform can be applied to the particle system as well.
\end{rem}

\paragraph{Free energy and invariant measure.}
We focus on the long-term behavior of the MFL \cref{eq:mf-fp}
and the corresponding particle system \cref{eq:ps-fp},
where invariant measures play a key role.
Define \emph{mean field free energy} \(\mathcal F : \mathcal P_2(\mathbb
R^d) \to (-\infty,+\infty]\) by
\begin{equation}
\label{eq:def-mf-free-energy}
\mathcal F(m) = F(m) + H(m).
\end{equation}
Given the assumptions \cref{eq:convex,eq:lip-in-m,eq:lsi,eq:first-der},
we can show the existence of a unique minimizer of \(\mathcal F\),
denoted by \(m_\infty\).
Furthermore, this measure \(m_\infty\) satisfies the \emph{first-order
condition}:
\begin{equation}
\label{eq:mf-foc}
m_\infty(\dd x) = \hat m_\infty(\dd x)
= \frac 1{Z(\hat m_\infty)}
\exp \biggl( - \frac{\delta F}{\delta m}(m_\infty,x) \biggr) \dd x.
\end{equation}
The precise statement and proof is given in \cref{prop:mf-invariant-measure}.
Differentiating both sides of the first-order condition,
we obtain \(\Delta m_\infty + \nabla
\cdot \bigl(D_m F(m_\infty,x) m_\infty\bigr) = 0\),
which implies that \(m_\infty\) is an invariant measure to mean field Fokker--Planck
equation \cref{eq:mf-fp}.
Conversely,
we will show in \cref{cor:mf-invariant-measure-unique} that under our
conditions every invariant measure satisfies the first-order condition
and, therefore, we get the uniqueness of invariant measure as well.

The \(N\)-particle system \cref{eq:ps-sde} is a classical Langevin dynamics because
the equation \cref{eq:ps-fp} is linear.
We define the \emph{\(N\)-particle free energy}
\(\mathcal F^N : \mathcal P_2(\mathbb R^{dN}) \to (-\infty,+\infty]\) by
\begin{equation}
\label{eq:def-ps-free-energy}
\mathcal F^N(m^N) = N \int F(\mu_{\mathbf x}) m^N(\dd\mathbf x) + H(m^N).
\end{equation}
We will prove in \cref{prop:ps-invariant-measure} that
under our assumptions \cref{eq:convex,eq:lip-in-m,eq:lsi}
the minimizer \(m^N_\infty\) of \(\mathcal F^N\) exists, and has the density
\begin{equation}
\label{eq:def-ps-invariant-measure}
m^N_\infty(\dd\mathbf x)
\propto \exp \bigl( -NF (\mu_\mathbf{x}) \bigr) \dd\mathbf x,
\end{equation}
which is invariant to the \(N\)-particle Fokker--Planck equation \cref{eq:ps-fp}.
By the definition of free energy we have \(\mathcal F^N(m^N) = H(m^N |
m^N_\infty) + \textnormal{constant}\),
so \(m^N_\infty\) also minimizes the \(N\)-particle free energy \(\mathcal
F^N\).

\paragraph{\(L^p_+\) space for all \(p \in \mathbb R\).}
We investigate the convergence of the marginal distributions of the mean field
dynamics in the \(L^p(m_\infty)\)-norm for all \(p \in \mathbb R\)
and take special care when \(p < 1\).
Let \(\mu\) be a probability measure on \(\mathbb R^d\)
and \(h : \mathbb R^d \to [0,+\infty]\) be a measurable function.
For \(p \neq 0\) define
\begin{align*}
\Vert h \Vert_{L^p(\mu)} &= \biggl( \int h(x)^p \mu(\dd x) \biggr)^{\!1/p}, \\
\intertext{and for \(p = 0\) define}
\Vert h \Vert_{L^0(\mu)} &= \exp \biggl( \int \log h(x) \mu(\dd x) \biggr). \\
\intertext{We say \(h \in L^p_+(\mu)\) if}
\Vert h \Vert_{L^p(\mu)} &\begin{cases}
< +\infty& \text{if}~p > 0, \\
\in (0,+\infty)& \text{if}~p = 0, \\
> 0& \text{if}~p < 0.
\end{cases}
\end{align*}
It is well-known that \(p \mapsto \Vert h \Vert_p\) is increasing,
which ensures that the \(0\)-norm is well defined once \(\Vert h \Vert_p <
+\infty\) for some \(p > 0\) or \(\Vert h \Vert_q > 0\) for some \(q < 0\).
In this paper we will only work with \(\mu\) equal to \(m_\infty\),
the mean field invariant measure.
In this case we write \(\Vert h \Vert_p = \Vert h \Vert_{L^p(m_\infty)}\)
for simplicity.
We also say \(h \in L^{1+} (m_\infty)\) or \(h\) is \(L^{1+}\)-integrable
if there exists a number \(p_0 > 1\) such that \(h \in L^{p_0} (m_\infty)\).

\paragraph{Statement of main results.}

Recall that \(m_t\) and \(m^N_t\) are the respective marginal distributions of
the mean field and the \(N\)-particle system \cref{eq:mf-sde,eq:ps-sde}.
We slightly improve the exponential energy dissipation result
for the MFL dynamics \cref{eq:mf-sde}.

\begin{thm}[Energy dissipation of MFL]
\label{thm:mf-entropy-convergence}
Assume \(F\) satisfies \cref{eq:convex,eq:lip-in-m,eq:lsi,eq:first-der}.
If \(m_{t_0}\) has finite entropy and finite second moment
for some \(t_0 \geq 0\),
then for every \(t \geq t_0\),
\begin{equation}
\label{eq:mf-entropy-convergence}
H(m_t | m_\infty) \leq \mathcal F(m_t) - \mathcal F(m_\infty)
\leq \bigl(\mathcal F(m_{t_0}) - \mathcal F(m_\infty)\bigr) e^{-4\rho(t-t_0)}.
\end{equation}
\end{thm}

\begin{rem}
The theorem stated here differs slightly from the previous results
(\cite[Theorem 3.2]{chizat2022mean} and \cite[Theorem 1]{nitanda2022convex}),
in that we have removed the technical restriction
that \(x \mapsto D_m F(m, x)\) is infinitely differentiable.
This is achieved by using
the differential calculus in the Wasserstein space developed
in the monograph \cite{gf}.
\end{rem}

The proof of the theorem is postponed to
\cref{sec:proof-thm:mf-entropy-convergence}.

We also study the MFL system's convergence beyond the entropic sense.
In particular, we show
that the system converges in the $L^2$ sense given $L^2$-initial values
(\cref{prop:l2-convergence}),
and that the system is also hypercontractive
and reverse-hypercontractive (\cref{prop:hypercontractivity}).

Denote
\[
h_t(x) \coloneqq \frac{\dd m_t}{\dd m_\infty} (x)
\]
for the solution $m_t$ of the MFL dynamics \cref{eq:mf-fp},
where $m_\infty$ is the unique invariant measure to the MFL,
satisfying \cref{eq:mf-foc}.

\begin{prop}[\(L^2\)-convergence]
\label{prop:l2-convergence}
Assume \(F\) satisfies
\cref{eq:convex,eq:lip-in-m,eq:lsi,eq:first-der,eq:higher-der}.
Let \(m_t \in C\bigl([0,+\infty); (\mathcal P_2, W_2)\bigr)\)
be a solution to \cref{eq:mf-fp}.
If \(h_{t_0} \in L^2(m_\infty)\),
then \(h_t \in L^2(m_\infty)\) for all \(t \geq t_0\).
Moreover, for all \(\rho' \in (0,\rho)\), we have
\begin{equation}
\label{eq:l2-convergence}
\forall t \geq t_0,\qquad
\Vert h_{t} - 1 \Vert_2^2 \leq M e^{-4\rho' (t - t_0)},
\end{equation}
for the constant $M$ defined by
\[
M = \exp \biggl( \frac{\Delta(t_0)}{4\rho} \biggr)
\biggl( \lVert h_{t_0} - 1 \rVert_2^2 + \frac{\Delta(t_0)}{4(\rho - \rho')}
\biggr),
\]
where
\[
\Delta(t_0) =
\frac{(M^F_{mm})^2}{\rho - \rho'}
\biggl( 1 + \frac{M^F_{mm}}{\rho} + \frac{(M^F_{mm})^2}{2\rho^2} \biggr)
\log \lVert h_{t_0} \rVert_2.
\]
\end{prop}

\begin{prop}[Hypercontractivity]
\label{prop:hypercontractivity}
Assume \(F\) satisfies
\cref{eq:convex,eq:lip-in-m,eq:lsi,eq:first-der,eq:higher-der}.
Suppose \(h_{t_0} \in L^{q_0}(m_\infty)\) for some \(q_0 \neq 1\).
Let \(\varepsilon \in (0,1]\) and \(q(t)\) solve the ODE \(\dot q = 4 (1 -
\varepsilon) \rho (q - 1)\) with the initial condition \(q(t_0) = q_0\).
Then \(h_t \in L^{q(t)} (m_\infty)\) for \(t \geq t_0\).
Moreover, we have for \(q_0 > 1\),
\begin{equation}
\label{eq:hypercontractivity-q>1}
\log \lVert h_t \rVert_{q(t)} \leq \log \lVert h_{t_0} \rVert_{q_0}
+ \int_{t_0}^t \delta(s) \dd s,
\end{equation}
and for \(q_0 < 1\),
\begin{equation}
\label{eq:hypercontractivity-q<1}
\log \lVert h_t \rVert_{q(t)} \geq \log \lVert h_{t_0} \rVert_{q_0}
+ \int_{t_0}^t \delta(s) \dd s,
\end{equation}
where \(\delta(t) = \frac{1}{4\varepsilon} \bigl(q(t) - 1\bigr)
(M^F_{mm})^2 W_1^2(m_t, m_\infty)\).
\end{prop}

\begin{rem}[Optimality of exponent's growth]
\label{rem:hypercontractivity-exp-growth}
In the case where the mean field interaction is absent,
Nelson's theorem \cite[Théorème 2.3.1]{lsi}
shows the optimality of the exponent's growth in \cref{prop:hypercontractivity}.
\end{rem}

The proofs of
\cref{prop:l2-convergence,prop:hypercontractivity}
are given in \cref{sec:l2-convergence-hypercontractivity}.

By combining the $L^2$-convergence and the hypercontractivity,
we can obtain the $L^p$-convergence of the MFL dynamics.

\begin{thm}[\(L^p\)-convergence of MFL]
\label{thm:lp-convergence}
Assume \(F\) satisfies
\cref{eq:convex,eq:lip-in-m,eq:lsi,eq:first-der,eq:higher-der}.
Suppose \(h_0 \in L^{p_0}(m_\infty)\) for some \(p_0 > 1\).
For \(\rho' \in (0,\rho)\) and \(p \in \mathbb R\),
we set
\[
\tau_p = \begin{cases}
\frac{1}{4\rho'} \log \frac{(p-1)\vee 1}{(p_0-1)\wedge 1},
&\text{if $p \geqslant 0$}, \\
\frac{1}{4\rho'} \log \frac{2(1-p)}{(p_0 - 1)\wedge 1},
&\text{if $p < 0$}.
\end{cases}
\]
Then for all \(t \geq \tau_p\),
we have that \(h_t\) belongs to \(L^p(m_\infty)\) and its norm satisfies
\begin{align}
\bigl\lvert\log \lVert h_t \rVert_p\bigr\rvert
&\leq \biggl(\frac{2(1-p)}{p}\1_{p \in (0,1)} + \1_{p \not\in (0,1)}\biggr)
\nonumber \\
&\mathrel{\hphantom{\leq}}
\qquad \biggl( 1 + \frac{P(\alpha)}{8\varepsilon^2} \biggr)
H_1^{P(\alpha)/4\varepsilon}
\bigl( H_1^2 - 1\bigr) e^{-4(1-\varepsilon)\rho(t-\tau_p)} \nonumber \\
& \mathrel{\hphantom{\leq}} \quad \negmedspace {}
+ (p-2)_+ \1_{p > 0} \frac{p_0P(\alpha)\log\lVert h_0\rVert_{p_0}}
{16(p_0-1)\varepsilon(1-\varepsilon)} \cdot e^{(1-\varepsilon)^{-1}((p-1)\vee 1 )} e^{-4\rho t} \nonumber \\
& \mathrel{\hphantom{\leq}} \quad \negmedspace {}
+ (1/2 - p) \1_{p \leq 0} \frac{p_0P(\alpha)\log\lVert h_0\rVert_{p_0}}
{16(p_0-1)\varepsilon(1-\varepsilon)} \cdot e^{(1-\varepsilon)^{-1}(2(1-p) )} e^{-4\rho t} ,
\label{eq:lp-convergence}
\end{align}
where $\alpha = M^F_{mm} / \rho$,
$P(\alpha) = \alpha^2 + \alpha^3 + \alpha^4\!/2$,
and
\[
\log H_1 =
\biggl( 1 + \frac{p_0(2-p_0)_+P(\alpha)}{16(p_0-1)\varepsilon(1-\varepsilon)}
\biggr)
\log \lVert h_0\rVert_{p_0}.
\]
\end{thm}

\begin{rem}[Necessity of $L^{1+}$-initial condition]
\label{rem:necessity-lp0}
We here explain why it is necessary to assume \(m_{0} \in L^{p_0} (m_\infty)\)
for some $p_0 > 1$ in \cref{thm:lp-convergence}.
Let \(m_0 (\dd x) \propto \exp \bigl(- \sum_{\nu=1}^d |x^\nu|\bigr) \dd x\),
i.e., the \(d\)-tensorized exponential distribution
and \(F(m) = \frac 12 \int |x|^2 m(\dd x)\).
The Langevin dynamics \cref{eq:mf-sde} is nothing but Ornstein--Uhlenbeck:
\[
\dd X_t = - X_t \dd t + \sqrt 2 \dd W_t.
\]
The SDE is solved explicitly by
\[
X_t = e^{-t} X_0 + \sqrt 2 \int_0^t e^{-(t - s)} \dd W_s
\overset{d}{=} e^{-t} X_0 + \sqrt {1 - e^{-2t}} \mathcal N,
\]
where \(\mathcal N \sim \mathcal N(0,1)\) is a standard normal independent from
\(X_0\).
The Langevin has unique invariant measure \(m_\infty \propto \exp (- |x|^2/2)\),
i.e., the standard normal distribution in \(\mathbb R^d\).
The initial condition \(m_0\) lies in all \(\mathcal P_p\) for all \(p \geq 1\)
but \(m_0 / m_\infty\) does not belong to \(L^{p_0}\) for any \(p_0 > 1\).
And so is \(m_t\).
Indeed, for all \(\varepsilon > 0\),
\begin{align*}
\Expect [\exp (\varepsilon |X_t|^2)]
&= \Expect \bigl[\exp \bigl(\varepsilon (e^{-t}|X_0| + \sqrt{1 - e^{-2t}}
\mathcal N )^2\bigr)\bigr] \\
&\geq \Expect \biggl[\exp \biggl(\frac \varepsilon 2 (e^{-2t}|X_0|^2
- 2(1 - e^{-2t}) \mathcal N^2) \biggr)\biggr] \\
&= \Expect \biggl[\exp \biggl(\frac \varepsilon 2 e^{-2t}|X_0|^2\biggr)\biggr]
\Expect \bigl[\exp \bigl(-\varepsilon (1 - e^{-2t}) \mathcal N^2 \bigr)\bigr]
= +\infty.
\end{align*}
Here we used \((a+b)^2 \geq \frac 12 a^2 - b^2\) and the independence between
\(X_0\) and \(\mathcal N\).
This implies \(\int m_t(x) m_\infty(x)^{-\varepsilon}\dd x = +\infty\)
for all \(\varepsilon > 0\).
Let \(p > 1\). By Hölder's inequality we have
\begin{multline*}
\biggl(\int m_t(x)^p m_\infty(x)^{-(p-1)} \dd x \biggr)^{\!1/p}
\biggl(\int m_\infty(x)^{1-\varepsilon}\dd x\biggr)^{\!1 - 1/p} \\
\geq \int m_t(x) m_\infty(x)^{-\varepsilon (1 - 1/p)} \dd x = +\infty.
\end{multline*}
Hence \(\int m_t(x)^p m_\infty(x)^{-(p-1)} \dd x = +\infty\).
\end{rem}

As a by-product of our $L^p$-convergence result above,
we can use the transport method to show the following uniform-in-time
concentration of measure result.

\begin{thm}[Uniform-in-time concentration of measure]
\label{thm:mf-uniform-concentration}
Under the hypotheses of \cref{thm:lp-convergence},
for all \(\rho' \in (0, \rho)\)
there exist constants
\[
C_{\rho'} = C_{\rho'} (\rho, M^F_{mm}, p_0, \Vert h_0 \Vert_{p_0}),\quad
\tau_{\rho'} = \tau_{\rho'} (\rho, p_0)
\]
such that
for every \(1\)-Lipschitz function \(f : \mathbb R^d \to \mathbb R\),
every \(t \geq \tau_{\rho'}\) and every \(r \geq 0\),
\begin{equation}
\label{eq:mf-uniform-concentration}
m_t [ |f - \Expect_{m_t} f| \geq r ]
\leq 2 \exp \Bigl( - \rho' r^2 + C_{\rho'} e^{- 4 \rho' t} (r + 1) \Bigr).
\end{equation}
\end{thm}

The proofs of \cref{thm:lp-convergence,thm:mf-uniform-concentration}
are postponed to \cref{sec:proof-thm:lp-convergence}.

We further study the system of $N$ particles,
and show that its marginal distributions approximate $m_\infty^{\otimes N}$,
the $N$-tensorized mean field invariant measure,
at a uniform-in-$N$ exponential rate with a uniform-in-$N$ ``bias'',
whose precise meaning will be given below.

\begin{thm}[Uniform-in-\(N\) energy dissipation of particle systems]
\label{thm:ps-entropy-convergence}
Assume \(F\) satisfies \cref{eq:convex,eq:lip-in-m,eq:lsi,eq:first-der}.
If \(m^N_{t_0}\) belongs to \(\mathcal P_2(\mathbb R^{dN})\) and has finite
entropy
for some \(N \geq 2\) and \(t_0 \geq 0\),
then for all \(\rho' \in (0,\rho)\),
we have
\begin{align}
H(m^N_t | m_\infty^{\otimes N}) &\leq
\mathcal F^N(m^N_t) - N\mathcal F(m_\infty) \nonumber \\
&\leq \bigl(\mathcal F^N(m^N_{t_0}) - N\mathcal F(m_\infty)\bigr)
e^{ - (4\rho' - C_1N^{-1} ) (t - t_0) } \nonumber \\
&
\hphantom{\leq \bigl(\mathcal F^N(m^N_{t_0}) - N\mathcal F(m_\infty)\bigr)
e^{ - (4\rho' - C_1N^{-1} ) (t - t_0) }} \hspace{-5em}
+ \frac{C_2}{4\rho' - C_1N^{-1}},
\label{eq:ps-entropy-convergence}
\end{align}
for every $t \geq t_0$ and every $N > C_1 / 4\rho'$,
where the constants $C_1$, $C_2$ are defined by
\begin{align*}
C_1 &= M^F_{mm} \biggl( 16 + \frac{6 M^F_{mm} \rho'}{\rho(\rho -
\rho')}\biggr),
\\
C_2 &= d M^F_{mm} \biggl( 10 + \frac{3 M^F_{mm} \rho'}{\rho(\rho -
\rho')}\biggr).
\end{align*}
\end{thm}

The proof of \cref{thm:ps-entropy-convergence} is postponed
to \cref{sec:proof-thm:ps-entropy-convergence}.

\begin{rem}[Sharpness of the size of bias]
\label{rem:ps-entropy-convergence-sharp}
Let the initial condition $m^N_0$ of the $N$-particle system
be equal to $m^N_\infty$, the system's invariant measure.
By sending $t$ to infinity in \cref{eq:ps-entropy-convergence}, we have
\[
H(m^N_\infty | m_\infty^{\otimes N}) \leq \frac{C_2}{4\rho' - C_1 N^{-1}},
\]
provided that $\mathcal F^N(m^N_\infty) < +\infty$ and $N > C_1 / 4\rho'$.
Drawing an analogy with statistics, we will refer to the relative entropy
$H(m^N_\infty | m_\infty^{\otimes N})$ as the `bias'.
Then, the $O(1)$ order of the bias when $N \to +\infty$ is sharp
and we give an example attaining this order in the following.
Consider the mean field functional
\[
F(m) = \frac 12 \int x^2 m(\dd x) + \frac{\alpha}{2}
\biggl( \int x m(\dd x) \biggr)^{\!2}
\]
with $\alpha \geq 0$.
We can easily verify all our assumptions on $F$.
The mean field invariant measure is nothing but the $d$-dimensional standard
Gaussian variable:
\[
m_\infty (\dd x) = (2\pi)^{-d/2} \exp\biggl( -\frac{|x|^2}{2} \biggr) \dd x,
\]
and the invariant measure of the $N$-particle system reads
\[
m^N_\infty (\dd\mathbf x) = (2\pi)^{-dN/2} (\det A_N)^{1/2}
\exp \Biggl( - \frac 12 \sum_{i=1}^N |x^i|^2
- \frac{\alpha}{2N} \biggl( \sum_{j=1}^N x^i \biggr)^{\!2} \Biggr)\dd\mathbf x,
\]
where $A_N$ is the $Nd \times Nd$ matrix whose $d \times d$ blocks read
\[
(A_N)_{ij} = \begin{cases}
\bigl(1 + \frac{\alpha}{N}\bigr)\mathbf 1_{d\times d} & \text{if}~i=j, \\
\frac{\alpha}{N}\mathbf 1_{d \times d} & \text{if}~i\neq j.
\end{cases}
\]
Especially, we have $\mathcal F^N(m^N_\infty) < +\infty$.
By diagonalizing $A_N$, we can obtain $\det A_N = (1 + \alpha)^d$.
Hence, the relative density between $m^N_\infty$ and $m_\infty^{\otimes N}$
reads
\[
\frac{\dd m^N_\infty}{\dd m_\infty^{\otimes N}} (\mathbf x)
= (1 + \alpha)^{d/2}
\exp \Biggl(- \frac{\alpha}{2N} \biggl( \sum_{j=1}^N x^i \biggr)^{\!2} \Biggr),
\]
and the relative entropy satisfies
\begin{align*}
H(m^N_\infty | m_\infty^{\otimes N})
&= \Expect^{\mathbf X \sim m^N_\infty}
\biggl[ \log \frac{dm^N_\infty}{dm_\infty^{\otimes N}} (\mathbf X) \biggr] \\
&= \frac{d}{2} \log (1 + \alpha) - \frac{\alpha}{2N}
\Expect^{\mathbf X \sim m^N_\infty}
\Biggl[ \biggl( \sum_{i=1}^N X^i \biggr)^{\!2} \Biggr] \\
&= \frac{d}{2} \log (1 + \alpha) - \frac{d\alpha}{2(1+\alpha)}.
\end{align*}
So the $O(1)$ order in $N$ of the bias in
\cref{eq:ps-entropy-convergence} is sharp.
\end{rem}

Finally, we study the propagation of chaos phenomenon.
On finite horizon we use the classical arguments
of synchronous coupling and Girsanov's theorem
to show that the distance between the particle system $m^N_t$
and the tensorized mean field system $m_t^{\otimes N}$
grows at most exponentially,
in the sense of Wasserstein distance and relative entropy.
On the other hand, for large time,
we control the distance using the long time behavior proved in
\cref{thm:mf-entropy-convergence,%
thm:lp-convergence,thm:ps-entropy-convergence}.

\begin{thm}[Wasserstein and entropic propagation of chaos]
\label{thm:poc}
Assume \(F\) satisfies \cref{eq:convex,eq:lip-in-m,eq:lsi,eq:first-der}.
Suppose \(m_0\) belongs to \(\mathcal P_2(\mathbb R^d)\),
\(m^N_0\) belongs to \(\mathcal P_2(\mathbb R^{dN})\)
and they both have finite entropy
for some \(N \geq 2\).

\begin{itemize}
\item Then for all \(\rho' \in (0,\rho)\), we have
\begin{multline}
\label{eq:w-poc}
\rho W_2^2(m^N_t, m^{\otimes N}_t)
\leq
2N \bigl(\mathcal F(m_0) - \mathcal F(m_\infty)\bigr) e^{-4\rho t} \\
+ 2\bigl(\mathcal F^N(m_0^N) - N\mathcal F(m_\infty)\bigr)
e^{-(4\rho' - C_1N^{-1}) t}
+ \frac{2C_2}{4\rho' - C_1N^{-1}},
\end{multline}
for every $t \geq 0$ and every $N > C_1 / 4\rho' $,
where the constants $C_1$, $C_2$ are the same
as in \cref{thm:ps-entropy-convergence}.
If additionally $m_0 \in \mathcal P_6(\mathbb R^d)$,
then we have
\begin{multline}
\label{eq:w-poc-finite}
W_2^2 (m^N_t, m_t^{\otimes N})
\leq e^{C_4 t} W_2^2 (m^N_0, m_0^{\otimes N}) \\
+ N C_5 (e^{C_4 t} - 1)
\bigl(v_6(m_0)^{1/3} + 1\bigr)\delta_d (N),
\end{multline}
for every $t \geq 0$,
where $C_4 = \max \bigl( 1 + 3(M^F_{mx})^2 + 3(M^F_{mm})^2,
2 M^F_{mx} + 4d/3 + 16/3\bigr)$,
and $C_5$ is a constant depending only on $M^F_{mx}$, $M^F_{mm}$ and $d$,
the term $v_6(m_0)$ is defined by
$v_6(m_0) \coloneqq \int \bigl| x - \int x' m_0(\dd x') \bigr|^6 m_0(\dd x)$
and the term $\delta_d(N)$ is defined by
\[
\delta_d (N) \coloneqq \begin{cases}
N^{-1/2} & \text{if}~d < 4, \\
N^{-1/2} \log (1 + N) & \text{if}~d = 4, \\
N^{-2/d} & \text{if}~d > 4.
\end{cases}
\]

\item If additionally \cref{eq:higher-der} holds and
\(h_0 \in L^{p_0} (m_\infty)\) for some \(p_0 > 1\),
then we have
\begin{multline}
\label{eq:h-poc}
H(m^N_t | m^{\otimes N}_t)
\leq
N C_3 e^{-4\rho't} \\
+ 2\bigl(\mathcal F^N(m_0^N) - N\mathcal F(m_\infty)\bigr)
e^{-(4\rho' - C_1N^{-1}) t}
+ \frac{2C_2}{4\rho' - C_1N^{-1}},
\end{multline}
for every $t \geq \tau$ and every $N > C_1 / 4\rho'$,
for some constants $C_3$, $\tau \geq 0$
depending only on $\rho$, $\rho'$, $M^F_{mm}$, $p_0$ and
$\Vert h_0 \Vert_{L^{p_0}(m_\infty)}$.
If additionally $m_0 \in \mathcal P_6(\mathbb R^d)$
and $H(m_0^N | m_0^{\otimes N})$ is both finite,
then we have
\begin{multline}
\label{eq:h-poc-finite}
H(m^N_t | m_t^{\otimes N})
\leq H(m^N_0 | m_0^{\otimes N}) \\
+ N C_5 (e^{C_4 t} - 1)
\bigl(v_6(m_0)^{1/3}+ 1\bigr)\delta_d (N),
\end{multline}
for every $t \geq 0$,
for possibly different constants $C_4$, $C_5 > 0$ depending
on $M^F_{mx}$, $M^F_{mm}$ and $d$.
\end{itemize}
\end{thm}

If the initial error is zero, i.e., $m_0^N = m_0^{\otimes N}$,
we obtain the following result by combining the finite-time and long-time
estimates, as in the proof of Corollary 5 of \cite{guillin2021uniform}.

\begin{cor}
\label{cor:poc}
Assume \(F\) satisfies \cref{eq:convex,eq:lip-in-m,eq:lsi,eq:first-der}.
Suppose $m_0 \in \mathcal P_6(\mathbb R^d)$, $m_0$ has finite entropy,
and $m_0^N = m_0^{\otimes N}$.
Then there exist constants $C$, $N_0 > 0$,
depending on $\rho$, $M^F_{mm}$, $M^F_{mx}$, $m_0$ and $d$,
such that for all $N \geq N_0$,
\begin{align}
\label{eq:w-poc-unif}
\sup_{t \in [0,\infty)}
\frac 1N W_2^2(m^N_t, m_t^{\otimes N})
&\leq \frac{C}{N^\kappa}
\intertext{where $\kappa = \min (2\rho/C_4, 1) / (d\vee 4)$
with $C_4$ being the constant in the Wasserstein case
of \cref{thm:poc}.
If additionally $F$ satisfies \cref{eq:higher-der},
we have as well}
\label{eq:h-poc-unif}
\sup_{t \in [0,\infty)}
\frac 1N H(m^N_t | m_t^{\otimes N})
&\leq \frac{C}{N^\kappa}
\end{align}
for every $N \geq N_0$, with the constants
$C$, $\kappa$, $N_0 > 0$ redefined accordingly.
\end{cor}

The proofs of \cref{thm:poc,cor:poc} are postponed to \cref{sec:poc}.
The rate $\kappa$ obtained in the corollary above
seems to be highly optimal compared to the $O(1/N)$ rate
in \cref{thm:ps-entropy-convergence}.
This is due to the fact that, for finite time, we do not exploit
at all the coercive structure of the MFL.
We note that it is recently shown in \cite{delarue2021uniform} that
if we consider a weaker distance and work under stronger regularity conditions,
then the optimal $O(1/N)$ rate can be achieved even
when the supremum over all time is taken.

\paragraph{Comments on the assumptions.}
The conditions \cref{eq:lip-in-m,eq:first-der} ensure that the drift is jointly
Lipschitz continuous in measure and space,
which guarantees the wellposedness of the mean field and the particle system
dynamics \cref{eq:mf-sde,eq:ps-sde}.
This also implies that the flow is \(AC^2\) in \(L^2\)-Wasserstein space (refer
to \Cref{defn:ac2}),
which coincides with the type of curves studied in \cite[Chapter 8]{gf}.
In particular, the ``chain rule'' holds true,
which yields immediately the energy dissipation
\cref{eq:mf-energy-decrease,eq:ps-energy-decrease}.

The assumptions \cref{eq:convex,eq:lsi},
which have already appeared in the previous works
\cite{chizat2022mean,nitanda2022convex},
are key to the exponential convergence of relative entropy of the MFL.
They are also used in this work, along with \cref{eq:lip-in-m}, to show the
exponential entropic convergence of the particle system in
\cref{thm:ps-entropy-convergence}.

The condition \cref{eq:higher-der} is technical in that it does
not contribute to any constants in our results.
This condition allows us to obtain
a simple ``standard algebra'' of the time-dependent semigroup induced by the MFL
and to justify easily the computations in \(L^p\) spaces needed to prove
\cref{thm:lp-convergence},
which is then used to show \cref{thm:poc,cor:poc}.
It is possible that our results can also be obtained without the higher-order
bounds (for example, by an approximation argument).
We, however, choose to work in this setting to avoid excessive technicalities.

\section{Applications}
\label{sec:app}

\subsection{Sufficient conditions for functional convexity}
\label{sec:functional-convexity}

We propose two criteria for the convexity of mean field functionals.
The first criterion treats translationally invariant two-body interactions,
i.e., energy functionals of the form:
\begin{equation}
\label{eq:2-body}
F_\textnormal{Int} (m) = \frac 12 \iint V(x - y) m(\dd x) m(\dd y).
\end{equation}
We have the following modified version of Bochner's theorem.

\begin{thm}[Bochner]
\label{thm:bochner}
Let \(V : \mathbb R^d \to \mathbb R\) be a bounded, continuous
and even function. Then, the following conditions are equivalent:
\begin{enumerate}[label=(\roman*)]
\item The functional $F_\textnormal{Int}$, defined by \cref{eq:2-body},
is convex on $\mathcal P(\mathbb R^d)$.
\item For all signed measure $\mu$ on $\mathbb R^d$ with zero net mass, i.e.,
$\int\dd\mu = 0$, we have $\iint V(x - y) \mu(\dd x) \mu(\dd y) \geq 0$.
\item The Fourier transform $\hat V$ of $V$
is the sum of a finite and positive measure on $\mathbb R^d \setminus \{0\}$
and a scalar multiple of the Dirac mass $\delta_0$ at zero.
\end{enumerate}
\end{thm}

The proof of this modified version of Bochner's theorem
is postponed to \cref{sec:bochner}.

\begin{exm}[Regularized Coulomb]
\label{exm:reg-coulomb}
It is well-known that in dimension \(d \geq 3\) the Coulomb potential
\(V_\textnormal{C}(x) = 1\big/\bigl(d(d-2)c_d |x|^{d-2}\bigr)\)
is the fundamental solution to Laplace's equation, that is to say,
\begin{equation}
\label{eq:laplace}
\Delta V_\textnormal{C} = - \delta_0.
\end{equation}
Hence its Fourier transform \(\hat V_\textnormal{C}\) verifies \(\hat
V_\textnormal{C}(k) = (2\pi)^{-d/2}|k|^{-2} \geq 0\).
However \(\hat V_\textnormal{C} \not\in L^1(\mathbb R^d)\) and
\cref{thm:bochner} does not apply (which is consistent with the singularity of
\(V_\textnormal{C}\) at \(0\)).
To solve this problem, we propose the regularization
\[
\hat V_\textnormal{RC}(k) = \frac{e^{-r_0 |k|}}{(2\pi)^{d/2}|k|^2}
\]
for some \(r_0 > 0\).
Its Fourier inverse \(V_\textnormal{RC} : \mathbb R^d \to \mathbb R\) is then
indeed a bounded continuous function
and has the explicit expression for \(d = 3\):
\[
V_\textnormal{RC}(x) = \int \frac{e^{-r_0 |k|} e^{ik\cdot x}}{(2\pi)^3 |k|^2}
\dd^3 k
= \begin{cases} \arctan (|x|/r_0) (2\pi^2 |x|)^{-1} &\text{if}~x \neq 0, \\
(2\pi^2 r_0)^{-1} &\text{if}~x = 0.
\end{cases}
\]
Note that when \(r_0 \to 0\), we have \(V_\textnormal{RC}(x) \to
V_\textnormal{C}(x)\) for every \(x \in \mathbb R^d\).
The functional
\begin{multline}
\label{eq:reg-coulomb}
F_\textnormal{RC}(m)
= \frac 12 \iint V_\textnormal{RC} (x - y) m(\dd x) m(\dd y) \\
= \frac 12 \iint \frac{1}{2\pi^2} \frac{\arctan(|x - y|/r_0)}{|x - y|} m(\dd x)
m(\dd y)
\end{multline}
is well defined and convex on \(\mathcal P(\mathbb R^3)\) by \cref{thm:bochner}.
\end{exm}

\begin{rem}[Exclusion of two notions of convexity]
\label{rem:exclusion-two-convexities}
If the functional \(F_\textnormal{Int}\) satisfies the conditions of
\cref{thm:bochner}, we know
\[
2V(0) - V(s) - V(-s) = \frac{2}{(2\pi)^{d/2}}
\int_{\mathbb R^d} \bigl( 1 - \cos (k \cdot s) \bigr) \hat V(\dd k) \geq 0.
\]
If the function $V$ is not constant, then there exists some
$s_0 \in \mathbb R^d$ such that $V(s_0) \neq V(0)$.
The evenness of $V$ implies $V(-s_0) = V(s_0)$ and, therefore,
$V(s_0) = V(-s_0) < V(0)$.
In particular, \(V\) is not convex,
and the functional \(F_\textnormal{Int}\) cannot be geodesically convex.
In other words,
the only functionals of form \cref{eq:2-body} with
continuous, bounded and even \(V\)
that are both functionally and geodesically convex
are constant functionals.
\end{rem}

\begin{rem}
Other regularizations preserving the positivity of the Coulomb potential can
also be possible.
For example we can convolute Laplace's equation \cref{eq:laplace} with a heat
kernel
\(\rho^\varepsilon : x \mapsto (2\pi\varepsilon)^{-d/2}
\exp \bigl(- (2\varepsilon)^{-1} x^2\bigr)\) to obtain
\[
\Delta V'_\textnormal{RC} = \Delta (V_\textnormal{C} \star \rho^\varepsilon) =
- \rho^\varepsilon.
\]
The Fourier transform of \(V'_\textnormal{RC}\) reads
\[
\hat V'_\textnormal{RC} (k) = \frac{\hat \rho^\varepsilon (k)}{|k|^2}
= \frac{e^{-2\pi^2\varepsilon |k|^2}}{(2\pi)^{d/2} |k|^2},
\]
which is positive and \(L^1\)-integrable.
The main reason for choosing the regularization in \cref{exm:reg-coulomb} is
that it allows for the simple expression given in \cref{eq:reg-coulomb}
in three dimensions.
\end{rem}

The second criterion is an analogue of the property of convex functions under
composition.

\begin{prop}
\label{prop:convexity-composition}
Let \(X\) be a Banach space.
If \(V : \mathbb R^d \to X\) is a function of quadratic growth and \(g : X \to
\mathbb R\) is convex,
then the functional \(F : \mathcal P_2(\mathbb R^d) \to \mathbb R\) defined by
\[
F(m) = g \biggl( \int V(x) m(\dd x) \biggr)
\]
is convex.
\end{prop}

\begin{proof}
Immediate.
\end{proof}

\begin{exm}[\(L^2\)-loss of two-layer neural networks]
\label{exm:nn}
We first explain the structure of two-layer neural networks and then introduce
the mean field model for it.
Consider an \emph{activation function} \(\varphi : \mathbb R \to \mathbb R\)
satisfying
\begin{equation}
\label{eq:nn-activation}
\begin{gathered}
\text{$\varphi$ is continuous and non-decreasing,} \\
\lim_{x\to -\infty} \varphi(x) = 0, \quad \lim_{x\to +\infty} \varphi(x) = 1,
\end{gathered}
\end{equation}
Define \(S = \mathbb R \times \mathbb R^d \times \mathbb R\), where the
\emph{neurons} take values.
For each neuron \(\theta = (c,a,b) \in S\)
we define the \emph{feature map}:
\begin{equation}
\label{eq:nn-feature}
\mathbb R^d \ni z \mapsto \Phi(\theta; z) \coloneqq \ell(c) \varphi (a \cdot z
+ b) \in \mathbb R,
\end{equation}
where \(\ell : \mathbb R \to [-L, L]\) is a \emph{truncation function} with the
\emph{truncation threshold} \(L \in (0,+\infty]\).
Such truncation has been considered in recent papers
\cite{HRSS19,nitanda2022convex}.
The two-layer neural network is nothing but
the averaged feature map parameterized by
\(N\) neurons
\(\theta^1, \ldots, \theta^N \in S\):
\begin{equation}
\label{eq:nn-ps-feature}
\mathbb R^d \ni z \mapsto \Phi^N (\theta^1,\ldots,\theta^N ; z) = \frac 1N
\sum_{i=1}^N \Phi(\theta^i; z) \in \mathbb R.
\end{equation}
The training of neural network aims to minimize the distance between the
averaged output \cref{eq:nn-ps-feature} and a (only empirically known)
\emph{label} function \(f : \mathbb R^d \to \mathbb R\),
i.e.
\begin{equation}
\label{eq:nn-ps-optimize}
\inf_{\mathbf (\theta^1,\ldots,\theta^N) \in S^N}
\mathbf d\bigl( f , \Phi^N(\theta^1,\ldots,\theta^N; \cdot)\bigr)
\end{equation}
for some loss functional \(\mathbf d\).
In this paper, we use the \(L^2(\mu)\)-norm as the loss functional
where \(\mu \in \mathcal P(\mathbb R^d)\) represents the \emph{feature}
distribution.
In this way, the objective function of the minimization reads
\begin{equation}
\label{eq:nn-ps-loss}
F_\textnormal{NNet}^N (\theta^1,\ldots,\theta^N) =
\frac N2 \int
\bigl\lvert f(z) - \Phi^N(\theta^1,\ldots,\theta^N; z)\bigr\rvert^2 \mu(\dd z).
\end{equation}
To fit the problem to our theoretical framework,
we assume that
the feature map \(\Phi : S \times \mathbb R^d \to \mathbb R\) satisfies
\begin{align*}
\forall \theta \in S,\qquad &\Phi(\theta; \cdot) \in L^2(\mu), \\
\exists C > 0,~\forall \theta \in S,\qquad & \Vert\Phi(\theta;
\cdot)\Vert_{L^2(\mu)} \leq C (1 + |\theta|^2).
\end{align*}

Now we present the mean field formulation of two-layer neural networks.
Let \(\mathcal P_2(S)\) be the space of probability measures
on \(S\) of finite second moment
and define the class of functions representable by the mean field neural
network by:
\begin{equation}
\label{eq:nn-mf-expressiveness}
\mathcal N_{\varphi,\ell} = \{ h : \mathbb R^d \to \mathbb R : \exists m \in
\mathcal P_2(S),~\forall x \in \mathbb R^d,\quad h(x) = \Expect^{\Theta \sim m}
[\Phi(\Theta; x)]\}.
\end{equation}
In particular the \(N\)-neuron output functions defined in
\cref{eq:nn-ps-feature} belong to this class since
\[
\Phi^N(\theta^1,\ldots,\theta^N; \cdot)
= \Expect^{\Theta \sim \frac{1}{N}\sum_{i=1}^N \delta_{\theta^i}} [\Phi
(\Theta; \cdot)].
\]
Instead of the finite-dimensional optimization \cref{eq:nn-ps-optimize},
we consider the following mean field optimization:
\begin{equation}
\label{eq:nn-mf-loss}
\begin{aligned}
\inf_{\mathcal{P}_2(S)} &F_\textnormal{NNet}(m), \\
\text{where } &F_\textnormal{NNet}(m) \coloneqq
\mathbf d\bigl(f, \Expect^{\Theta\sim m} [\Phi(\Theta; \cdot)]\bigr)
= \frac 12 \int \bigl\lvert f(z) - \Expect^{\Theta \sim m} [\Phi(\Theta; z)]
\bigr\rvert^2 \mu(\dd z).
\end{aligned}
\end{equation}
The functional \(F_\textnormal{NNet}\) is convex by
\cref{prop:convexity-composition}
since \(F_\textnormal{NNet} (m) = g \bigl( \int V(\theta) m(\dd\theta) \bigr)\)
with \(V : S \ni \theta \mapsto
\bigl(z \mapsto \Phi(\theta; z)\bigr) \in L^2(\mu)\)
and \(g : L^2(\mu) \ni h \mapsto \Vert f - h \Vert_{L^2(\mu)}^2 \in \mathbb R\).
\end{exm}

\begin{rem}[Motivation of mean field formulation]
The \(N\)-neuron problem \cref{eq:nn-ps-loss} is non-convex due to the
non-linear activation function \(\varphi\).
Inspired by the fact that the width \(N\) of two-layer neural networks is
usually large in practice,
the authors of \cite{mei2019mean,chizat2018global,rotskoff2018neural,HRSS19}
consider the mean field formulation of neural networks
which convexifies the original problem.
\end{rem}

\begin{rem}[Absence of geodesic convexity]
\label{rem:nn-geodesic-convexity}
We highlight here that
if \(F_\textnormal{NNet}\) is geodesically convex and regular enough,
then the \(N\)-neuron problem \(F^N_\textnormal{NNet}\) is convex,
which is not true.
Hence by contradiction \(F_\textnormal{NNet}\) has no geodesic convexity.
Indeed, suppose \(F_\textnormal{NNet}\) is geodesically convex.
Note that \(t \mapsto \frac 1N \sum_{i=1}^N \delta_{\theta^i + tv^i}\)
is a geodesic in \((\mathcal P_2, W_2)\) in a neighborhood of \(t = 0\)
if \(\theta_i\) are distinct from each other
(as the pairing \((\theta^i, \theta^i + tv^i),~i=1,\ldots,N\) verifies cyclical
monotonicity for \(t\) small enough).
By the geodesic convexity of \(F_\textnormal{NNet}\)
and the relation
\(F^N_\textnormal{NNet} (\theta^1,\ldots,\theta^N)
= N F_\textnormal{NNet}\bigl( \frac 1N \sum_{i=1}^N \delta_{\theta^i}\bigr)\),
we obtain the local convexity of \(F^N_\textnormal{NNet}\) on the set
\[
S^N \setminus \Delta^N \coloneqq S^N \setminus \{ (\theta_1,\ldots,\theta_N)
\in S^N : \exists i \neq j,\quad \theta_i = \theta_j \}.
\]
If \(F^N_\textnormal{NNet}\) is additionally \(C^2\),
the local convexity implies \(\nabla^2 F^N_\textnormal{NNet} \geq 0\) on \(S^N
\setminus \Delta^N\) and by density \(\nabla^2 F^N_\textnormal{NNet} \geq 0\)
everywhere.
Therefore \(F^N_\textnormal{NNet}\) is convex on \(S^N\).
\end{rem}

\begin{rem}[Expressiveness of truncated networks]
\label{rem:nn-expressiveness}
It is well known that two-layer neural networks are universal approximators,
that is,
they can approximate any continuous function on \(\mathbb R^d\) arbitrarily well
with respect to the compact-open topology (\cite[Theorem 2.4]{hor91}).
This implies that the infimum in \cref{eq:nn-mf-loss} is zero if \(\mu\) is
compactly supported and no truncation is present
(that is, \(L = +\infty\) and \(\ell\) is the identity function).
However, if a truncation with \(L < +\infty\) is applied,
all functions \(h \in \mathcal N_{\varphi,\ell}\) satisfy the bound \(\Vert h
\Vert_\infty \leq L\)
and therefore cannot approximate well functions that exceed \(L\).
However, Barron's theorem \cite[Theorem 2]{Barron} says that
if a function $f$ satisfies
\[
f(x) = f(0) + \int (e^{i \omega \cdot x} - 1) F(\dd\omega)
\]
for every $x \in B(0,R)$, for some complex-valued measure $F$,
and if there exists $c_+$, $c_- \in \mathbb R$ such that
$\ell(c_+) = L$ and $\ell(c_-) = -L$,
and that
\[
L \geq R \int \lvert \omega \rvert \lvert F(\dd\omega) \rvert
+ \lvert f(0)\rvert,
\]
then the best approximation error
\[
\inf_{\Phi \in \mathcal N_{\varphi, \ell}} \lVert f - \Phi \rVert_{L^2(\mu)} = 0
\]
for every probability measure $\mu$ supported in $B(0,R)$.
\end{rem}

\subsection{Examples of MFL dynamics}
\label{sec:exm}

We construct MFL dynamics for the two examples discussed earlier
and demonstrate that our theorems are applicable in both cases.
To verify the LSI condition \cref{eq:lsi} we will use the following results.

\begin{prop}
\label{prop:lsi}
Let \(\mu (\dd x) = e^{-V(x)}\dd x\) be a probability measure in \(\mathbb R^d\)
for some \(V \in C^2(\mathbb R^d)\).
\begin{itemize}
\item (Bakry--Émery \cite{bakryemery}) If \(\nabla^2 V \geq \kappa\) then
\(\mu\) satisfies a \(\kappa/2\)-LSI.
\item (Holley--Stroock \cite{holleystroock}) If \(V = V_1 + V_2\),
where \(e^{-V_1}\)
is the density of a probability measure satisfying an \(\rho\)-LSI and \(V_2\) is
bounded with oscillation \(\osc V_2\),
then \(\mu\) satisfies a \(\rho \exp(-\osc V_2)\)-LSI.
\item (Aida--Shigekawa \cite{aidashigekawa}) If \(V_2\) in the previous
statement is Lipschitz-continuous instead of bounded, then \(\mu\) satisfies an
LSI as well.
\end{itemize}
\end{prop}

\begin{exm}[MFL for regularized Coulomb system]
\label{exm:gd-reg-coulomb}
Let \(\lambda > 0\).
Define
\begin{equation}
\label{eq:ext}
F_\textnormal{Ext} (m) = \frac{\lambda}{2} \int |x|^2 m(\dd x).
\end{equation}
We consider the functional
\(F = F_\textnormal{RC} + F_\textnormal{Ext}\)
where \(F_\textnormal{RC}\) is defined in \cref{eq:reg-coulomb}.
By the discussions in \cref{exm:reg-coulomb} the functional \(F\) satisfies the
convexity condition \cref{eq:convex}.
Its linear functional derivative reads
\[
\frac{\delta F}{\delta m} (m,x) = \int V_\textnormal{RC} (x - y) m(\dd y)
+ \frac 12 \lambda |x|^2
\]
and its intrinsic derivative reads
\(D_m F(m,x) = \int \nabla V_\textnormal{RC} (x - y) m(\dd y) + \lambda x\).
The conditions \cref{eq:lip-in-m,eq:first-der,eq:higher-der} are satisfied
because
\[
\Vert \nabla^n V_\textnormal{RC} \Vert_\infty
\leq \frac{1}{(2\pi)^{d/2}}\int |k|^n \hat V_\textnormal{RC} (\dd k)
= \int |k|^n \frac{e^{- r_0 |k|}}{(2\pi)^d |k|^2} \dd^d k < +\infty
\]
for all \(n \geq 0\) (and \(d \geq 3\)).
In particular, the bound in \cref{eq:lip-in-m} is verified by \(M^F_{mm} =
\Vert \nabla^2 V_\textnormal{RC} \Vert_\infty\).
For the uniform LSI, we can apply Holley--Stroock or Aida--Shigakawa,
since the first term in \(\frac{\delta F}{\delta m}\) is uniformly bounded and
uniformly Lipschitz
and the second term verifies the Bakry--Émery condition.
The LSI constant given by Holley--Stroock has the simple expression in three
dimensions
\(\rho = \lambda \exp (- \osc V_\textnormal{RC}) / 2
= \lambda \exp (- 1/2\pi^2 r_0) / 2\).
The \(L^{1+}\)-integrability of the initial condition, needed by
\cref{thm:lp-convergence} and the second part of \cref{thm:poc},
is verified once we have
\begin{equation}
\label{eq:init-necessary-l1+}
\exists C,\varepsilon > 0,~\forall x \in \mathbb R,\qquad m_0(x) \leq C e^{-
\varepsilon |x|^2}.
\end{equation}
However, as the regularization parameter \(r_0\) approaches \(0\),
we observe \(\rho \to 0\) and \(M^F_{mm} \to +\infty\),
suggesting our method is not suitable for the unregularized Coulomb interaction.
We refer readers to
\cite{bresch2019modulated,bresch2019mean,rosenzweig2023global,%
decourcel2023sharp}
for recent developments on the noised gradient flow of Coulomb (and more
generally, Riesz) particle systems,
where the \emph{modulated free energy} is used to tackle the singularity in the
interactions.
\end{exm}

\begin{exm}[MFL for two-layer neural networks]
\label{exm:nn-ngd}
Recall the mean field two-layer neural networks in \cref{exm:nn}.
Suppose
\begin{itemize}
\item the truncation \(L\) is finite;
\item the activation and truncation functions $\varphi$, $\ell$
have bounded derivatives of up to fourth order;
\item the feature distribution \(\mu\) has finite second moment;
\item the label function \(f\) belongs to \(L^2(\mu)\).
\end{itemize}
On top of the mean field optimization problems \cref{eq:nn-mf-loss},
we add the quadratic regularizer \(F_\textnormal{Ext}\) in \cref{eq:ext} to the
loss, as for the Coulomb system.
Then the function and the functional to optimize read
\begin{align*}
F^N (\theta^1,\ldots,\theta^N)
&= \frac N2 \int \biggl|f(z) - \frac 1N \sum_{i=1}^N \Phi(\theta^i; z)\biggr|^2
\mu(\dd z) + \frac{\lambda}{2} \sum_{i=1}^N |\theta^i|^2, \\
F (m) &= \frac 12 \int \bigl\lvert f(z) - \Expect^{\Theta \sim m}
[\Phi(\Theta; z)] \bigr\rvert^2
\mu(\dd z) + \frac{\lambda}{2} \int |\theta|^2 m(\dd\theta).
\end{align*}
The \(N\)-neuron loss can be recover from the mean field loss by
\(F^N(\theta^1,\ldots,\theta^N) = N F\bigl( \frac 1N \sum_{i=1}^N
\delta_{\theta^i}\bigr)\).
We verify the assumptions of our theorems one by one.
The functional convexity of \(F = F_\textnormal{NNet} + F_\textnormal{Ext}\) is
already proved in \cref{exm:nn}.
The linear functional derivative of \(F\) reads
\[
\frac{\delta F}{\delta m}(m, \theta) =
- \int \bigl(f(z) - \Expect^{\Theta \sim m} [\Phi (\Theta; z)]\bigr)
\Phi(\theta; z) \mu(\dd z)
+ \frac{\lambda}{2} |\theta|^2.
\]
The first term on the right hand side is uniformly bounded:
for every \(m \in \mathcal P_2(S)\) and every \(\theta \in S\),
\[
\biggl\lvert\int \bigl(f(z) - \Expect^{\Theta \sim m} [\Phi (\Theta; z)]\bigr)
\Phi(\theta; z) \mu(\dd z)\biggr\rvert
\leq (\Vert f\Vert_{L^1(\mu)} + \Vert \ell \Vert_\infty) \Vert \ell
\Vert_\infty.
\]
Hence by Holley--Stroock the uniform LSI condition \cref{eq:lsi} is satisfied
with the constant
\[
\rho = \frac{\lambda}{2} \exp\bigl( - 2 (\Vert f \Vert_{L^1(\mu)}
+ \Vert\ell\Vert_\infty) \Vert\ell\Vert_\infty\bigr).
\]
The intrinsic derivative of \(F\) reads
\[
D_m F(m, \theta) = - \int
\bigl(f(z) - \Expect^{\Theta \sim m} [\Phi (\Theta; z)]\bigr)
\frac{\partial \Phi}{\partial \theta}(\theta; z) \mu(\dd z) + \lambda \theta,
\]
where the partial derivative of the feature map \(\Phi\), defined in
\cref{eq:nn-feature}, reads
\[
\frac{\partial \Phi}{\partial c}(\theta; z) = \ell' (c) \varphi (a \cdot z +
b),~
\frac{\partial \Phi}{\partial a}(\theta; z) = \ell (c) \varphi' (a \cdot z + b)
z,~
\frac{\partial \Phi}{\partial b}(\theta; z) = \ell (c) \varphi' (a \cdot z + b)
\]
for \(\theta = (c,a,b) \in S\).
Similarly we obtain the second order intrinsic derivative:
\(D_m^2 F(m,\theta,\theta') = \int \frac{\partial \Phi}{\partial \theta}
(\theta; z) \otimes \frac{\partial \Phi}{\partial \theta} (\theta'; z)
\mu(\dd z)\).
Its \(2\)-norm satisfies the bound
\(|D_m^2 F(m,\theta,\theta')|^2_2
\leq \Vert \ell' \Vert_\infty^2 + \Vert \ell \Vert_\infty^2 \Vert \varphi'
\Vert_\infty^2 \bigl(1 + M_2(\mu)\bigr)\),
where \(M_2 (\mu) = \int |z|^2 \mu(\dd z)\) is the second moment of \(\mu\).
Thanks to the Kantorovich duality and the Cauchy--Schwarz inequality,
the \(W_1\)-Lipschitz constant of \(m \mapsto D_m F(m,x)\) can be given by
\[
M^F_{mm} = \Bigl(\Vert \ell' \Vert_\infty^2
+ \Vert \ell \Vert_\infty^2 \Vert \varphi' \Vert_\infty^2
\bigl(1 + M_2(\mu)\bigr) \Bigr)^{\!1/2}.
\]

So \(D_m F\) satisfies the condition \cref{eq:lip-in-m}.
Since $\ell$, $\varphi$ have bounded derivatives of up to fourth order, the
derivatives \(\nabla^k D_m F(m, \theta)\)
for $k = 1$, $2$, $3$ are also uniformly bounded.
Thus the technical conditions \cref{eq:first-der,eq:higher-der} are also
satisfied.
Finally, the \(L^{1+}\)-integrability of the initial value \(m_0\)
is verified once we require the pointwise Gaussian bound
\cref{eq:init-necessary-l1+} on the density of \(m_0\).
\end{exm}

\begin{rem}[Link to practice]
In the training of neural networks, the measure \(\mu\) is an empirical measure
\(\frac 1K \sum_{k=1}^K \delta_{z_k}\)
and on the feature points \(\{z_k\}_{k=1}^K\) the labels are known \(f(z_k) =
y_k\).
This collection of pairs \(\{z_k,y_k\}_{k=1}^K\)
are the available training data.
In practice, instead of the mean field dynamics,
we can only simulate the corresponding \(N\)-particle system.
In other words, we calculate the \(N\)-neuron SDE
\begin{equation}
\label{eq:nn-ps-sde}
\dd \Theta^i_t
= \frac{1}{K} \sum_{k=1}^K
\bigl(y_k - \Phi^N( \Theta^1_t, \ldots, \Theta^N_t; z_k)\bigr)
\frac{\partial \Phi}{\partial \theta} (\Theta^i_t; z_k) \dd t
- \lambda \Theta^i_t \dd t
+ \sigma \dd W^i_t,
\end{equation}
for \(i = 1,\ldots,N\).
The first drift term of the diffusion
is the gradient \(\nabla_{\theta^i} F^N(\Theta^1_t, \ldots, \Theta^N_t)\),
so the time-discretization of this diffusion is nothing but the \emph{noisy
gradient descent} (NGD) algorithm for training neural networks.
We refer readers to \cite{wu2020noisy,
zhu2018anisotropic,liu2021noisy,zhou2019toward, neelakantan2015adding} for its
applications.
The second drift term \(- \lambda \Theta^i_t\),
coming from our quadratic regularization,
is called \emph{weight decay}
in the field of machine learning.
It is believed to lead to better generalizations of the trained neural network
(see \cite{krogh1991simple,loshchilov2017decoupled}).
\end{rem}

\begin{rem}[Noised data]
\label{rem:noised-data}
In the previous remark we suppose the data available \(\{z_k,y_k\}_{k=1}^N\)
are precise: \(y_k = f(z_k)\),
while in practice they may be subject to errors:
\(y'_k = f(z_k) + \varepsilon_k\).
The new collection of points \(\{z_k, y'_k\}_{k=1}^N\) induces another mean
field functional \(F'_\textnormal{NNet}\) defined by
\[
F'_\textnormal{NNet} (m) = \frac{1}{2K} \sum_{k=1}^K
\bigl(y'_k - \Expect^{\Theta \sim m} [\Phi(\Theta; z_k)]\bigr)^2.
\]
From the triangle inequality for the \(L^2\)-distance we deduce
\[
|F'_\textnormal{NNet}(m) - F_\textnormal{NNet}(m)|
\leq \biggl( \frac 1K\sum_{k=1}^K \varepsilon_k^2 \biggr)^{\!1/2}
F_\textnormal{NNet}(m)^{1/2}
+ \frac{1}{2K} \sum_{k=1}^K \varepsilon_k^2.
\]
The actual \(N\)-neuron training process is therefore the noised gradient
descent for the functional \(F' \coloneqq F'_\textnormal{NNet} +
F_\textnormal{Ext}\)
and approximately converges to \((m'_\infty)^{\otimes N}\)
where \(m'_\infty\) minimizes \(\mathcal F' = F' + \frac{\sigma^2}{2}H\).
The difference between respective minima can be bounded as follows:
\begin{align*}
\mathcal F'(m'_\infty) - \mathcal F(m_\infty)
&\leq \mathcal F'(m_\infty) - \mathcal F(m_\infty)
= F'_\textnormal{NNet} (m_\infty) - F_\textnormal{NNet} (m_\infty) \\
&\leq \biggl( \frac 1K\sum_{k=1}^K \varepsilon_k^2 \biggr)^{\!1/2}
F_\textnormal{NNet}(m_\infty)^{1/2}
+ \frac {1}{2K}\sum_{k=1}^K \varepsilon_k^2.
\end{align*}
Hence the additional error converges to zero as the noise in the data
\((\varepsilon_k)_{k=1}^K\) tends to zero.
\end{rem}

\begin{rem}[Advantages over other approaches]
\label{rem:nn-convexity-advantage}
Our \cref{thm:ps-entropy-convergence,thm:poc} establish the exponential
convergence of the \(N\)-neurons training process \cref{eq:nn-ps-sde}
without supposing the truncation satisfies the regularity conditions such as
\(\Vert \nabla^k \ell \Vert_\infty < c\) for some small constant \(c\).
This stands in contrast to many previous studies on uniform-in-time propagation
of chaos
relying on the smallness of the mean field interaction
(e.g.\ \cite{durmus2020elementary} and the first setting
of \cite{delarue2021uniform}).
Yet the smallness approach does not apply to general neural networks:
in our setting, the smallness requires the Lipschitz constants \(M^F_{mm}\)
to be smaller than a constant times \(\rho\),
which we denote by \(M^F_{mm} \lesssim \rho\),
and the relation is difficult to verify.
Indeed, using the constants \(M^F_{mm}, \rho\) obtained in \cref{exm:nn-ngd},
we need
\[
\Bigl(\Vert \ell' \Vert_\infty^2
+ \Vert \ell \Vert_\infty^2 \Vert \varphi' \Vert_\infty^2
\bigl(1 + M_2(\mu)\bigr) \Bigr)^{\!1/2}
\lesssim \frac{\lambda}{2}
\exp\bigl( - 2 (\Vert f \Vert_{L^1(\mu)} + \Vert \ell \Vert_\infty) \Vert \ell
\Vert_\infty\bigr).
\]
This forces either the regularization \(\lambda\) to be very large or the
truncation \(\Vert \ell \Vert_\infty\) to be very small.
In conclusion, our approach based on the functional convexity
offers the advantage of obtaining the exponential convergence,
albeit at a very slow rate,
without such restrictions on \(\lambda\) or \(\ell\).
\end{rem}

\subsection{Numerical experiments}
\label{sec:exp}

As explained in \cref{exm:nn,exm:nn-ngd},
the MFL dynamics for training two-layer neural networks
verifies all the conditions of our theorems,
so its particle systems satisfy the uniform exponential energy dissipation
\cref{eq:ps-entropy-convergence}.
We now present our numerical experiments.

\paragraph{Setup.}

We aim to train a neural network to approximate the elementary function
\(z \mapsto f(z) = \sin 2\pi z_1 + \cos 2\pi z_2\)
on \([0,1]^2\).
We uniformly sample \(K\) points \(\{z_i\}_{k=1}^K\) from \([0,1]^2\)
and calculate the corresponding labels \(y_k = f(z_k)\)
to prepare our training data \(\{ z_k, y_k \}_{k=1}^K\).
These points are plotted in \cref{fig:samples}.
We fix the truncation function \(\ell\) by \(\ell(x) = (x \wedge 100) \vee
-100\)
and the sigmoid activation function \(\varphi\) by
\(\varphi(x) = 1/\bigl(1+\exp(-x)\bigr)\).
The Brownian noise has volatility \(\sigma\),
and it is necessary to apply the scaling transform in \cref{rem:vol-scaling}
before comparing to the theoretical results.
Additionally, the quadratic regularization constant \(\lambda\) is fixed
in our experiments.
The initial values \((\Theta^i_0)_{i=1}^N = (c^i_0, a^i_0, b^i_0)_{i=1}^N\) of
the \(N\) neurons are sampled independently
from a normal distribution \(m_0\) in four dimensions.
The training process \cref{eq:nn-ps-sde} is discretized with time step \(\Delta
t\)
and terminated at time \(T\).
The values of the hyperparameters $K$, $\sigma$, $m_0$, $\Delta t$, $T$
are listed in \cref{tab:hypara}
and the training algorithm is shown in \cref{alg:nn-ngd}.
We take the number of neurons \(N\) to be \(2^P\) for \(P = 6,\ldots,10\) and
repeat the training \(10\) times for each \(N\).

\begin{table}
\begin{minipage}{0.5\linewidth}
\centering
\includegraphics[width=\linewidth]{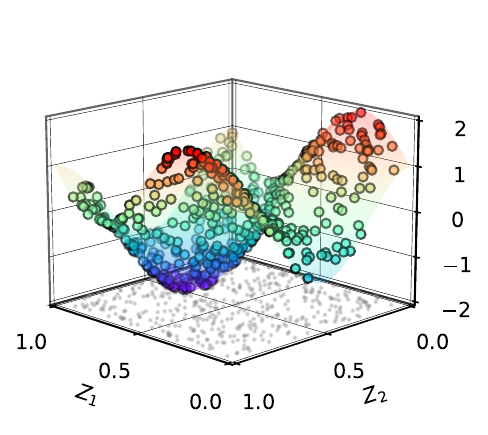}
\captionof{figure}{Data samples \(\{z_k, y_k\}_{k=1}^K\) (schematic).}
\label{fig:samples}
\end{minipage}
\hfill
\begin{minipage}{0.35\linewidth}
\centering
\begin{tabular}{cc}
\toprule
Parameters & Value \\
\midrule
$\Delta t$ & $0.2$ \\
$T$	       & $4000$ \\
$K$        & $1000$ \\
$m_0$      & $\mathcal N(0,5^2)$ \\
$\sigma$   & 1 \\
$\lambda$  & $10^{-5}$ \\
\bottomrule
\end{tabular}
\caption{Hyperparameters of neural network training.}
\label{tab:hypara}
\end{minipage}
\end{table}

\begin{algorithm}
\caption{Noised gradient descent for training a two-layer neural network}
\label{alg:nn-ngd}
\KwIn{number of particles \(N\),
activation \(\varphi\),
truncation \(\ell\),
data set \((z_k,y_k)_{k=1}^K\),
noise \(\sigma\),
initial distribution \(m_0\),
time step \(\Delta t\),
time horizon \(T\)}
\KwOut{$(\Theta^i_{T})_{i=1}^N$}
generate i.i.d.\ $\Theta_0^i = (A^i_0,B^i_0,C^i_0) \sim m_0, i = 1,\ldots,N$\;
\For{$t=0, \Delta t, 2\Delta t,\ldots,T-\Delta t$}{
generate i.i.d.\ $\mathcal{N}^i_{t} \sim \mathcal N(0,1), i=1,\ldots,N$\;
\tcp{update particles according to discretized Langevin}
\For{$i = 1,\ldots,N$}{
$\Theta^i_{t+\Delta t} \leftarrow \Theta^i_{t}
- \Bigl( D_m F_\textnormal{NNet}
\bigl(\frac{1}{N}\sum_{j=1}^N \delta_{\Theta^j_{t}},\Theta^i_{t}\bigr)
+ \lambda \Theta^i_{t} \Bigr) \Delta t
+ \sigma \sqrt{\Delta t}\,\mathcal{N}^i_{t}$\;
\tcc{where
\(D_m F_\textnormal{NNet} \bigl(\frac 1N \sum_{j=1}^N \delta_{\Theta^j_t},
\Theta^i_t\bigr)
= \frac{1}{K} \sum_{k=1}^K
\bigl(y_k - \Phi^N( \Theta^1_t, \ldots, \Theta^N_t; z_k) \bigr)
\frac{\partial \Phi}{\partial \theta} (\Theta^i_t; z_k)\)}}}
\end{algorithm}

\paragraph{Results.}

\begin{figure}
\begin{minipage}{0.45\linewidth}
\centering
\includegraphics[width=\linewidth]{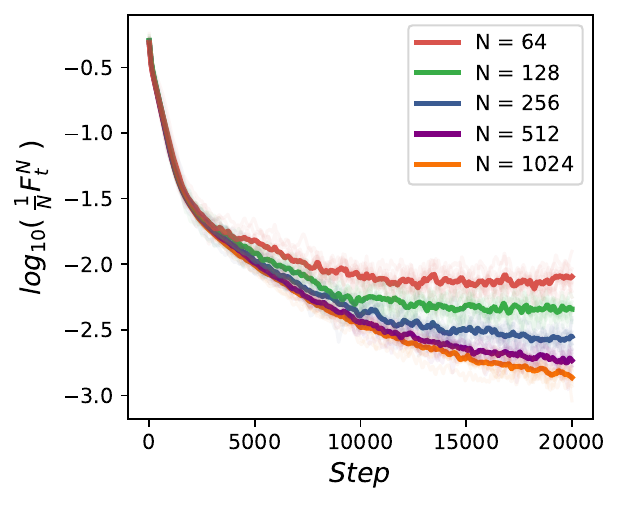}
\captionof{figure}{Individual (shadowed) and \(10\)-averaged (bold) losses
versus time steps.}
\label{fig:loss}
\end{minipage}
\hfill
\begin{minipage}{0.45\linewidth}
\centering
\includegraphics[width=\linewidth]{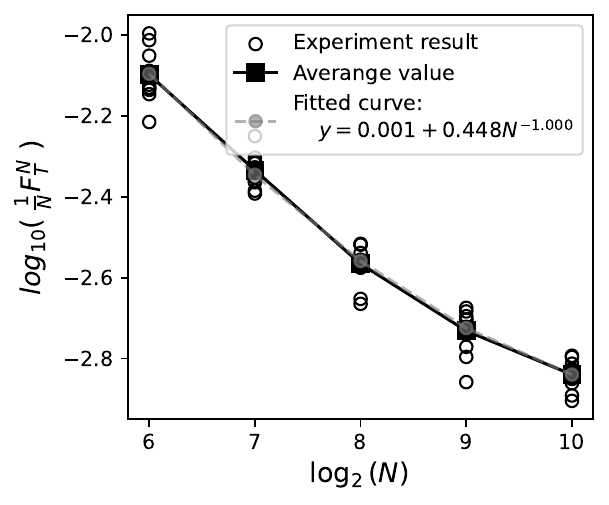}
\captionof{figure}{Average losses of last \(500\) steps for individual
trainings (shadowed) and its \(10\)-average (bold).}
\label{fig:error-gap}
\end{minipage}
\end{figure}

We compute the sum of the \(N^{-1}\)-scaled loss
\(\frac 1N F_\textnormal{NNet}^N (\Theta^1_t,\ldots,\Theta^N_t)\)
at each time \(t\)
and plot its evolution in \cref{fig:loss}.
We observe the value of \(\frac 1N F^N_\textnormal{NNet}\) first decreases
exponentially
and then decreases more slowly or even stabilizes.
To explore the relationship between this residual error and the number of
neurons,
for each value of \(N\) we calculate the average value of \(\frac
1NF^N_\textnormal{NNet}\)
during the last \(500\) training steps
and take the average of these values over the \(10\) independent runs.
The results are plotted in \cref{fig:error-gap}.

\paragraph{Discussions.}

Our truncation function \(\ell\) does not have bounded derivatives of up to
fourth order as required in \cref{exm:nn-ngd}
and we can work around this by taking a sequence of regular \(\ell_n\)
approximating \(\ell\) since the constants \(M^F_{mm}, \rho\) depends only on
\(\Vert \ell \Vert_\infty, \Vert \ell' \Vert_\infty\).
In our experiment we also ignore the time-discretization error and the
difference between training and validation data sets.
As shown in \cref{fig:loss} the losses first decrease exponentially
at a uniform rate for different numbers of neurons, \(N\).
This is consistent with the convergence rate \(\rho' - \frac{C_1}{N}\)
predicted by \cref{thm:ps-entropy-convergence,thm:poc}.
However, the LSI constant obtained in \cref{exm:nn-ngd} by Holley--Stroock is
excessively small
and fails to predict the actual convergence rate.
Given that the Holley--Stroock method relies solely on the boundedness of
neural networks,
this phenomenon suggests the internal structure of neural networks allows for a
faster convergence rate
that is not captured by the perturbation lemma.

We fit the residual losses with the curve \(\frac \alpha N + \beta\)
in \cref{fig:error-gap}.
We choose this parametrization for two reasons:
the first term \(\frac \alpha N\) corresponds to the error term in the
convergence result \cref{eq:ps-entropy-convergence} of the free energy \(\frac
1N \mathcal F^N(m^N_t)\);
the second term \(\beta\) accounts for the facts that
\(\mathcal F(m_\infty) \neq 0\) and that
the free energy differs from the neural network's loss by
\[
\frac 1N \mathcal F^N(m^N_t) - \frac 1NF^N_\textnormal{NNet}(m^N_t)
= \frac{\lambda}{2N} \int |\boldsymbol\theta|^2 m^N_t(\dd\boldsymbol\theta)
+ \frac{\sigma^2}{2N} H(m^N_t).
\]
In particular the relative entropy \(H(m^N_t)\) can not be directly calculated.

\section{Mean field system}
\label{sec:mf}

\subsection{Existence of the measures
\texorpdfstring{$\hat m$, $m_\infty$, $m^N_\infty$}{hat m, m∞, mN∞}}

Our assumptions differ from those in the earlier works,
such as \cite{HRSS19}.
Specifically, we do not require the coercivity condition of type
\[
\forall m\in \mathcal P_2(\mathbb R^d),~\forall x \in \mathbb R^d,\qquad
D_m F(m, x) \cdot x \geq C (|x|^2 - 1).
\]
Instead we only assume the condition \cref{eq:first-der} on \(D_m F(m,x)\).
As a result, the existence of the measures \(\hat m, m_\infty, m_\infty^N\),
introduced in \cref{sec:main-results}, is not obvious.
In this subsection we show that thanks to the conditions
\cref{eq:convex,eq:lip-in-m,eq:lsi},
these measures are indeed well defined.

First we sketch a proof that regular enough measures satisfying an LSI in
\(\mathbb R^d\) have finite moments.

\begin{lem}
\label{lem:lsi-implies-moments}
Let \(\mu(\dd x) = e^{-\Psi}\dd x\)
be a probability measure in \(\mathbb R^d\) where
\(\Psi\) is twice differentiable with the bound \(|\nabla^2\Psi| \leq C\).
If \(\mu\) satisfies an LSI, i.e.\ \cref{eq:lsi} holds
when \(\hat m\) is replaced by \(\mu\) for some \(\rho > 0\),
then \(\mu \in \cap_{p \geq 1} \mathcal P_p (\mathbb R^d)\) and \(\int
e^{\alpha|x|} \mu(\dd x) < +\infty\) for all \(\alpha \geq 0\).
\end{lem}

\begin{proof}
Here we repeat the argument of Otto and Villani in \cite{ov}.
Suppose \(\mu\) satisfies a \(\rho\)-LSI (but we do not suppose \(\mu \in
\mathcal P_2(\mathbb R^d)\) a priori).
For every measure \(\nu \in \mathcal P_2(\mathbb R^d)\)
of finite entropy (e.g.\ the Gaussians),
the heat flow
\[\partial_t \nu_t = \Delta \nu_t + \nabla \cdot (\nu_t \nabla \Psi),\qquad
\nu_0 = \nu\]
is well defined and is an absolutely continuous curve in \((\mathcal P_2,
W_2)\) thanks to
the bound \(|\nabla^2 \Psi| \leq C\) and \cite[Theorem 7.4.1]{fpkeq}.
Hence by the argument of \cite[Proposition 1']{ov},
we can obtain \(H(\nu_t | \mu) \leq H(\nu | \mu) e^{-4\rho t}\) and
\begin{equation}
\label{eq:ov}
W_2 (\nu, \nu_t) \leq \frac 1{\sqrt{\rho}}
\Bigl(\sqrt{H(\nu | \mu)} - \sqrt{H(\nu_t | \mu)}\Bigr).
\end{equation}
The sequence \(\nu_t\) are tight in the weak topology of \(\mathcal P\) since
we have \(\rho W_2(\nu, \nu_t)^2 \leq H(\nu | \mu) = \int (\log \nu + \Psi) \nu
< +\infty\) (recall that \(\Psi\) is of quadratic growth).
By the lower-semicontinuity of \(H(\cdot | \mu)\)
we must have \(\nu_t \to \mu\) in \(\mathcal P\) weakly when \(t \to \infty\).
Then we take \(\liminf_{t \to \infty}\) on both side of \cref{eq:ov}
and use the lower-semicontinuity of \(W_2\) with respect to the weak topology
of \(\mathcal P\)
to obtain Talagrand's inequality
\[
\rho W_2^2(\nu,\mu) \leq H(\nu | \mu).
\]
Hence \(\mu \in \mathcal P_2\).
Finiteness of higher moments and exponential moments then follows from
concentration inequalities via Herbst's argument
(see e.g.\ the proof of \cite[Theorem 5.5]{concentration}).
\end{proof}

We give a sufficient condition to the existence of \(\hat m\) for every \(m \in
\mathcal P_2(\mathbb R^d)\) so that the condition \cref{eq:lsi} makes sense.

\begin{prop}
\label{prop:exist-m-hat}
Assume \(F\) satisfies \cref{eq:lip-in-m}.
If there exists a measure \(m_0 \in \mathcal P_2(\mathbb R^d)\) such that
\(\hat m_0\) is well defined
(i.e.\ \(Z(\hat m_0) < +\infty\))
and \(m_0\) satisfies LSI \cref{eq:lsi},
then \(\hat m\) are well defined
(i.e.\ \(Z(\hat m) < +\infty\)) for all \(m \in \mathcal P_2(\mathbb R^d)\).
\end{prop}

\begin{proof}
By definition we have
\begin{multline*}
Z(\hat m) = \int \exp \biggl(- \frac{\delta F}{\delta m}(m,x) \biggr)\dd x \\
= Z(\hat m_0) \int \exp \biggl( \frac{\delta F}{\delta m}(m_0,x)
-\frac{\delta F}{\delta m}(m,x) \biggr) \hat m_0(\dd x),
\end{multline*}
where the term on the exponential is of linear growth since its derivative
is uniformly bounded:
$\bigl\lvert\nabla\bigl(\frac{\delta F}{\delta m}(m_0,x)
- \frac{\delta F}{\delta m}(m,x)\bigr)\bigr\rvert
= |D_m F(m_0,x) - D_m F(m,x)| \leq M^F_{mm} W_2(m_0,m)$.
But by \cref{lem:lsi-implies-moments}, all exponential moments of \(\hat m_0\)
are finite.
Thus \(Z(\hat m) < +\infty\) and \(\hat m\) is well defined.
\end{proof}

We now show that the \(N\)-particle invariant measure is also well defined.

\begin{prop}
\label{prop:ps-invariant-measure}
Assume \(F\) satisfies \cref{eq:convex,eq:lsi}.
Then the measure \(m^N_\infty\) in \cref{eq:def-ps-invariant-measure} is well
defined and has finite exponential moments for all \(N \geq 2\).
\end{prop}

\begin{proof}
Fix \(m_0 \in \mathcal P_2(\mathbb R^d)\).
Using convexity we obtain
\begin{multline*}
NF(\mu_{\mathbf x}) \geq NF(m_0) + N \int \frac{\delta F}{\delta m}(m_0,y)
(\mu_{\mathbf x} - m_0)(\dd y) \\
= NF(m_0) - N \int \frac{\delta F}{\delta m}(m_0,y) m_0(\dd y) + \sum_{i=1}^N
\frac{\delta F}{\delta m} (m_0,x^i).
\end{multline*}
The integral \(\int \frac{\delta F}{\delta m}(m_0,y)m_0(\dd y)\) is finite
thanks to \cref{lem:lsi-implies-moments}.
Hence
\[
\int \exp\bigl(-NF(\mu_{\mathbf x})\bigr) \dd\mathbf x
\leq C \int \exp\biggl(- \sum_{i=1}^N \frac{\delta F}{\delta m} (m_0,x^i)
\biggr)\dd\mathbf x
= C \bigl(Z(\hat m_0)\bigr)^N < +\infty.
\]
Apply the same argument to
\(\int \exp\bigl(\alpha\sum_{i=1}^N |x^i|\bigr)
\exp\bigl(-NF(\mu_\mathbf{x})\bigr)\dd\mathbf x\) we obtain the finiteness of
exponential moments.
\end{proof}

\begin{prop}
\label{prop:mf-invariant-measure}
Assume \(F\) satisfies \cref{eq:convex,eq:lsi,eq:lip-in-m,eq:first-der}.
Then the mean field free energy \(\mathcal F\), defined in
\cref{eq:def-mf-free-energy},
has a unique minimizer \(m_\infty\).
The minimizer \(m_\infty\) is also the unique solution to the first-order
equation \cref{eq:mf-foc}
and an invariant measure to the MFL dynamics \cref{eq:mf-fp}.
\end{prop}

\begin{proof}
Recall that \(\mathcal F(m) = F(m) + H(m)\) where
the absolute entropy \(H(m)\) is well defined for \(m \in \mathcal P_2\) and
has value in \(( - \infty, +\infty]\) thanks to the decomposition
\begin{multline}
\label{eq:entropy-decomposition}
H(m) = \int \log m(x) m(x) \dd x \\
= \int \log \frac{m(x)}{(2\pi)^{-d/2} e^{-x^2/2}} m(x)\dd x
+ \int \biggl(\log (2\pi)^{-d/2} - \frac{x^2}{2}\biggr) m(x) \dd x.
\end{multline}
The first term, which is the relative entropy between \(m\) and a normalized
Gaussian, is always nonnegative
and the second term is finite.
Moreover the free energy \(\mathcal F\) satisfies
\begin{multline}
\mathcal F(m) - F(m_0)
\geq \int \frac{\delta F}{\delta m} (m_0,x) (m - m_0) (\dd x) + H(m) \\
= - \int \log \hat m_0(x) (m - m_0)(\dd x) + H(m)
= H(m | \hat m_0) + \int \log \hat m_0(x) m_0(\dd x)
\label{eq:free-energy-lower-bound}
\end{multline}
for all \(m, m_0 \in \mathcal P_2\) such that \(m_0\) has finite entropy.
Since the LSI \cref{eq:lsi} implies the \(T_2\) inequality \cref{eq:t2},
the functional \(\mathcal F\) has \(\mathcal P_2\)-coercivity:
\[
\rho W_2^2(m, \hat m_0) \leq H(m | \hat m_0) \leq \mathcal F(m) - \int \log
\hat m_0(x) m_0(\dd x) - F(m_0).
\]
The conditions \cref{eq:convex,eq:first-der} imply also the \(\mathcal
P_2\)-lower-continuity of \(F\):
if \((m_n)_{n \in \mathbb N}\) is a sequence convergent to \(m\) in the weak
topology of \(\mathcal P_2\),
then we have
\begin{align*}
&\hspace{-1em}\liminf_{n} F(m_n) - F(m) \\
&\geq \liminf_{n} \int \frac{\delta F}{\delta m}(m,x) (m_n - m)(\dd x) \\
&= \liminf_{n} \int \biggl(\frac{\delta F}{\delta m}(m,x)
- \frac{\delta F}{\delta m}(m,0)\biggr) (m_n - m)(\dd x) \\
&\geq \liminf_{n} \int \biggl( D_m F(m,0) \cdot x
- \frac {M^F_{mx}}{2}|x|^2 \biggr) (m_n - m)(\dd x) \\
&= 0.
\end{align*}
Here the second inequality follows from Taylor's formula and \(M^F_{mx}\)
denotes the constant in the condition \cref{eq:first-der}.
The entropy \(H\) is also \(\mathcal P_2\)-lower-semicontinuous by the previous
decomposition \cref{eq:entropy-decomposition}.
The free energy \(\mathcal F\) is then lower-bounded, coercive,
lower-semicontinuous and convex, so there exists unique minimizer in \(\mathcal
P_2\) which we denote by \(m_\infty\).

Now we show the equivalence between the minimizing property of the free energy
\(\mathcal F\) and the first-order condition \cref{eq:mf-foc}.
If \(m_0\) satisfies \cref{eq:mf-foc} then \(\hat m_0 = m_0\) and from
\cref{eq:free-energy-lower-bound} we deduce \(\mathcal F(m) \geq \mathcal
F(m_0)\) for all \(m \in \mathcal P_2\),
i.e.\ \(m_0\) is the minimizer of \(\mathcal F\).
For the reverse implication we refer readers to the necessary part of the proof
of \cite[Proposition 2.5]{HRSS19}.

Finally since \(m_\infty\) satisfies \cref{eq:mf-foc} we have
\[
\Delta m_\infty + \nabla \cdot (D_m F(m_\infty,x) m_\infty)
= \nabla \cdot \Biggl(m_\infty
\nabla \biggl(\frac{\delta F}{\delta m} (m_\infty,x) + \log m_\infty\biggr)
\Biggr) = 0,
\]
and \(m_\infty\) is invariant to \cref{eq:mf-fp}.
\end{proof}

\begin{rem}
\label{rem:converse-invariant-measure}
We will establish the uniqueness of the invariant measure of the MFL
in \cref{cor:mf-invariant-measure-unique}
after deriving the free energy dissipation formula \cref{eq:mf-energy-decrease}.
\end{rem}

\subsection{Proof of \cref{thm:mf-entropy-convergence}}
\label{sec:proof-thm:mf-entropy-convergence}

First we recall the definition of \(AC^2\) curves in \cite{gf}.

\begin{defn}
\label{defn:ac2}
Let \((X,d)\) be a complete metric space
and \(x : [a,b] \to X\) be a continuous mapping.
We say \(x\) is \emph{absolutely continuous} (a.c.)
and write \(x \in AC\bigl([a,b]; (X,d)\bigr)\) if there exists \(m \in
L^1([a,b])\) such that
\[
\forall a \leq s < t \leq b,\qquad d\bigl(x(s), x(t)\bigr) \leq \int_s^t m(u)
du.
\]
We say \(x \in AC^2\bigl([a,b]; (X,d)\bigr)\) if additionally
\(m \in L^2([a,b])\).
For a globally defined curve \(x : [t_0, +\infty) \to X\)
we say \(x\) belongs to the class \(AC_\textnormal{loc}^2\)
and denote \(x \in AC_\textnormal{loc}^2\bigl([t_0, +\infty; (X, d)\bigr)\),
if $x \in AC_\textnormal{loc}^2 \bigl( [t_0, T]; (X,d) \bigr)$
for every $T \geq t_0$.
\end{defn}

Now we state the wellposedness and regularity result.

\begin{prop}[Existence, uniqueness and regularity of MFL]
\label{prop:mf-exist-unique-regular}
Assume \(F\) satisfies \cref{eq:lip-in-m,eq:first-der}.
Then
\begin{enumerate}
\item for all \(m_0 \in \mathcal P_2(\mathbb R^d)\)
there exists a unique continuous flow \(m : [0,+\infty) \to \mathcal
P_2(\mathbb R^d)\)
solving weakly the Fokker--Planck equation \cref{eq:mf-fp};
\item moreover, this solution has density and finite entropy for positive time:
\[\forall t > 0,\qquad\int\lvert\log m_t(x)\rvert m_t(x)\dd x < +\infty;\]
\item if additionally \(m_{t_0}\) has finite entropy for some \(t_0 \geq 0\),
then the integral
\begin{equation}
\int_{t_0}^t \int \frac{|\nabla m_s(x)|^2}{m_s(x)}\dd x\dd s
\label{eq:integrated-fisher}
\end{equation}
is finite for every \(t \geq t_0\);
therefore \((m_t)_{t \geq t_0} \in AC^2_\textnormal{loc}
\bigl([t_0,+\infty); (\mathcal P_2,W_2)\bigr)\)
and has tangent vector
\(v_t(x) = - D_m F(m_t,x) - \nabla \log m_t(x)\) for \(t \geq t_0\)
a.e.\ in the sense of \cite[Proposition 8.4.5]{gf}.
\end{enumerate}
\end{prop}

Due to the technical nature of this proposition its proof is postponed to
\cref{sec:mf-technical}.
Using the results of \cref{prop:mf-exist-unique-regular} and applying the
formalism of \cite{gf},
we establish the free energy dissipation formula,
which is crucial to our studies on the dynamics of gradient flow.

\begin{prop}[Energy dissipation]
\label{prop:mf-energy-decrease}
Assume \(F\) satisfies \cref{eq:lip-in-m,eq:first-der}.
If \(m_{t_0}\) is a measure of finite entropy and finite second moment for some
\(t_0 \geq 0\),
then the free energy \(\mathcal F\), defined in \cref{eq:def-mf-free-energy},
is absolutely continuous along the flow \((m_t)_{t \geq t_0}\) constructed in
\cref{prop:mf-exist-unique-regular}.
Moreover it has derivative
\begin{equation}
\label{eq:mf-energy-decrease}
\frac{\dd\mathcal F(m_t)}{\dd t}
= - \int | D_m F(m_t, x) + \nabla \log m_t(x) |^2
m_t(\dd x),\qquad \textnormal{for $t \geq t_0$ a.e.}
\end{equation}
\end{prop}

\begin{proof}
We will apply the chain rule result of \cite[Proposition 10.3.18]{gf}
and we verify its conditions,
namely, the differentiability of the free energy \(\mathcal F = F + H\)
and of the flow of measures \(m_t\).
Firstly under the conditions \cref{eq:lip-in-m,eq:first-der} we can apply the
argument of \cite[Lemma A.2]{chizat2022mean} to show that \(F : \mathcal
P_2(\mathbb R^d) \to \mathbb R\) is \(-\lambda\)-geodesically-convex for some
\(\lambda > 0\)
and it has differential \(D_m F (m_t,\cdot)\) at \(m_t\).
Secondly the entropy \(H : \mathcal P_2(\mathbb R^d) \to (-\infty, +\infty]\)
is also \(0\)-geodesically-convex by the result of \cite[Proposition 9.3.9]{gf}
and for \(t \geq t_0 \) a.e.\ has subdifferential \(\nabla \log m_t\) at
\(m_t\) by \cite[Theorem 10.4.6]{gf},
thanks to the regularity bounds in the previous
\cref{prop:mf-exist-unique-regular}.
Hence the free energy \(\mathcal F = F + H\) is
\(-\lambda\)-geodesically-convex and has differential \(D_m F(m_t, \cdot) +
\nabla \log m_t\) at \(m_t\).
For the flow of measures \(m_t\)
we have already obtained its \(AC^2\)-regularity in the previous proposition
and its tangent vector reads \(v_t = - D_m F(m_t, \cdot) - \nabla \log m_t\) at
\(m_t\) for \(t \geq t_0\) a.e.
Then we can apply the chain rule to obtain the absolute continuity of \(t
\mapsto \mathcal F(m_t)\) and
\[
\forall T > t_0,\quad
\mathcal F(m_T) - \mathcal F(m_{t_0})
= \int_{t_0}^T\bigl(D_m F(m_t, x) + \nabla \log m_t(x)\bigr)
\cdot v_t(x) m_t(\dd x)\dd t
\]
which is the desired result.
\end{proof}

\begin{cor}[Uniqueness of the invariant measure]
\label{cor:mf-invariant-measure-unique}
Under \cref{eq:convex,eq:lip-in-m,eq:lsi,eq:first-der} there exists a unique
invariant measure in \(\mathcal P_2(\mathbb R^d)\) to the mean field dynamics
\cref{eq:mf-fp}.
\end{cor}

\begin{proof}
The existence part is already shown in \cref{prop:mf-invariant-measure}.
Let \(m_* \in \mathcal P_2(\mathbb R^d)\) be an invariant measure.
We let the initial condition \(m_0\) be equal to \(m_*\)
and construct according to \cref{prop:mf-exist-unique-regular}
the MFL solution \((m_t)_{t \geq 0}\).
By the invariance of \(m_*\) we have \(m_t = m_*\) for all \(t \geq 0\),
so \(m_*\) must have density and finite entropy.
We then apply the energy dissipation formula \cref{eq:mf-energy-decrease} and
obtain
\[
\text{for $x \in \mathbb R^d$ a.e.,}\qquad
D_m F(m_*, x) + \nabla \log m_*(x) = 0.
\]
Integrating this equation, we obtain \(m_*\) solves the first-order condition
\cref{eq:mf-foc}
which has unique solution by \cref{prop:mf-invariant-measure}.
\end{proof}

Now we show the close relation between the free energy and the relative
entropies.

\begin{lem}[Entropy sandwich]
\label{lem:mf-entropy-sandwich}
Assume \(F\) satisfies \cref{eq:convex,eq:lip-in-m,eq:lsi,eq:first-der}.
Then for every \(m \in \mathcal P_2(\mathbb R^d)\)
we have
\begin{multline}
\label{eq:entropy-sandwich}
H(m | m_\infty) \leq \mathcal F(m) - \mathcal F(m_\infty)
\leq H(m | \hat m) \\
\leq \biggl(1 + \frac{M^F_{mm}}{\rho} + \frac{(M^F_{mm})^2}{2\rho^2}
\biggr) H(m | m_\infty).
\end{multline}
\end{lem}

\begin{proof}
The first two inequalities are proved in \cite[Lemma 3.4]{chizat2022mean}.
We show the rightmost one.
Recall that \(Z(\hat m)\) is the normalization constant defined in
\cref{eq:def-Z-m-hat}.
We have
\begin{multline*}
H(m | \hat m) - H(m | m_\infty)
= \int \biggl( \log \frac{m}{\hat m} - \log \frac{m}{m_\infty} \biggr) m
= \int \log \frac{m_\infty}{\hat m} m \\
= \int \biggl( \frac{\delta F}{\delta m}(m, x)
- \frac{\delta F}{\delta m}(m_\infty, x) \biggr) m(x) \dd x
+ \log Z(\hat m) - \log Z(m_\infty).
\end{multline*}
By Jensen's inequality, the difference between
\(\delta \coloneqq \log Z(\hat m) - \log Z(\hat m_\infty)\)
satisfies
\begin{align*}
\delta
&= \log Z(\hat m)
- \log \int \exp\biggl(- \frac{\delta F}{\delta m}(m_\infty,x)\biggr)\dd x \\
&= \log Z(\hat m)
- \log \int \exp\biggl(- \frac{\delta F}{\delta m}(m_\infty,x) - \log \hat
m(x)\biggr) \hat m(x)\dd x \\
&\leq \log Z(\hat m)
+ \int \biggl(\frac{\delta F}{\delta m}(m_\infty,x) + \log \hat m(x)\biggr)
\hat m(x)\dd x \\
&\leq \log Z(\hat m)
+ \int \biggl(\frac{\delta F}{\delta m}(m_\infty,x) - \frac{\delta F}{\delta
m}(m,x) - \log Z(\hat m)\biggr) \hat m(x)\dd x \\
&= \int \biggl(\frac{\delta F}{\delta m}(m_\infty,x) - \frac{\delta F}{\delta
m}(m,x)\biggr) \hat m(x)\dd x.
\end{align*}
Then we have by Kantorovich duality and \(W_1\)-Lipschitzianity in
\cref{eq:lip-in-m}
\begin{align*}
H(m|\hat m) - H(m|m_\infty)
&\leq \int \biggl(\frac{\delta F}{\delta m}(m,x)
- \frac{\delta F}{\delta m}(m_\infty,x)\biggr)\bigl(m(x) - \hat m(x)\bigr)\dd x
\\
&\leq \Vert D_m F(m,x) - D_m F(m_\infty,x) \Vert_\infty W_1(m, \hat m) \\
&\leq M^F_{mm} W_1(m,m_\infty) W_1(m,\hat m) \\
&\leq M^F_{mm} W_1(m,m_\infty) \bigl(W_1(m,m_\infty) + W_1(\hat m,m_\infty)\bigr).
\end{align*}
Note that, for the first term in the bracket above, we have
\(W_1(m,m_\infty) \leq W_2(m,m_\infty) \leq \sqrt{\rho^{-1} H(m | m_\infty)}\)
by the \(T_2\) and log-Sobolev inequalities, \cref{eq:t2,eq:lsi},
and for the second term, we have
\begin{align*}
W_1^2(\hat m,m_\infty)
&\leq W_2^2(\hat m,m_\infty)
\leq \frac 1{\rho} H(\hat m|m_\infty)
\leq \frac 1{4\rho^2} \int \biggl\lvert\nabla \log \frac{\hat m}{m_\infty}
\biggr\rvert^2 \hat m \\
&= \frac 1{4\rho^2} \int \lvert D_mF(m,x) - D_mF(m_\infty,x)\rvert^2
\hat m(x)\dd x
\\
&\leq \frac {(M^F_{mm})^2}{4\rho^2} W_1^2(m,m_\infty)
\leq \frac {(M^F_{mm})^2}{4\rho^3} H(m | m_\infty),
\end{align*}
which concludes.
\end{proof}

The proof of \cref{thm:mf-entropy-convergence}
is nothing but the combination of the previous two results.

\begin{proof}[Proof of \cref{thm:mf-entropy-convergence}]
By \cref{prop:mf-energy-decrease} we have
\begin{align*}
\frac{d\mathcal F(m_t)}{dt} &= - \int |D_m F(m_t,x) + \nabla \log m_t(x)|^2
m_t(\dd x)
= - I(m_t | \hat m_t) \\
&\leq - 4\rho H(m_t | \hat m_t)
\leq - 4\rho\bigl(\mathcal F(m_t) - \mathcal F(m_\infty)\bigr),
\qquad\text{for $t \geq t_0$ a.e.}
\end{align*}
The first inequality is due to the uniform log-Sobolev inequality \cref{eq:lsi}
and the second to the entropy sandwich \cref{eq:entropy-sandwich}.
The second inequality in \cref{eq:mf-entropy-convergence}
is then obtained by Grönwall's lemma,
and the first inequality has already been proved
in \cref{lem:mf-entropy-sandwich}.
\end{proof}

\subsection{\texorpdfstring{\(L^2\)}{L2}-convergence and hypercontractivity}
\label{sec:l2-convergence-hypercontractivity}

\subsubsection{Standard algebra}
\label{sec:std-alg}

We first work on dense set of sufficiently regular functions
that will be necessary our proofs.

For notational simplicity, define
\(b_t(x) \coloneqq - D_m F(m_t,x)\),
\(b_\infty(x) \coloneqq - D_m F(m_\infty,x)\)
and recall that \(h_t (x) \coloneqq \frac{dm_t}{dm_\infty} (x)\).
The relative density \(h_t\) then solves
\begin{equation}
\label{eq:h}
\partial_t h = \Delta h + (2 b_\infty - b_t) \cdot \nabla h
- \bigl(\nabla \cdot (b_t - b_\infty) + (b_t - b_\infty) \cdot b_\infty\bigr) h.
\end{equation}
In this subsection we will fix the flow of measures \(m_t\)
to be that constructed in \cref{prop:mf-exist-unique-regular}
and let \(h\) change independently from \(m_t\).
We will also only consider solutions
in \(L^\infty\bigl([t_0, T]; L^1(m_\infty)\bigr)\)
with initial value \(h_{t_0} \in L^1(m_\infty)\)
to the evolution equation \cref{eq:h}
(in the sense of \cite[(6.1.3)]{fpkeq}).
We then know that the solution is then unique
by applying \cite[Theorem 9.6.3]{fpkeq} to \(h m_\infty\).

\begin{defn}[Standard algebra]
\label{defn:std-alg}
The \emph{standard algebra} \(\mathcal A_+\)
is the set of positive and \(C^2\) functions \(h : \mathbb R^d \to (0,\infty)\)
satisfying the following conditions:
\begin{itemize}
\item there exists a constant \(M > 0\) such that
for every \(x \in \mathbb R^d\), \(\lvert\log h(x)\rvert \leq M(1+|x|)\);
\item for $k = 1$, $2$, there exist constants \(M_k > 0\) such that
for every \(x \in \mathbb R^d\),
\(|\nabla^k h(x)| \leq \exp\bigl(M_k(1+|x|)\bigr)\).
\end{itemize}
For a collection of functions \((h_i)_{i \in I}\) we say that
\(h_i \in \mathcal A_+\) uniformly for \(i \in I\)
or \((h_i)_{i \in I} \subset \mathcal A_+\) uniformly,
if there exist constants $M$, $M_1$, $M_2$ such that
the previous bounds holds for all $h_i$, $i \in I$.
\end{defn}

\begin{rem}
\label{rem:std-alg}
The word ``standard algebra'' is the terminology in \cite{lsi}.
Readers may have noticed \(\mathcal A_+\) is not an algebra in the usual sense,
as it contains only positive functions and is not closed under scalar
multiplication by \(-1\).
To remedy this we can define \(\mathcal A = \mathcal A_+ - \mathcal A_+\)
and \(\mathcal A\) is truely an algebra.
We introduce this unusual set of functions in order to do \(L^p\)-computations
for \(p < 1\).
\end{rem}

Then we can state the density and stability of \(\mathcal A_+\).

\begin{prop}[Density of \(\mathcal A_+\)]
\label{prop:std-alg-density}
Let \(p \geq 1\), \(q < 1\),
\(h : \mathbb R^d \to [0,+\infty]\) be a measurable function
and \(\mu\) be a probability measure on \(\mathbb R^d\) having a density
with respect to the Lebesgue measure.
If \(h \in L^p(\mu)\), then there exists a sequence \((h_n)_{n \in \mathbb N}\)
in \(\mathcal A_+\) such that \(h_n \to h\) in \(L^p(\mu)\);
if \(h \in L^q(\mu)\), then there exists a sequence \((h_n)_{n \in \mathbb N}\)
in \(\mathcal A_+\) such that \(\Vert h_n \Vert_{q} \to \Vert h \Vert_q\);
and if \(h \in L^p \cap L^q (\mu)\), then the sequence in \(\mathcal A_+\) can
be chosen such that both convergences hold.
\end{prop}

\begin{prop}[Stability of \(\mathcal A_+\) under flow]
\label{prop:std-alg-stability}
Assume \(F\) satisfies
\cref{eq:convex,eq:lip-in-m,eq:lsi,eq:first-der,eq:higher-der}.
For every \(t_0 \geq 0\) and \(h' \in \mathcal A_+\),
there exists a solution \(h : [t_0,+\infty) \to \mathcal A_+\) to \cref{eq:h}
with initial value \(h(t_0,\cdot) = h'\).
Moreover the temporal weak derivative \(\partial_t h\) exists
and \(h_t\) belongs to \(\mathcal A_+\) locally uniformly,
i.e., \((h_t)_{t \in K} \subset \mathcal A_+\) uniformly
for every compact subset \(K \subset [t_0, +\infty)\).
\end{prop}

The proofs of \cref{prop:std-alg-density,prop:std-alg-stability}
are postponed to \cref{sec:mf-technical} due to their technical nature.

\subsubsection{Proof of \cref{prop:l2-convergence}}

First, by working in \(\mathcal A_+\),
we obtain the following \(L^p\)-norm growth formula.

\begin{prop}[\(L^p\)-norm growth]
\label{prop:lp-growth}
Assume \(F\) satisfies
\cref{eq:convex,eq:lip-in-m,eq:lsi,eq:first-der,eq:higher-der}.
Let \(p \neq 0\)
and \(h : [a,b] \to \mathcal A_+\) be a solution to the evolution \cref{eq:h}.
Then the growth of \(p\)-norm \(t \mapsto \int h_t(x)^p m_\infty(\dd x)\) is
absolutely continuous and has derivative
\begin{multline}
\label{eq:lp-growth}
\frac{\dd}{\dd t} \int h_t(x)^p m_\infty(\dd x)
= p(p-1)\biggl(-\int h_t(x)^{p-2} |\nabla h_t(x)|^2 m_\infty(\dd x) \\
+ \int h_t(x)^{p-1} \nabla h_t(x) \cdot \bigl(b_t(x) - b_\infty(x)\bigr)
m_\infty(\dd x)\biggr)
\end{multline}
for \(t \in [a,b]\) a.e.
\end{prop}

\begin{proof}
We first suppose \(t \mapsto h(t,x)\) is \(C^1\) instead of only absolutely
continuous.
Notice that the evolution equation \cref{eq:h} of \(h\) can be rewritten as
\[
\partial_t h = (\Delta + b_\infty \cdot \nabla) h
- (b_t - b_\infty) \cdot \nabla h
- \frac{\nabla \cdot\bigl(m_\infty (b_t - b_\infty)\bigr)}{m_\infty} h,
\]
where the first term corresponds to the symmetric operator
\(\Delta + b_\infty \cdot \nabla\) in \(L^2(m_\infty)\).
We then have
\begin{align*}
&\hspace{-1em}\frac{\dd}{\dd t} \int h_t(x)^p m_\infty(\dd x) \\
&= p \int h_t(x)^{p-1} \bigl(\Delta + b_\infty(x) \cdot \nabla\bigr)
h_t(x) m_\infty(\dd x) \\
&\phantom{={}}\quad
- p \int h_t(x)^{p-1} \bigl(b_t(x) - b_\infty(x)\bigr)
\cdot \nabla h_t(x) m_\infty(\dd x) \\
&\phantom{={}}\quad
- p\int \nabla \cdot \bigl(m_\infty (b_t - b_\infty)\bigr)(x) h_t(x)^p\dd x \\
&= - p(p-1) \int h_t(x)^{p-2} |\nabla h_t(x)|^2 m_\infty(\dd x) \\
&\phantom{={}}\quad
- p \int h_t(x)^{p-1} \bigl(b_t(x) - b_\infty(x)\bigr) \cdot \nabla h_t(x)
m_\infty(\dd x) \\
&\phantom{={}}\quad
+ p \int \nabla h_t(x)^p \cdot \bigl(b_t(x) - b_\infty(x)\bigr) m_\infty(\dd x) \\
&= p(p-1) \biggl(
\begin{aligned}[t]
&- \int h_t(x)^{p-2} |\nabla h_t(x)|^2 m_\infty(\dd x) \\
&\quad+ \int h_t(x)^{p-1} \nabla h_t(x) \cdot
\bigl(b_t(x) - b_\infty(x)\bigr) m_\infty(\dd x) \biggr).
\end{aligned}
\end{align*}
We can justify the first equality by the dominated convergence theorem and
the two integrations by parts in the second one by an approximating sequence of
functions,
thanks to the fact that \(h_t \in \mathcal A_+\) locally uniformly.

Then, for the general case where $t \mapsto h_t(x)$ is only absolutely
continuous, thanks to the fact that $h_t$ belongs to $\mathcal A_+$
locally uniformly, we have for every $s$, $t \in [a,b]$ with $s \leq t$,
\[
\int h_t(x)^p m_\infty(\dd x) - \int h_s(x)^p m_\infty(\dd x)
= p \int_s^t \int h_u(x)^{p-1} \partial_u h_u(x) m_\infty(\dd x)\dd u,
\]
where $\partial_u h_u(x)$ is the weak derivative that exists only a.e.
Then we plug in the evolution equation \cref{eq:h} and compute as before.
\end{proof}

\begin{rem}
\label{rem:entropy-growth}
By dividing \cref{eq:lp-growth} by \(p-1\) and taking the limit \(p \to 1\),
one formally obtains
\begin{multline}
\label{eq:entropy-growth}
\frac{\dd}{\dd t} \int h_t(x)\log h_t(x) m_\infty(\dd x)
= - \int \frac{|\nabla h_t(x)|^2}{h_t(x)} m_\infty(\dd x) \\
+ \int \nabla h_t(x) \cdot \bigl(b_t(x) - b_\infty(x)\bigr) m_\infty(\dd x).
\end{multline}
This entropy growth formula is one of the key ingredients
of the method of Jabin and Wang \cite{JabinWang}
and has also been used in \cite{guillin2021uniformvortex}.
A weak version of this formula under weak regularity of $b$
has been rigorously proved in the Appendix A
of the first arXiv version of \cite{lacker2023sharp}.
In our case, the formula can be first rigorously proved for $h$ taking value
in $\mathcal A_+$, as is done in the proposition above,
and then we treat the general case by the density of $\mathcal A_+$.
\end{rem}

The \(L^p\)-norm growth formula implies the existence of a strongly continuous
semigroup in \(L^p(m_\infty)\) for all \(p \in [1, +\infty)\).

\begin{cor}[\(L^p\)-continuity of flow]
\label{cor:lp-cont}
Under the hypotheses of \cref{prop:lp-growth},
for every \(p \geq 1\) and every \(a \leq s \leq t \leq b\)
there exists a constant \(C_{s,t,p} > 0\) such that
\[
\int h_t(x)^p m_\infty(\dd x) \leq C_{s,t,p} \int h_s(x)^p m_\infty(\dd x)
\]
holds for every solutions to \cref{eq:h} in \(\mathcal A_+\).
Therefore the evolution equation \cref{eq:h} determines a strongly continuous
(and positive) semigroup \((P_s^t)_{s\leq t}\) in \(L^p_+(m_\infty)\) for \(p
\in [1,+\infty)\).
\end{cor}

\begin{proof}
For \(h_s \in \mathcal A_+\) define \(h_t = h(t,\cdot) \in \mathcal A_+\)
where \(h\) is the unique solution of \cref{eq:h} in \(\mathcal A_+\).
The mapping \(h_s \mapsto h_t\) is linear (when the multiplying scalar is
positive).
For \(p \ge 1\),
the growth of \(L^p\)-norm satisfies
\begin{align*}
\frac d{du} \int h_u(x)^p m_\infty(\dd x)
&\leq \frac{p(p-1)}{4} \int h_u(x)^p |b_u(x) - b_\infty(x)|^2 m_\infty(\dd x) \\
&\leq \frac{p(p-1)}{4} (M^F_{mm})^2 W_1^2(m_u, m_\infty) \int h_u(x)^p
m_\infty(\dd x)
\end{align*}
for \(u \in [s,t]\) a.e.,
by \cref{prop:lp-growth} and by Cauchy--Schwarz inequality
The existence of the stated constant \(C_{s,t,p}\) then follows from an
application of Grönwall's lemma.
For \(p \geq 1\), the mapping \(h_s \mapsto h_t \eqqcolon P_s^t h_s\) extends
uniquely to a continuous linear one by the density of \(\mathcal A_+\) in
\(L^p_+(m_\infty)\).
By the dominated convergence theorem we have
\(\lim_{t \to s} \int |h_t(x) - h_s(x)|^p m_\infty(\dd x) = 0\)
when \(h_s \in \mathcal A_+\),
using the fact that \((h_u)_{u \in [s,t]} \subset \mathcal A_+\) uniformly.
This property extends to general \(h_s \in L^p_+(m_\infty)\) by the density
in \cref{prop:std-alg-density}.
Hence \(P_s^t\) is a strongly continuous semigroup on \(L^p_+(m_\infty)\).
To recover the usual definition of strongly continuous semigroup we note that
\(L^p = L^p_+ - L^p_+\) and define \(P^t_s h \coloneqq P^t_s h_+ - P^t_s h_-\)
for \(h \in L^p(m_\infty)\).
\end{proof}

\begin{proof}[Proof of \cref{prop:l2-convergence}]
First suppose \(h_{t_0} \in \mathcal A_+\).
Thanks to \cref{prop:lp-growth} with \(p = 2\), we have
\begin{align*}
&\hspace{-1em}\frac{d}{dt} \int h(x)_t^2 m_\infty(\dd x) \\
&= - 2 \int |\nabla h_t(x)|^2 m_\infty(\dd x)
+ 2 \int h_t(x) \nabla h_t(x) \cdot \bigl(b_t(x) - b_\infty(x)\bigr)
m_\infty(\dd x) \\
&\leq - 2(1 - \varepsilon) \int |\nabla h_t(x)|^2 m_\infty(\dd x)
+ \frac{1}{2\varepsilon} \int h_t(x)^2 |b_t(x) - b_\infty(x)|^2 m_\infty(\dd x) \\
&\leq - 4(1 - \varepsilon) \rho
\biggl( \int h_t^2(x) m_\infty(\dd x) - 1 \biggr)
+ \frac{(M^F_{mm})^2}{2\varepsilon} W_1^2(m_t, m_\infty) \Vert h_t \Vert_2^2 \\
&= - 4(1-\varepsilon) \rho \Vert h_t - 1 \Vert_{2}^2
+ \frac{(M^F_{mm})^2}{2\varepsilon} W_1^2 (m_t, m_\infty) \Vert h_t \Vert_2^2,
\end{align*}
where we first use the Cauchy--Schwarz inequality before
applying the Poincaré inequality \cref{eq:poincare} satisfied by \(m_\infty\)
and the Lipschitz bound on
\(|b_t(x) - b_\infty(x)| = |D_m F(m_t,x) - D_m F(m_\infty,x)|\).
By the \(T_2\) inequality \cref{eq:t2} we have
\(W_1^2(m_t,m_\infty) \leq W_2^2(m_t, m_\infty)
\leq \rho^{-1} H(m_t | m_\infty)\).
Thanks to \cref{lem:mf-entropy-sandwich,thm:mf-entropy-convergence} we have
\begin{multline*}
H(m_t | m_\infty) \leq \mathcal F(m_t) - F(m_\infty)
\leq e^{-4\rho(t - t_0)} (\mathcal F(m_{t_0}) - \mathcal F(m_\infty)) \\
\leq \left(1 + \frac{M^F_{mm}}{\rho} + \frac{(M^F_{mm})^2}{2\rho^2} \right) e^{-4\rho(t - t_0)}H(m_{t_0} | m_\infty).
\end{multline*}
Finally note that the relative entropy satisfies, for $p > 1$,
\begin{equation}
\label{eq:entropy-lp}
H(m_{t_0} | m_\infty) \leq
\log \lVert h_{t_0} \rVert_p^{p/(p-1)}
\end{equation}
since by Jensen's inequality we have
\[
\exp \biggl( \int \log \bigl(h_{t_0}^{p-1}\bigr)\dd m_{t_0}\biggr)
\leq \int h_{t_0}^{p-1}\dd m_{t_0}
= \int h_{t_0}^p\dd m_{\infty}.
\]
Chaining up the previous three inequalities we obtain
\begin{multline*}
\frac{(M^F_{mm})^2}{2\varepsilon} W_1^2(m_t,m_\infty)
\leq \frac{(M^F_{mm})^2}{2\varepsilon} W_2^2(m_t,m_\infty) \\
\leq \frac{\rho \alpha^2}{2\varepsilon}
\biggl(1 + \alpha+\frac{\alpha^2}{2}\biggr)
\log \lVert h_{t_0} \rVert_2^2 e^{-4\rho(t-t_0)} \eqqcolon \Delta(t),
\end{multline*}
where we define \(\alpha \coloneqq M^F_{mm} / \rho\).
The decrease of \(L^2\)-norm then satisfies
\[\frac{d}{dt} \Vert h_t \Vert_2^2 \leq -\bigl(4\rho' - \Delta(t)\bigr)
\Vert h_t - 1\Vert_2^2 + \Delta(t)\]
with \(\rho' \coloneqq (1-\varepsilon)\rho\).
Thanks to Grönwall's lemma
and the fact that \(\int_s^{+\infty}\Delta(u)\dd u \leq \Delta(s)/4\rho\),
we obtain
\begin{align*}
&\hspace{-1em}\Vert h_t - 1 \Vert_2^2 \\
&\leq e^{ - 4\rho' (t-t_0) + \int_{t_0}^{t} \Delta(s)\dd s} \Vert h_{t_0} - 1
\Vert_2^2
+ \int_{t_0}^t e^{ - 4\rho' (t-s) + \int_s^t \Delta(u)du}
\Delta (s) \dd s \\
&\leq e^{\Delta(t_0) / 4\rho} \biggl(
e^{-4\rho'(t-t_0)} \Vert h_{t_0} - 1 \Vert_2^2
+ \int_{t_0}^t e^{ - 4\rho' (t - s)}
\Delta (s) \dd s
\biggr) \\
&\leq e^{\Delta(t_0)/4\rho} \biggl(
e^{-4\rho'(t-t_0)} \Vert h_{t_0} - 1 \Vert_2^2
+ \Delta(t_0) \int_{t_0}^t e^{ - 4\rho' (t-s)}
e^{-4\rho(s-t_0)} \dd s
\biggr) \\
&\leq e^{\Delta(t_0)/4\rho} \biggl(
e^{-4\rho'(t-t_0)} \Vert h_{t_0} - 1 \Vert_2^2
+ \frac{\Delta(t_0)}{4(\rho - \rho')} (e^{-4\rho'(t-t_0)} - e^{-4\rho(t - t_0)})
\biggr) \\
&\leq e^{\Delta(t_0)/4\rho} \biggl( \Vert h_{t_0} - 1 \Vert_2^2 +
\frac{\Delta(t_0)}{4\varepsilon \rho} \biggr) e^{-4\rho'(t - t_0)}.
\end{align*}

For general \(h_{t_0} \in L^2(m_\infty)\), we take an approximating sequence
\((h_{t_0}^n)_{n \in \mathbb N}\) in \(\mathcal A_+\) such that \(h^n_{t_0} \to
h_{t_0}\) in \(L^2(m_\infty)\)
according to \cref{prop:std-alg-stability}.
We have established that \(\Vert h_t^n - 1\Vert_2 \leq C e^{-\gamma t}\) where
\(h_t^n = P^t_{t_0} h^n_{t_0}\).
By the continuity shown in \cref{cor:lp-cont}, we have \(h_t^n \to h_t\) in
\(L^2(m_\infty)\).
Therefore, the inequality \cref{eq:l2-convergence} holds
for general \(h_{t_0}\in L^2(m_\infty)\).
\end{proof}

\subsubsection{Proof of \cref{prop:hypercontractivity}}

\begin{proof}[Proof of \cref{prop:hypercontractivity}]
First assume \(h_{t_0} \in \mathcal A_+\) so that \(h_t \in \mathcal A_+\) for
all \(t \geq t_0\) and that \(h_t\in \mathcal A_+\) uniformly on compact sets
of \([t_0, +\infty)\) thanks to \Cref{prop:std-alg-stability}.
Define the function \(\varphi(t) = \log \lVert h_t \rVert_{q(t)}\).
In particular, if \(q(t) = 0\), then
\(\varphi(t) = \int \log h_t(x) m_\infty(\dd x)\).
By the definition of the stable algebra \(\mathcal A_+\) we know \(\varphi(t)\)
is well defined for \(t \geq t_0\).
Moreover, it follows from Fubini's theorem that \(t \mapsto \varphi(t)\) is
absolutely continuous for \(t \geq t_0\) and its weak derivative reads
\begin{align*}
&\dot\varphi(t) \\
&= \frac{\dot q(t)}{q(t)^2 \int h_t(x)^{q(t)} m_\infty(\dd x) }
\biggl(\begin{aligned}[t]
&\int h_t(x)^{q(t)} \log h_t(x)^{q(t)} m_\infty(\dd x) \\
&- \int h_t(x)^{q(t)} m_\infty(\dd x)
\log \int h_t(x)^{q(t)} m_\infty(\dd x)\biggr)
\end{aligned} \\
&\phantom{={}}
+ \frac{q(t)-1}{\int h_t(x)^{q(t)} m_\infty(\dd x)}
\biggl(\begin{aligned}[t]
&- \int h_t(x)^{q(t) - 2} |\nabla h_t(x)|^2 m_\infty(\dd x) \\
&+ \int h_t(x)^{q(t) - 1} \nabla h_t(x) \cdot
\bigl(b_t(x) - b_\infty(x)\bigr) m_\infty(\dd x)\biggr).
\end{aligned}
\end{align*}
We recognize the term on the first line as the entropy,
\begin{multline*}
\int h_t(x)^{q(t)} \log h_t^{q(t)} m_\infty(\dd x)
- \int h_t(x)^{q(t)} m_\infty(\dd x) \log \int h_t(x)^{q(t)} m_\infty(\dd x) \\
= \Ent_{m_\infty} (h_t^{q(t)}),
\end{multline*}
which, by LSI \cref{eq:lsi}, has upper bound
\[
\Ent_{m_\infty} (h_t^{q(t)})
\leq \frac{1}{\rho} \Expect_{m_\infty}\bigl[ |\nabla h^{q(t)/2} |^2\bigr]
\leq \frac{q(t)^2}{4\rho} \int h_t(x)^{q(t) - 2}
|\nabla h_t(x)|^2 m_\infty(\dd x).
\]
By Cauchy--Schwarz, the second term on the second line satisfies
\begin{align*}
&\hspace{-1em}\int h_t(x)^{q(t) - 1} \nabla h_t(x)
\cdot \bigl(b_t(x) - b_\infty(x)\bigr) m_\infty(\dd x) \\
&\leq \varepsilon \int h_t(x)^{q(t) - 2} |\nabla h_t(x)|^2 m_\infty(\dd x)
+ \frac{1}{4\varepsilon} \biggl(\int h_t(x)^{q(t)} m_\infty(\dd x)\biggr)
\Vert b_t - b_\infty \Vert_\infty^2 \\
&\leq \varepsilon \int h_t(x)^{q(t) - 2} |\nabla h_t(x)|^2 m_\infty(\dd x)
+ \frac{(M^F_{mm})^2 W_1^2(m_t, m_\infty)}{4\varepsilon}
\int h_t(x)^{q(t)} m_\infty(\dd x).
\end{align*}
Therefore, for \(q_0 > 1\) (so that \(q(t) > 1, \dot q(t) > 0\)), we have
\(\dot\varphi(t) \leq \delta(t)\) while for \(q_0 < 1\) (so that \(q(t) < 1,
\dot q(t) < 0\)) we have \(\dot\varphi(t) \geq \delta(t)\).
To deal with the case \(q(t) = 0\) we use the continuity of \(t \mapsto
\varphi(t)\).
We have thus shown \cref{eq:hypercontractivity-q>1,eq:hypercontractivity-q<1}
for \(h_{t_0} \in \mathcal A_+\).

Now consider general \(h_{t_0} \in L^{q_0}_+(m_\infty)\).
In the case \(q_0 > 1\),
we use the density of \(\mathcal A_+\) (\cref{prop:std-alg-density}) to find a
sequence \((h^n_{t_0})_{n \in \mathbb N}\) in \(\mathcal A_+\)
with \(h^n_{t_0} \to h_{t_0} \) in \(L^{q_0}\).
To each \(h^n_{t_0}\) there exists a flow \(t \mapsto h^n_t\) in \(\mathcal
A_+\) satisfying \cref{eq:hypercontractivity-q>1}.
For \(t \geq t_0\), we also have \(h^n_t \to h_t\) in \(L^{q_0}\)
by the semigroup property in \cref{cor:lp-cont}
so that along a subsequence \(h^n_t \to h_t\) a.e.
By Fatou's lemma we obtain
\begin{multline*}
\log \biggl( \int h_t^{q(t)}(x) m_\infty(\dd x) \biggr)^{\!1/q(t)}
\leq \liminf_{n \to \infty} \biggl(
\int h^n_t(x)^{q(t)} m_\infty(\dd x) \biggr)^{\!1/q(t)} \\
\leq \liminf_{n \to \infty} \log \lVert h^n_{t_0} \rVert_{q_0}
+ \int_{t_0}^t \delta (s)\dd s
= \log \lVert h_{t_0} \rVert_{q_0} + \int_{t_0}^t \delta (s)\dd s.
\end{multline*}
So \cref{eq:hypercontractivity-q>1} is proved for general \(h_{t_0} \in
L^{q_0}\).
In the case \(q_0 < 1\),
we choose again by \cref{prop:std-alg-density} a sequence \((h^n_{t_0})_{n \in
\mathbb N}\) in \(\mathcal A_+\) such that \(h^n_{t_0} \to h_{t_0}\) in \(L^1\)
and \(\lim_{n \to \infty} \Vert h^n_{t_0} \Vert_{q_0} = \Vert h_{t_0}
\Vert_{q_0}\).
By the \(L^1\)-continuity, \(h^n_t \to h_t\) in \(L^1\) so that along a
subsequence \(h^n_t \to h_t\) pointwise \(m_\infty\)-a.e.
For \(q(t) > 0\) we have by Fatou's lemma
\[
\liminf_{n\to\infty}
\int \bigl(|h^n_t(x)| + 1 - |h^n_t(x)|^{q(t)}\bigr) m_\infty(\dd x)
\geq \int \bigl(|h_t(x)| + 1 - |h_t(x)|^{q(t)}\bigr) m_\infty(\dd x).
\]
Thus \(\limsup_{n\to\infty} \int |h^n_t(x)|^{q(t)} m_\infty(\dd x)
\leq \int |h_t(x)|^{q(t)} m_\infty(\dd x)\).
So taking \(\limsup\) on both sides of the inequality
\[
\log \lVert h^n_t \rVert_{q(t)} \geq \log \lVert h^n_{t_0} \rVert_{q_0} +
\int_{t_0}^t \delta(s)\dd s
\]
gives us \cref{eq:hypercontractivity-q<1}.
For \(q(t) < 0\) we have directly by Fatou
\[
\liminf_{n \to \infty} \int h^n_t(x)^{q(t)} m_\infty(\dd x)
\geq \int h_t(x)^{q(t)} m_\infty(\dd x)
\]
so that
\begin{multline*}
\log \lVert h_t \rVert_{q(t)} \geq \limsup_{n\to\infty}
\log \lVert h^n_t \rVert_{q(t)} \geq \limsup_{n\to\infty}
\log \lVert h^n_{t_0} \rVert_{q_0} +
\int_{t_0}^t \delta(s)\dd s \\
= \log \lVert h_{t_0} \rVert_{q_0} + \int_{t_0}^t \delta(s)\dd s.
\end{multline*}
To conclude we treat \(q(t) = 0\) by a continuity argument.
Take \(\varepsilon' \in (0, \varepsilon)\)
and let \(q'\) be the solution to \(\dot q' = 4(1 - \varepsilon') \rho (q'-1)\)
with \(q'(t_0) = q(t_0) = q_0 < 1\)
and \(\delta'(t) = \frac{1}{4\varepsilon'}(q'(t) - 1)
(M^F_{mm})^2 W_1^2(m_t,m_\infty)\).
We have \(q'(t) < q(t) = 0\) so that by previous discussions
\[
\log \lVert h_t \rVert_{q'(t)} \geq \log \lVert h_{t_0} \rVert_{q_0}
+ \int_{t_0}^t \delta'(s)\dd s,
\]
whereas \(\log \lVert h_t \rVert_{q(t)} \geq \log \lVert h_t \rVert_{q'(t)}\) by
the monotonicity of \(p\)-norm.
We take the limit \(\varepsilon' \to \varepsilon\) to obtain
\cref{eq:hypercontractivity-q<1}.
\end{proof}

\begin{rem}
\label{rem:defective-lsi}
The computations are similar to that for the hypercontractivity of a diffusion
process whose invariant measure \(m\) satisfies a \emph{defective LSI}, i.e.
for some \(c, \delta\geq 0\),
\[
\forall f \in C^1_b(\mathbb R^d),\qquad
\Ent_m (f^2) \leq c \Expect_m [|\nabla f|^2] + \delta \Expect_m [|f|^2].
\]
See \cite[Chapter 5]{bgl} and \cite[Chapter 2]{lsi} for the link between
defective LSI and hypercontractivity.
\end{rem}

\subsection{Proofs of \cref{thm:lp-convergence,thm:mf-uniform-concentration}}
\label{sec:proof-thm:lp-convergence}

After showing the \(L^2\)-convergence and the hypercontractivity, we are
finally ready to give the proof of \cref{thm:lp-convergence}.

\begin{proof}[Proof of \cref{thm:lp-convergence}]
We will first use \cref{prop:hypercontractivity} to show that after a finite
time \(h\) lies in \(L^2(m_\infty)\),
then use \cref{prop:l2-convergence} to show that its \(L^2(m_\infty)\)-norm
diminishes exponentially
and finally apply \cref{prop:hypercontractivity} again to extend this result to
all \(L^p\).

To this end, let $\rho' \in (0,\rho)$
be arbitrary and set \(\varepsilon = 1 - \rho' / \rho\).
Define
\(\dot q_1(t) = 4(1 - \varepsilon) \rho\bigl(q_1(t) - 1\bigr)\)
with \(q_1(0) = p_0\), and we know
\[
q_1(s) = (p_0-1) \exp\bigl(4(1-\varepsilon)\rho s\bigr) + 1.
\]
Since \(p_0 > 1\), \(q_1\) is exponentially increasing.
If \(p_0\in(1,2)\) we set \(t_1 = (4(1 - \varepsilon) \rho)^{-1} \log
\frac{1}{p_0-1}\).
This definition ensures that \(q_1(t_1) = 2\).
Otherwise if \(p_0 \geq 2\), we simply set \(t_1 = 0\).
Thus, in both cases, we have
\[
t_1 = \frac{1}{4(1-\varepsilon)\rho} \log \frac{1}{(p_0 - 1) \wedge 1}.
\]
By the hypercontractivity \cref{eq:hypercontractivity-q>1} in
\cref{prop:hypercontractivity},
we have
\[
\Vert h_{t_1} \Vert_2 \leq \exp \biggl(\int_0^{t_1} \delta_1(s)\dd s\biggr)
\Vert h_{0} \Vert_{p_0},
\]
where \(\delta_1(s) = \frac{1}{4\varepsilon} (q_1(s) - 1)
(M^F_{mm})^2 W_1^2(m_s, m_\infty)\). On the other hand,
we can control the Wasserstein distance $W_1^2(m_s, m_\infty)$ as follows:
\begin{align*}
W_1^2(m_s, m_\infty) \leq W_2^2(m_s, m_\infty)
&\leq \rho^{-1} H(m_s | m_\infty) \\
&\leq \rho^{-1} \bigl( \mathcal F(m_s) - \mathcal F(m_\infty) \bigr) \\
&\leq \rho^{-1} \bigl( \mathcal F(m_0) - \mathcal F(m_\infty) \bigr)e^{-4\rho s} \\
&\leq \rho^{-1}
\biggl( 1 + \frac{M^F_{mm}}{\rho} + \frac{(M^F_{mm})^2}{2\rho^2} \biggr)
H(m_0 | m_\infty) e^{-4\rho s}\\
&\leq \rho^{-1}
\biggl( 1 + \frac{M^F_{mm}}{\rho} + \frac{(M^F_{mm})^2}{2\rho^2} \biggr)
\log \lVert h_0 \rVert_{p_0}^{p_0/(p_0 - 1)}e^{-4\rho s},
\end{align*}
thanks to the $T_2$ inequality \cref{eq:t2}, \cref{thm:mf-entropy-convergence},
\cref{lem:mf-entropy-sandwich} and the inequality \cref{eq:entropy-lp}.
Setting $\alpha \coloneqq M^F_{mm} / \rho$
and $P(\alpha) = \alpha^2 + \alpha^3 + \alpha^4\!/2$, we get
\begin{align*}
\int_0^{t_1} \delta_1 (s)\dd s
&\leq \frac{M^F_{mm}p_0}{4\varepsilon(p_0-1)}
\biggl( \alpha + \alpha^2 + \frac{\alpha^3}{2} \biggr)
\log \lVert h_0 \rVert_{p_0}
\int_0^{t_1} (q_1(s) - 1)\dd s \\
&\leq \frac{M^F_{mm}p_0}{4\varepsilon(p_0-1)}
\biggl( \alpha + \alpha^2 + \frac{\alpha^3}{2} \biggr)
\log \lVert h_0 \rVert_{p_0}
\frac{1}{4(1-\varepsilon)\rho} (2 - p_0)_+ \\
&\leq \frac{p_0(2-p_0)_+}
{16(p_0-1)\varepsilon(1-\varepsilon)}
P(\alpha)
\log \lVert h_0\rVert_{p_0} \eqqcolon M \log \lVert h_0 \rVert_{p_0}.
\end{align*}
And thus,
$\Vert h_{t_1} \Vert_2 \leq \Vert h_0 \Vert_{p_0}^{1 + M}$.
By \cref{prop:l2-convergence} we know that for all \(t \in [t_1,+\infty)\),
\begin{align*}
\Vert h_t \Vert_2^2  - 1
&\leq \exp \biggl( \frac{P(\alpha)}{4\varepsilon}
\log \lVert h_{t_1} \rVert_2 \biggr)
\biggl( \lVert h_{t_1} \rVert_2^2 - 1
+ \frac{P(\alpha)}{4\varepsilon^2} \log \lVert h_{t_1} \rVert_2 \biggr)
e^{-4(1-\varepsilon)\rho(t-t_1)} \\
&\leq \lVert h_{t_1} \rVert_{2}^{P(\alpha)/4\varepsilon}
\biggl( 1 + \frac{P(\alpha)}{{8}\varepsilon^2} \biggr)
\bigl(\lVert h_{t_1}\rVert_2^2 - 1\bigr) e^{-4(1-\varepsilon)\rho(t-t_1)} \\
&\leq \biggl( 1 + \frac{P(\alpha)}{8\varepsilon^2}\biggr)
H_1^{P(\alpha)/4\varepsilon} \bigl( H_1^2 - 1 \bigr)
e^{-4(1-\varepsilon)\rho(t-t_1)},
\end{align*}
for $H_1$ being the upper bound of $\lVert h_{t_1}\rVert_2$ defined by
\[
\log H_1 =
\biggl( 1 + \frac{p_0(2-p_0)_+P(\alpha)}{16(p_0-1)\varepsilon(1-\varepsilon)}
\biggr)
\log \lVert h_0\rVert_{p_0}.
\]

Now we define \(\tau_p\) by
\begin{align*}
\tau_p &= \begin{cases}
t_1 + \frac{1}{4(1-\varepsilon)\rho} \log\bigl((p-1)\vee 1\bigr)
& \text{if}~p > 1, \\
t_1& \text{if}~p \in (0,1) \\
t_1 + \frac{1}{4(1-\varepsilon)\rho} \log\bigl(2(1-p)\bigr)
& \text{if}~p \leq 0
\end{cases} \\
&= \begin{cases}
\frac{1}{4(1-\varepsilon)\rho} \log \frac{(p-1)\vee 1}{(p_0 - 1)\wedge 1}
& \text{if}~p \geq 0, \\
\frac{1}{4(1-\varepsilon)\rho} \log \frac{2(1-p)}{(p_0-1)\wedge 1}
& \text{if}~p < 0, \\
\end{cases}
\end{align*}
In the case \(p > 1\),
for \(t \geq \tau_p\) we set
\(t_2 = t - (4(1-\varepsilon)\rho)^{-1} \log\bigl((p-1)\vee 1\bigr) \geq t_1\)
and let \(q_2\) solves \(\dot q_2(t) = 4(1 - \varepsilon)
\rho \bigl(q_2(t) - 1\bigr)\) with \(q_2(t_2) = 2\).
Our choice ensures \(q_2(t) = 2 \vee p \geq p\).
By the hypercontractivity \cref{eq:hypercontractivity-q>1} we have
\[
\Vert h_t \Vert_{q_2(t)} \leq \exp \biggl(\int_{t_2}^t \delta_2(s)\dd s\biggr)
\Vert h_{t_2} \Vert_2,
\]
where $\delta_2(s) = \frac{1}{4\varepsilon} \bigl(q_2(s) - 1\bigr)
(M^F_{mm})^2 W_1^2(m_s, m_\infty)$.
The integral of $\delta_2$ can be controlled in the same way
as we did to push $p_0 \to 2$ by hypercontractivity:
\begin{align*}
\int_{t_2}^t \delta_2(s)\dd s
&\leq \frac{M^F_{mm}p_0}{4\varepsilon(p_0-1)}
\biggl( \alpha + \alpha^2 + \frac{\alpha^3}{2} \biggr)
\log \lVert h_0 \rVert_{p_0}
\int_{t_2}^{t} \bigl(q_2(s) - 1\bigr)\dd s \cdot e^{-4\rho t_2} \\
&\leq \frac{p_0P(\alpha)}{16(p_0-1)\varepsilon(1-\varepsilon)}
\log \lVert h_0\rVert_{p_0} (p - 2)_+ \cdot e^{(1-\varepsilon)^{-1}\log((p-1)\vee1)}e^{-4\rho t}.
\end{align*}
The \(p\)-norm then satisfies
\begin{align*}
\log \lVert h_t \rVert_p &\leq \log \lVert h_t \rVert_{q_2(t)} \\
&\leq \log\lVert h_{t_2} \rVert_2
+ \frac{p_0(p-2)_+P(\alpha)}
{16(p_0-1)\varepsilon(1-\varepsilon)} \log \lVert h_0 \rVert_{p_0} \cdot e^{(1-\varepsilon)^{-1}\log((p-1)\vee1)}e^{-4\rho t} \\
&\leq \frac 12 \bigl(\Vert h_{t_2} \Vert_2^2 - 1\bigr)
+ \frac{p_0(p-2)_+P(\alpha)}
{16(p_0-1)\varepsilon(1-\varepsilon)} \log \lVert h_0 \rVert_{p_0} \cdot e^{(1-\varepsilon)^{-1}\log((p-1)\vee1)}e^{-4\rho t} \\
&\leq
\frac 12\biggl( 1 + \frac{P(\alpha)}{8\varepsilon^2}\biggr)
H_1^{P(\alpha)/4\varepsilon} \bigl( H_1^2 - 1 \bigr)
e^{-4(1-\varepsilon)\rho(t_2 - t_1)} \\
&\mathrel{\hphantom{\leq}} \negmedspace {}
\quad + \frac{p_0(p-2)_+P(\alpha)}
{16(p_0-1)\varepsilon(1-\varepsilon)} \log \lVert h_0 \rVert_{p_0} \cdot e^{(1-\varepsilon)^{-1}\log((p-1)\vee1)}e^{-4\rho t} \\
&\leqslant \frac 12\biggl( 1 + \frac{P(\alpha)}{8\varepsilon^2}\biggr)
H_1^{P(\alpha)/4\varepsilon} \bigl( H_1^2 - 1 \bigr)
e^{-4(1-\varepsilon)\rho(t - \tau_p)} \\
&\mathrel{\hphantom{\leq}} \negmedspace {}
\quad + \frac{p_0(p-2)_+P(\alpha)}
{16(p_0-1)\varepsilon(1-\varepsilon)} \log \lVert h_0 \rVert_{p_0} \cdot e^{(1-\varepsilon)^{-1}\log((p-1)\vee1)}e^{-4\rho t}.
\end{align*}
So the upper bound in \cref{eq:lp-convergence} is established.
The lower bound follows from the monotonicity of $p$-norm:
we have $\log \lVert h_t \rVert_p \ge \log \lVert h_t \rVert_1 = 0$.

For \(p \in (0,1)\), we observe Hölder's inequality
\[
\biggl(\int h^p m_\infty\biggr)^{\!1/(2 - p)}
\biggl(\int h^2 m_\infty\biggr)^{\!(1 - p)/(2 - p)}
\geq \int h m_\infty = 1,
\]
so that for \(t \geq \tau_p = t_1\) we have
\(\log \lVert h_t \rVert_p \geq - \frac{2(1 - p)}{p} \log \lVert h_t \rVert_2\).
Thus we obtain the desired bound by inserting the upper bound
for $\lVert h_t\rVert_2$.

Finally we treat \(p \leq 0\).
Given \(t \geq \tau_p\),
set \(t_3 = t - \bigl(4(1-\varepsilon)\rho\bigr)^{-1}
\log\bigl(2(1-p)\bigr) \geq t_1\)
and let \(q_3\) solves
\(\dot q_3(t) = 4(1-\varepsilon)\rho\bigl(q_3(t)-1\bigr)\)
with \(q_3(t_3) = \frac 12\).
Our choice ensures \(q_3(t) = p\).
Define \(\delta_3(s) = \frac{1}{4\varepsilon}\bigl(q_3(s) - 1\bigr)
(M^F_{mm})^2 W_1^2(m_s,m_\infty)\).
It satisfies, as done in the previous steps,
\[
\int_{t_3}^t \delta_3(s)\dd s
\geq - \frac{p_0\bigl( \frac 12 - p\bigr)P(\alpha)}
{16(p_0-1)\varepsilon(1-\varepsilon)} \log \lVert h_0 \rVert_{p_0} \cdot e^{-4\rho t_3}.
\]
We obtain, by the reverse hypercontractivity \cref{eq:hypercontractivity-q<1},
\begin{align*}
\log \lVert h_t \rVert_p
&\geq \log \lVert h_{t_3} \rVert_{\frac 12}
+ \int_{t_3}^t \delta_3(s)\dd s \\
&\geq - 2 \log \lVert h_{t_3} \rVert_2
- \frac{p_0\bigl( \frac 12 - p\bigr)P(\alpha)}
{16(p_0-1)\varepsilon(1-\varepsilon)} \log \lVert h_0 \rVert_{p_0}\cdot e^{(1-\varepsilon)^{-1}\log(2(1-p))}e^{-4\rho t} \\
&= - \log\bigl(1 + \Vert h_{t_3} - 1\Vert_2^2\bigr) \\
&\mathrel{\hphantom{\geq}} \negmedspace {}
\quad - \frac{p_0\bigl( \frac 12 - p\bigr)P(\alpha)}
{16(p_0-1)\varepsilon(1-\varepsilon)} \log \lVert h_0 \rVert_{p_0} \cdot e^{(1-\varepsilon)^{-1}\log(2(1-p))}e^{-4\rho t}\\
&\geq - \lVert h_{t_3} - 1 \rVert_2^2
- \frac{p_0\bigl( \frac 12 - p\bigr)P(\alpha)}
{16(p_0-1)\varepsilon(1-\varepsilon)} \log \lVert h_0 \rVert_{p_0} \cdot e^{(1-\varepsilon)^{-1}\log(2(1-p))}e^{-4\rho t}\\
&\geq - \biggl( 1 + \frac{P(\alpha)}{8\varepsilon^2}\biggr)
H_1^{P(\alpha)/4\varepsilon} \bigl( H_1^2 - 1 \bigr)
e^{-4(1-\varepsilon)\rho(t-t_1)} \\
&\mathrel{\hphantom{\geq}} \negmedspace {}
\quad - \frac{p_0\bigl( \frac 12 - p\bigr)P(\alpha)}
{16(p_0-1)\varepsilon(1-\varepsilon)} \log \lVert h_0 \rVert_{p_0}\cdot e^{(1-\varepsilon)^{-1}\log(2(1-p))}e^{-4\rho t}.
\end{align*}
Thus, we have established the lower bound in \cref{eq:lp-convergence},
for both $p \in (0,1)$ and $p \leq 0$. To conclude,
we compare again the $p$-norm with the $1$-norm and use the monotonicity.
\end{proof}

To conclude the discussions about the mean field dynamics we show a lemma which
uses \(L^p\)-norms to control a ``cross entropy''-like quantities
and use it to obtain the uniform-in-time concentration of measure result in
\cref{thm:mf-uniform-concentration}.
The lemma will also be used in the proof of \cref{thm:poc}.

\begin{lem}
\label{lem:lp-bounds-ce}
Let \(\mu\), \(\nu \in \mathcal P(\mathbb R^d)\) and \(h : \mathbb R^d \to
(0,+\infty)\) be a measurable function.
Then for all \(p > 0\),
\begin{equation}
\label{eq:lp-bounds-ce}
- \frac 1p H(\nu | \mu) + \log \lVert h \rVert_{L^{-p}(\mu)}
\leq \int \log h \dd \nu
\leq \frac 1p H(\nu | \mu) + \log \lVert h \rVert_{L^p(\mu)}.
\end{equation}
\end{lem}
\begin{proof}
Let \(X\) be a measurable space,
\(\mu, \nu\) be probability measures on \(X\)
and \(U : X \to \mathbb R\) be a random variable.
We have the convex duality inequality
(see e.g. \cite[Corollary 4.14]{concentration})
\begin{equation}
\label{eq:duality-entropy}
\Expect_\nu [U] \leq H(\nu | \mu) + \log \Expect_\mu [e^U].
\end{equation}
The right hand side of the inequality is always well defined in \((-\infty,
+\infty]\).
Take \(U = p \log h\).
For \(p > 0\) we obtain
\[
\int \log h\dd \nu
\leq \frac 1p H(\nu | \mu) + \frac 1p \log \int e^{p\log h} \dd \mu
= \frac 1p H(\nu | \mu) + \log \lVert h \rVert_{L^p(\mu)},
\]
and for \(p < 0\) we obtain
\[
\int \log h\dd \nu \geq \frac 1p H(\nu | \mu) + \log \lVert h \rVert_{L^p(\mu)}.
\qedhere
\]
\end{proof}

\begin{proof}[Proof of \cref{thm:mf-uniform-concentration}]
Let \(f : \mathbb R^d \to \mathbb R\) be \(1\)-Lipschitz continuous
and define for \(t \geq 0\) the moment-generating function
\(\psi_{t,f}(\lambda) = \log \Expect_{m_t} e^{\lambda (f - \Expect_{m_t} f)}\).
The equality in \cref{eq:duality-entropy} can be attained
and therefore we have
(see also \cite[Corollary 4.14]{concentration})
\[
\psi_{t,f}(\lambda) = \sup_{\mu \ll m_t} \lambda (\Expect_\mu f - \Expect_{m_t}
f) - H(\mu | m_t).
\]
For each \(\mu \ll m_t\),
the first term satisfies
\begin{align*}
\Expect_\mu f - \Expect_{m_t} f
&\leq W_1(\mu, m_t)
\leq W_1(\mu, m_\infty) + W_1(m_t, m_\infty) \\
&\leq \sqrt {\frac 1 \rho H(\mu | m_\infty)} + W_1(m_t, m_\infty)
\end{align*}
by Talagrand's transport inequality \cref{eq:t2} for \(m_\infty\).
The second term satisfies
\begin{align*}
H(\mu | m_t)
&= \int \log \frac{\dd \mu}{\dd m_t} \dd\mu
= \int \biggl( \log \frac{\dd\mu}{\dd m_\infty} - \log h_t \biggr)\dd \mu \\
&= H(\mu | m_\infty) - \int \log h_t\dd\mu \\
&\geq H(\mu | m_\infty) - \frac 1p H(\mu | m_\infty) - \log \lVert h_t \rVert_p
\end{align*}
for \(p > 1\) by the previous \cref{lem:lp-bounds-ce}.
Hence for \(\lambda \geq 0\) the moment-generating function \(\psi_{t,f}\)
satisfies
\begin{align*}
\psi_{t,f} (\lambda) &\leq \sup_{\mu \ll m_t} \lambda \sqrt {\frac 1\rho H(\mu
| m_\infty)} + \lambda W_1(m_t,m_\infty)
- ( 1 - p^{-1}) H(\mu | m_\infty) + \log \lVert h_t \rVert_p \\
&\leq \frac{\lambda^2}{4(1 - p^{-1}) \rho}
+ \lambda W_1(m_t, m_\infty) + \log \lVert h_t \rVert_p.
\end{align*}
For \(r,\lambda \geq 0\) we have by Markov's inequality
\begin{align*}
m_t [ f - \Expect f \geq r ]
&\leq e^{-\lambda r} \Expect_{m_t} e^{\lambda (f - \Expect_{m_t} f)} \\
&\leq \exp \biggl(-\lambda r + \frac{\lambda^2}{4(1-p^{-1})\rho}
+ \lambda W_1(m_t,m_\infty) + \log \lVert h_t \rVert_p\biggl).
\end{align*}
Take \(\lambda = 2 (1 - p^{-1})\rho\).
We obtain
\begin{align*}
&\hspace{-1em}m_t [ f - \Expect f \geq r] \\
&\leq \exp \Biggl( - \biggl(1 - \frac 1p\biggr) \rho r^2
+ 2 \biggl( 1 - \frac 1p\biggr) \rho W_1(m_t, m_\infty) r
+ \log \lVert h_t \rVert_p \Biggr).
\end{align*}
The bound on \(m_t [f - \Expect f \leq - r]\) is obtained by applying the
previous inequality to \(-f\).
Given \(\rho' \in (0, \rho)\), find \(p > 1\) such that
\((1 - p^{-1})\rho = \rho'\).
The desired result follows from
\cref{thm:mf-entropy-convergence,thm:lp-convergence}.
\end{proof}

\begin{rem}
\label{rem:concentration-transport-method}
Our proof is based on the standard transport method for concentration
inequalities
and we refer readers to \cite[Chapter 6]{ledouxconcentration} and \cite[Chapter
8]{concentration} for an introduction to it.
In fact, our method allows us to prove a more general perturbative result:
if \(m\) satisfies a \(T_1\) inequality, \(h \in L^p_+(m)\) for \(p > 1\) and
\(\int hm = 1\),
then \(hm\) also has Gaussian concentration (albeit with a weaker constant).
\end{rem}

\section{Particle system}
\label{sec:ps}

\subsection{Proof of \cref{thm:ps-entropy-convergence}}
\label{sec:proof-thm:ps-entropy-convergence}

Before giving the proof of \cref{thm:ps-entropy-convergence} we first show two
lemmas on entropies.

\begin{lem}[Information inequalities]
\label{lem:info-ineqs}
Let \(X_1, \ldots, X_N\) be measurable spaces,
\(\mu\) be a probability measure on the product space
\(X = X_1 \times \cdots \times X_N\)
and \(\nu = \nu^1 \otimes \cdots \otimes \nu^N\) be a \(\sigma\)-finite measure.
Then
\begin{equation}
\label{eq:info-ineqs}
\sum_{i=1}^N H(\mu^i | \nu^i) \leq H(\mu | \nu)
\leq \sum_{i=1}^N \int
H\Bigl(\mu^{i|-i}(\cdot|\mathbf x^{-i})\Big| \nu^i\Bigr)\mu^{-i}
(\dd\mathbf x^{-i}).
\end{equation}
Here we set the rightmost term to \(+\infty\) if the conditional distribution
\(\mu^{i|-i}\) does not exist \(\mu^{-i}\)-a.e.
\end{lem}

\begin{proof}
The inequality is non-trivial only if \(\mu \ll \nu\)
and in this case we denote the relative density by \(f = \dd\mu/\!\dd\nu\).
For \(I \subset \{1,\ldots,N\}\), we define the conditional densities by
\[
f^{I|-I} (\mathbf x^I | \mathbf x^{-I}) =
\begin{cases}
\displaystyle \frac{f(\mathbf x^I, \mathbf x^{-I})}{\int f(\mathbf x^I,\mathbf
x^{-I}) \nu^{-I} (\dd\mathbf x^{-I})}&
\text{if}~\displaystyle \int f(\mathbf x^I,\mathbf x^{-I}) \nu^{-I}(\dd \mathbf
x^{-I}) > 0, \\
0& \text{otherwise.}
\end{cases}
\]
The conditional measures are defined via densities
\[
\mu^{I|-I}(\dd\mathbf x^I) = f^{I|-I} (\mathbf x^I | \mathbf x^{-I}) \nu^{I} (\dd
\mathbf x^I).
\]
In particular, we do not need the regularity of the underlying spaces
\(X_1,\ldots,X_N\) in order to apply disintegration theorems.
Define \(I_i = \{1,\ldots,i \}\) for \(i = 1,\ldots,N\).
The relative entropy admits the decomposition
\[
H(\mu | \nu)
= \sum_{i=1}^N \int H\Bigr( \mu^{i | I_{i-1}} (\cdot | \mathbf x^{I_{i-1}})
\Big| \nu^i \Bigr)\mu^{I_{i-1}} (\dd\mathbf x^{I_{i-1}}).
\]
We conclude by applying Jensen's inequality to the convex mappings \(\lambda^i
\mapsto H(\lambda^i | \nu^i)\).
\end{proof}

\begin{lem}
\label{lem:ps-entropy-sandwich}
Assume that \(F\) satisfies \cref{eq:convex}
and there exists a measure \(m_\infty \in \mathcal P_2(\mathbb R^d)\) verifying
\cref{eq:mf-foc}.
Then for all \(m^N \in \mathcal P_2(\mathbb R^{dN})\) of finite entropy, we have
\begin{equation}
\label{eq:ps-free-energy-bounds-entropy}
H(m^N | m_\infty^{\otimes N})
\leq \mathcal F^N(m^N) - N \mathcal F(m_\infty).
\end{equation}
\end{lem}

\begin{proof}
Let \(\mathbf X\) be a random variable distributed as \(m^N\).
By the convexity of \(F\)
we have
\begin{align*}
\hspace{-1em}\mathcal F^N&(m^N) - N\mathcal F(m_\infty) \\
&= \Expect [NF(\mu_{\mathbf X}) - NF(m_\infty)] + H(m^N) - NH(m_\infty) \\
&\geq \Expect \biggl[
N\int \frac{\delta F}{\delta m} (m_\infty, x) (\mu_{\mathbf X} - m_\infty)(\dd x)
\biggr] + H(m^N) - NH(m_\infty) \\
&= - \Expect \biggl[
N\int \log m_\infty(x) (\mu_{\mathbf X} - m_\infty)(\dd x) \biggr] + H(m^N) -
NH(m_\infty) \\
&= - \Expect \biggl[
N\int \log m_\infty(x) \mu_{\mathbf X}(\dd x) \biggr] + H(m^N) \\
&= - \int \sum_{i=1}^N \log m_\infty (x^i) m^N(\dd\mathbf x) + H(m^N)
= H(m^N | m_\infty^{\otimes N}). \qedhere
\end{align*}
\end{proof}

\begin{proof}[Proof of \cref{thm:ps-entropy-convergence}]

Let \(t_0 \geq 0\) be such that \(m_{t_0}\) has finite entropy and finite
second moment.
Since \(\nabla_i NF (\mu_\mathbf{x}) = D_m F(\mu_\mathbf{x}, x^i)\)
corresponds to the drift of \cref{eq:ps-sde},
we recognize the particle system flow of measure \(m^N_t\) as a linear Langevin
flow with the invariant measure \(m^N_\infty\), defined in
\cref{eq:def-ps-invariant-measure}.
In particular, \cref{prop:mf-energy-decrease} applied to this dynamics yields
\begin{equation}
\label{eq:ps-energy-decrease}
\frac{d\mathcal F^N(m^N_t)}{dt} = - I(m^N_t | m^N_\infty)
\end{equation}
for \(t \geq t_0\) a.e.
In the following we establish a lower bound of the relative Fisher information
\(I_t \coloneqq I(m^N_t | m^N_\infty)\)
in order to obtain the desired result. We divide the proof into several steps.

\paragraph{Regularity of conditional distribution.}
By the elliptic positivity (see e.g. \cite[Theorem 8.2.1]{fpkeq}),
we know that for all \(t > t_0\) and \(\mathbf x \in \mathbb R^{dN}\),
\(m^N_t(\mathbf x) > 0\) with explicit lower bound.
Let \(i \in \{1,\ldots,N\}\).
Define marginal density \(m^{N,-i}_t (\mathbf x^{-i})= \int m^N_t(\mathbf x)\dd
x^i\).
It is strictly positive everywhere by the positivity of \(m^N_t\)
and is lower semicontinuous (in \(\mathbf x^{-i}\))
thanks to the continuity of \(\mathbf x \mapsto m^N_t(\mathbf x)\)
and Fatou's lemma.
Since Fubini gives \(\int m^{N,-i}_t (\mathbf x^{-i})\dd\mathbf x^{-i} = 1\),
we have \(m^{N,-i}_t(\mathbf x^{-i}) < +\infty\) everywhere.
We are therefore able to define the conditional probability density
\[
m^{N,i|-i}_t (x^i | \mathbf x^{-i} )
= \frac{m^{N}_t (\mathbf x)}{m^{N,-i}_t(\mathbf x^{-i})}
= \frac{m^{N}_t (\mathbf x)}{\int m^{N}_t (\mathbf x)\dd x^{i}}
\]
which has generalized derivative in \(x^i\) and is strictly positive everywhere.

\paragraph{Decomposing Fisher componentwise.}
Using the conditional distributions, we can decompose the relative Fisher
information by
\begin{align*}
I_t
&= \int \biggl| \nabla \log \frac{m^N_t(\mathbf x)}{m^N_\infty(\mathbf
x)}\biggr|^2 m^N_t (\dd\mathbf x)
= \Expect \biggl[ \biggl| \nabla \log \frac{m^N_t (\mathbf
X_t)}{m^N_\infty(\mathbf X_t)}\biggr|^2 \biggr] \\
&= \sum_{i=1}^N \Expect \biggl[ \biggl| \nabla_{x^i} \log
\frac{m^{N,i|-i}_t(X^i_t | \mathbf X^{-i}_t) m^{N,-i}_t (\mathbf X^{-i}_t)}
{m^N_\infty (\mathbf X_t)}\biggr|^2 \biggr] \\
&= \sum_{i=1}^N \Expect \biggl[ \biggl| \nabla_{x^i} \log
\frac{m^{N,i|-i}_t(X^i_t | \mathbf X^{-i}_t)}
{m^{N}_\infty (\mathbf X_t)}\biggr|^2 \biggr] \\
&
= \sum_{i=1}^N \Expect \Bigl[
\Bigl| \nabla_{x^i} \log m^{N,i|-i}_t(X^i_t|\mathbf X^{-i}_t) + D_m
F(\mu_{\mathbf X_t}, X^i_t)\Bigr|^2 \Bigr].
\end{align*}

\paragraph{Change of empirical measure and componentwise LSI.}
We replace the empirical measure \(\mu_{\mathbf x}\) in \(D_m F\) by
\(\mu_{\mathbf x^{-i}}\).
Define \(\delta^i_1 (\mathbf x; y) = D_m F(\mu_{\mathbf x}, y) - D_m
F(\mu_{\mathbf x^{-i}}, y)\).
Take \(\varepsilon \in (0,1)\).
The Fisher information satisfies
\begin{align*}
I_t&
=\sum_{i=1}^N \Expect \biggl[
\Bigl| \nabla_{x^i} \log m^{N,i|-i}_t (X^i_t|\mathbf X_t^{-i})+ D_m
F(\mu_{\mathbf X^{-i}_t}, X^i_t) + \delta^i_1(\mathbf X_t; X^i_t)\Bigr|^2
\biggr] \\
&\geq\sum_{i=1}^N \Expect\left[
\begin{multlined}
(1 - \varepsilon) \Bigl| \nabla_{x^i} \log m^{N,i|-i}_t(X^i_t|\mathbf X_t^{-i})
+ D_m F(\mu_{\mathbf X_t^{-i}}, X^i_t)\Bigr|^2 \\
\quad - (\varepsilon^{-1} -1) |\delta^i_1(\mathbf X_t; X^i_t)|^2
\end{multlined}
\right]\\
&= (1 - \varepsilon) \sum_{i=1}^N \Expect
\Bigl[ I \Bigl(m^{N,i|-i}_t(\cdot|\mathbf X_t^{-i}) \Big| \hat \mu_{\mathbf
X_t^{-i}}\Bigr) \Bigr]
- (\varepsilon^{-1} - 1) \sum_{i=1}^N \Expect [|\delta^i_1(\mathbf X_t;
X^i_t)|^2],
\end{align*}
where we used the elementary inequality \((a + b)^2 \geq (1 - \varepsilon)
|a|^2 - (\varepsilon^{-1} - 1) |b|^2\)
and \(\hat \mu_{\mathbf x^{-i}}\) is the probability of density proportional to
$\exp \bigl( - \frac{\delta F}{\delta m} (\mu_{\mathbf x^{-i}}, x) \bigr)\dd x$.
Define the first error
\begin{equation}
\label{eq:Delta_1}
\Delta_1 \coloneqq \sum_{i=1}^N \Expect [|\delta^i_1(\mathbf X_t; X^i_t)|^2]
\coloneqq \sum_{i=1}^N \Expect
\bigl[|D_m F(\mu_{\mathbf X_t},X_t^i) - D_m F(\mu_{\mathbf X_t^{-i}}, X_t^i)|^2
\bigr].
\end{equation}
The previous inequality writes
\begin{equation}
\label{eq:I-lower-bound-1}
I_t \geq
(1 - \varepsilon) \sum_{i=1}^N \Expect
\Bigl[ I \Bigl(m^{N,i|-i}_t(\cdot|\mathbf X_t^{-i}) \Big| \hat \mu_{\mathbf
X_t^{-i}}\Bigr) \Bigr] - (\varepsilon^{-1} - 1) \Delta_1.
\end{equation}
We apply the uniform log-Sobolev inequality for \(\hat \mu_{\mathbf X^i_t}\)
and obtain
\begin{multline*}
\frac{1}{4\rho}I\Bigl(m^{N,i|-i}_t(\cdot|\mathbf X^{-i}_t) \Big| \hat
\mu_{\mathbf X^{-i}_t}\Bigr)
\geq H\Bigl(m^{N,i|-i}_t(\cdot|\mathbf X^{-i}_t) \Big| \hat \mu_{\mathbf
X^{-i}_t}\Bigr) \\
= \int \biggl( \log m^{N,i|-i}_t(x^i | \mathbf X^{-i}_t) + \frac{\delta
F}{\delta m} (\mu_{\mathbf X^{-i}_t}, x^i)\biggr)
m^{N,i|-i}_t(\dd x^i|\mathbf X^{-i}_t) + \log Z(\hat \mu_{\mathbf X^{-i}_t}).
\end{multline*}
Then we apply Jensen's inequality to \(\log Z (\hat\mu_{\mathbf x^{-i}})\) to
obtain
\[
\log Z (\hat\mu_{\mathbf X^{-i}_t})
\geq - \int \frac{\delta F}{\delta m} (\mu_{\mathbf X^{-i}_t}, x^i)
m_\infty(\dd x^i) - \int m_\infty(x^i) \log m_\infty(x^i) \dd x^i.
\]
Chaining the previous two inequalities and summing over \(i\),
we have
\begin{multline}
\label{eq:I-i|-i-lower-bound}
\frac{1}{4\rho} \sum_{i=1}^N
I\Bigl(m^{N,i|-i}_t(\cdot|\mathbf X^{-i}_t) \Big| \hat \mu_{\mathbf
X^{-i}_t}\Bigr)
\geq \sum_{i=1}^N \biggl[ \int \frac{\delta F}{\delta m} (\mu_{\mathbf
X_t^{-i}}, x^i) \\
\Bigl(m^{N,i|-i}_t(\dd x^i|\mathbf X_t^{-i}) - m_\infty(\dd x^i)\Bigr)
+ H\Bigl(m^{N,i|-i}_t(\cdot|\mathbf X_t^{-i})\Bigr) - H(m_\infty)
\biggr].
\end{multline}

\paragraph{Another change of empirical measure.}
We wish to change back \(\mu_{\mathbf x^{-i}} \to \mu_{\mathbf x}\) in
\cref{eq:I-i|-i-lower-bound}.
Define
\(
\delta^i_2 (\mathbf x; y) \coloneqq \frac{\delta F}{\delta m} (\mu_{\mathbf
x^{-i}}, y)
- \frac{\delta F}{\delta m} (\mu_{\mathbf x}, y)
\)
and the second error
\begin{equation}
\Delta_2
\coloneqq \sum_{i=1}^N \int \delta^i_2 (\mathbf x; x^i) m^N_t(\dd\mathbf x)
- \sum_{i=1}^N \iint \delta^i_2 (\mathbf x; x') m_\infty(\dd x')
m^N_t(\dd\mathbf x). \label{eq:Delta_2}
\end{equation}
Then we obtain by taking expectations on both sides of
\cref{eq:I-i|-i-lower-bound}
\begin{multline}
\frac{1}{4\rho} \sum_{i=1}^N \Expect
\Bigl[I\Bigl(m^{N,i|-i}_t(\cdot|\mathbf X^{-i}_t) \Big| \hat \mu_{\mathbf
X^{-i}_t}\Bigr)\Bigr]
\geq N \Expect \biggl[ \int \frac{\delta F}{\delta m} (\mu_{\mathbf X_t}, y)
(\mu_{\mathbf X_t} - m_{\infty}) (\dd y) \biggr] \\
+ \sum_{i=1}^N \Expect H\Bigl(m^{N,i|-i}_t(\cdot|\mathbf X_t^{-i})\Bigr)
- N H(m_\infty)
+ \Delta_2.
\end{multline}
Thanks to the convexity of \(F\), the first term satisfies the tangent
inequality
\begin{multline}
\label{eq:ps-entropy-convergence-tangent}
N \Expect \biggl[ \int \frac{\delta F}{\delta m} (\mu_{\mathbf X_t}, y)
(\mu_{\mathbf X_t} - m_{\infty}) (\dd y) \biggr]
\geq N \Expect \bigl[ F (\mu_{\mathbf X_t}) - F(m_{\infty}) \bigr] \\
= F^N(m^N_t) - NF(m_{\infty}).
\end{multline}
For the second term we apply the information inequality \cref{eq:info-ineqs} to
obtain
\[
\sum_{i=1}^N \Expect^{-i}\Bigl[
H\Bigl(m^{N,i|-i}_t(\cdot|\mathbf X_t^{-i})\Bigr)
\Bigr]\geq H(m^N_t).
\]
Hence,
\begin{multline*}
\sum_{i=1}^N \Expect
\Bigl[I\Bigl(m^{N,i|-i}_t(\cdot|\mathbf X^{-i}_t) \Big| \hat \mu_{\mathbf
X^{-i}_t}\Bigr)\Bigr]
\\
\geq 4\rho\bigl(F^N(m^N_t) - NF(m_\infty) + H(m^N_t) - NH(m_\infty) + \Delta_2
\bigr).
\end{multline*}
Using \cref{eq:I-lower-bound-1}
and recalling the definition of free energies
\(\mathcal F(m) = F(m) + H(m)\),
\(\mathcal F^N(m^N) = F^N(m^N) + H(m^N)\),
we obtain
\begin{equation}
\label{eq:I-lower-bound-2}
I_t = I(m^N_t | m^N_\infty)
\geq 4\rho (1 - \varepsilon)
\bigl(\mathcal F^N (m^N_t) - N\mathcal F(m_\infty) + \Delta_2\bigr)
- (\varepsilon^{-1} - 1) \Delta_1.
\end{equation}

\paragraph{Estimate of the errors $\Delta_1$, $\Delta_2$.}
The transport plan between \(\mu_{\mathbf x}\) and \(\mu_{\mathbf x^{-i}}\)
\begin{equation}
\label{eq:transport-plan-empiricals}
\pi^i = \frac 1N \sum_{j\neq i} \delta_{(x^j,x^j)}
+ \frac 1{N(N-1)} \sum_{j \neq i} \delta_{(x^j,x^i)}
\end{equation}
gives the bound
\[
W_1 (\mu_{\mathbf x}, \mu_{\mathbf x^{-i}})
\leq \frac{1}{N(N-1)} \sum_{j\neq i}|x^j - x^i|.
\]
We use this transport plan to bound the errors $\Delta_1$, $\Delta_2$.

Let us treat the first error \(\Delta_1\).
Since \(m \mapsto D_m F(m,x)\) is \(M^F_{mm}\)-Lipschitz continuous in \(W_1\) metric,
we have
\[
|\delta^i_1(\mathbf x; y)| \leq M^F_{mm} W_1 (\mu_{\mathbf x}, \mu_{\mathbf x^{-i}}) \leq \frac{M^F_{mm}}{N(N-1)} \sum_{j=1,j\neq i}^N |x^j - x^i|.
\]
Under the \(L^2\)-optimal transport plan
\(\Law \bigl((X^i_t)_{i=1}^N, (\tilde X^i_\infty)_{i=1}^N\bigr)
\in \Pi(m^N_t, m^{\otimes N}_\infty)\)
we have
\begin{align*}
\Delta_1 &= \sum_{i=1}^N \Expect[ |\delta^i_1 (\mathbf X_t; X^i_t)|^2 ]
\leq (M^F_{mm})^2 \sum_{i=1}^N \Expect[W_1^2(\mu_{\mathbf X_t}, \mu_{\mathbf X^{-i}_t})] \\
&\leq \frac{(M^F_{mm})^2}{N(N-1)} \Expect \biggl[ \sum_{\substack{1 \leq i,j \leq N \\ i\neq j}}|X^j_t - X^i_t |^2 \biggr] \\
&\leq \frac{3(M^F_{mm})^2}{N(N-1)} \Expect \biggl[ \sum_{\substack{1 \leq i,j \leq N \\ i\neq j}}
\left(|X^i_t - \tilde X^i_\infty |^2 + |\tilde X^i_\infty - \tilde X^j_\infty|^2 + |X^j_t - \tilde X^j_\infty|^2\right) \biggr] \\
&\leq \frac{3(M^F_{mm})^2}{N(N-1)}
\biggl(
2(N-1)\Expect \biggl[
\sum_{i=1}^N|X^i_t - \tilde X^i_\infty|^2 \biggr]
+ N(N-1) \Expect [|\tilde X^1_\infty - \tilde X^2_\infty|^2]
\biggr).
\end{align*}
The first term
\(\Expect [\sum_{i=1}^N|X^i_t - \tilde X^i_\infty|^2 ]\)
is the Wasserstein distance \(W_2^2(m^N_t, m_\infty^{\otimes N})\),
while the second \(\Expect [|\tilde X^1_\infty - \tilde X^2_\infty|^2]\)
equals \(2 \Var m_\infty\).
Hence the first error satisfies the bound
\begin{equation}
\Delta_1 \leq 6(M^F_{mm})^2
\biggl( \frac 1N W_2^2(m^N_t, m^{\otimes N}_\infty) + \Var m_\infty \biggr).
\label{eq:bound-Delta_1}
\end{equation}

Now treat the second error \(\Delta_2\).
The Lipschitz constant of the mapping
\(y \mapsto \delta^i_2(\mathbf x; y) = \frac{\delta F}{\delta m}(\mu_{\mathbf
x^{-i}}, y) - \frac{\delta F}{\delta m}(\mu_{\mathbf x}, y)\)
is controlled by
\[
| \nabla_y \delta^i_2(\mathbf x; y) | = | D_m F(\mu_{\mathbf x}, y) - D_m
F(\mu_{\mathbf x^{-i}}, y)|
\leq M^F_{mm} W_1 (\mu_{\mathbf x}, \mu_{\mathbf x^{-i}}).
\]
Hence we have
\[
|\delta^i_2(\mathbf x; y) - \delta^i_2(\mathbf x; y')|
\leq M^F_{mm} W_1(\mu_{x}, \mu_{x^{-i}})|y - y'|.
\]
Use Fubini's theorem to first integrate \(x'\) in the definition of the second
error \cref{eq:Delta_2}
and let \(\tilde X'_\infty\) be independent from \(\mathbf X_t\).
Then we obtain
\begin{align*}
|\Delta_2|
&\leq \sum_{i=1}^N \int \biggl(
\int |\delta^i_2(\mathbf x; x^i) - \delta^i_2(\mathbf x; x') | m_\infty (\dd x')
\biggr) m^N_t(\dd\mathbf x) \\
&\leq \sum_{i=1}^N \iint \frac{M^F_{mm}}{N(N-1)}
\sum_{j =1\,j\neq i}^N |x^j - x^i| |x' - x^i| m_\infty(\dd x')
m^N_t(\dd\mathbf x) \\
&= \frac{M^F_{mm}}{N(N-1)}\sum_{\substack{i,j=1 \\ i\neq j}}^N
\Expect [ |X^j_t - X^i_t| |X^i_t - \tilde X'_\infty| ] \\
&\leq \frac{M^F_{mm}}{2N(N-1)}
\biggl(\sum_{\substack{i,j=1 \\ i\neq j}}^N
\Expect |X^i_t - X^j_t|^2
+ (N-1)\sum_{i=1}^N \Expect |X^i_t - \tilde X'_\infty|^2
\biggr).
\end{align*}
Using the same method we used for \(\Delta_1\),
we control the first term by
\[
\sum_{\substack{i,j=1 \\ i\neq j}}^N
\Expect |X^i_t - X^j_t|^2
\leq 6 N(N-1)
\biggl(\frac 1N W_2^2(m_t^N, m_\infty^{\otimes N}) + \Var m_\infty \biggr).
\]
For the second term we work again under the \(L^2\)-optimal plan
\(\Law\bigl((X^i_t)_{i=1}^N, (\tilde X^i_\infty)_{i=1}^N\bigr)
\in \Pi(m^N_t, m^{\otimes N}_\infty)\)
and let \(\tilde X'_\infty\) remain independent from the other variables.
We have
\begin{multline*}
\sum_{i=1}^N \Expect |X^i_t - \tilde X'_\infty|^2
\leq 2 \sum_{i=1}^N \Bigl(\Expect |X^i_t - \tilde X^i_\infty|^2 + |\tilde
X^i_\infty - \tilde X'_\infty|^2 \Bigr) \\
= 2 N \biggl(\frac 1NW_2^2(m_t^N, m_\infty^{\otimes N}) + 2 \Var
m_\infty\biggr).
\end{multline*}
As a result,
\begin{equation}
\label{eq:bound-Delta_2}
|\Delta_2| \leq M^F_{mm}
\biggl( \frac 4N W_2^2(m_t^N,m_\infty^{\otimes N}) + 5\Var m_\infty \biggr).
\end{equation}

\paragraph{Conclusion.}
Inserting the bounds on the errors \cref{eq:bound-Delta_1,eq:bound-Delta_2}
to the lower bound of Fisher information \cref{eq:I-lower-bound-2},
we obtain
\begin{align*}
I(m^N_t | m^N_\infty) &\geq 4\rho (1 - \varepsilon)
\bigl(\mathcal F^N(m^N_t) - N\mathcal F(m_\infty)\bigr) \\
&\quad- \bigl(16\rho M^F_{mm} + 6(\varepsilon^{-1} - 1) (M^F_{mm})^2\bigr)
\frac 1N W_2^2(m_t^N, m^{\otimes N}_\infty) \\
&\quad- \bigl(20\rho M^F_{mm} + 6(\varepsilon^{-1} - 1) (M^F_{mm})^2\bigr)
\Var m_\infty.
\end{align*}
Thanks to the Poincaré inequality \cref{eq:poincare} for \(m_\infty = \hat
m_\infty\),
its variance satisfies
\[2\rho \Var_{m_\infty} (x^i)
\leq \Expect_{m_\infty}\bigl[|\nabla x^i|^2\bigr]
= 1.\]
So \(\Var m_\infty = \sum_{i=1}^d \Var_{m_\infty} (x^i) \leq d/2\rho\).
Using the \(T_2\)-transport inequality \cref{eq:t2} for \(m^{\otimes
N}_\infty\) and the entropy sandwich \cref{lem:ps-entropy-sandwich}
we control the transport cost by
\[
W_2^2(m^N_t, m^{\otimes N}_\infty)
\leq \frac 1{\rho} H(m^N_t | m^{\otimes N}_\infty)
\leq \frac 1{\rho} \bigl(\mathcal F^N(m^N_t) - N \mathcal F(m_t)\bigr).
\]
In the end we obtain
\begin{multline*}
\frac{d\mathcal F^N(m^N_t)}{dt} = - I(m^N_t | m^N_\infty) \\
\leq
- \Biggl( 4(1-\varepsilon)\rho
- \frac {M^F_{mm}}N \biggl( 16 + 6(\varepsilon^{-1} - 1)
\frac{M^F_{mm}}{\rho} \biggl) \Biggl)
\bigl( \mathcal F^N(m^N_t) - N\mathcal F(m_\infty) \bigr) \\
+ dM^F_{mm}\biggl(10 + 3(\varepsilon^{-1} - 1) \frac{M^F_{mm}}{\rho}\biggr).
\end{multline*}
We conclude by applying Grönwall's lemma to the differential inequality above
and using the entropy inequality of \cref{lem:ps-entropy-sandwich}.
\end{proof}

\begin{rem}
If the initial condition \(m^N_0\) of the particle system is a tensor product
\((m_0)^{\otimes N}\),
one may expect the (non-uniform) convergence of the free energy \(\frac 1N
\mathcal F(m^N_t) \to \mathcal F(m_t)\) for all \(t \geq 0\).
If this is true, one can take the limit \(N \to \infty\) to recover the result
of \cref{thm:mf-entropy-convergence}.
However, while the convergence of the regular part \(\frac 1N F(m^N_t) \to
F(m_t)\) can be expected from the finite-time Wasserstein convergence
\(\frac 1N \sup_{t \in [0,T]} W_2(m^N_t , m^{\otimes N}_t) \to 0\),
the convergence of entropy \(H(m^N_t) \to H(m^{\otimes N}_t)\) is more
difficult to obtain.
\end{rem}

\begin{rem}
We used the convexity of \(F\) to achieve two things in the proof:
(i) the existence of mean field invariant measure \(m_\infty\);
and (ii) to derive
\cref{eq:ps-entropy-convergence-tangent,eq:ps-free-energy-bounds-entropy}.
Under mild assumptions (i) can also be obtained by a Schauder-type fixed point
theorem for the mapping \(m \mapsto \hat m\),
or by finding stationary points of the mean field free energy \(\mathcal F\).
For (ii), if \(F\) is only \(-\kappa\)-semi-convex around \(m_\infty\), in the
sense that
\[
F(m) - F(m_\infty) \geq \int \frac{\delta F}{\delta m}(m_\infty,x)
(m-m_\infty)(\dd x) - \frac{\kappa^2}{2} W_2^2(m,m_\infty),
\]
we can expect our method to apply as long as \(\kappa\) is sufficiently small.
\end{rem}

\subsection{Proofs of \cref{thm:poc,cor:poc}}
\label{sec:poc}

\begin{proof}[Proof of \cref{thm:poc}]
We separate the proof in two parts, each dealing with
the finite-time and long-time propagation of chaos respectively.
In each part, we shall first control the Wasserstein distance
$W_2(m^N_t, m_t^{\otimes N})$
between the particle system and the tensorized mean field system,
and then control their relative entropy
$H(m^N_t|m_t^{\otimes N})$.

\paragraph{Finite-time behavior.}

We shall use the synchronous coupling method
to control the Wasserstein distance between $m^N_t$ and $m_t^{\otimes N}$
and use Girsanov's theorem to control their relative entropy
on finite time intervals.
This may be considered folklore by specialists
and the method of proof has appeared
in the end of Chapter 6 of \cite{MasterEquation}.
We, however, include a proof for the sake of self-containedness.

First let us show the bound on the Wasserstein distance
$W_2(m^N_t, m_\infty^{\otimes N})$.
Recall that $\mathbf X_t = (X^i_t)_{i=1}^N$ is the solution
of the SDE \cref{eq:ps-sde}
with Brownian motions $(W^i)_{i=1}^N$.
Let $\tilde{\mathbf X}^i_t = (\tilde X^i_t)_{i=1}^N$ solve
\[
\dd\tilde X^i_t = - D_m F (m_t, \tilde X^i_t) \dd t + \sqrt 2 \dd W^i_t,
\quad i = 1,\ldots, N
\]
with the initial condition
$\Law (\tilde X^1_0, \ldots, \tilde X^N_t) = m_0^{\otimes N}$ and
\[
W_2^2 (m^N_0, m_t^{\otimes N})
= \sum_{i=1}^N \Expect \bigl[ \lvert X^i_0 - \tilde X^i_t \rvert^2 \bigr],
\]
i.e., the couple $(\mathbf X^i_0, \tilde{\mathbf X}^i_0)$ is distributed
as the $L^2$-optimal transport plan between $m^N_0$ and $m_0^{\otimes N}$.
Then, by subtracting the dynamical equations of $\mathbf X_t$
and $\tilde{\mathbf X}_t$, we have
\begin{align*}
\dd \biggl( \sum_{i=1}^N \lvert X^i_t - \tilde X^i_t \rvert^2 \biggr)
&= - 2 \sum_{i=1}^N (X^i_t - \tilde X^i_t)
\cdot \bigl( D_m F (\mu_{\mathbf X_t}, X^i_t)
- D_m F(m_t, \tilde X^i_t) \bigr) \\
&\leq \sum_{i=1}^N \lvert X^i_t - \tilde X^i_t \rvert^2
+ \sum_{i=1}^N \lvert D_m F (\mu_{\mathbf X_t}, X^i_t)
- D_m F(m_t, \tilde X^i_t) \rvert^2,
\end{align*}
where the difference between the drifts satisfies
\begin{align*}
&\hspace{-1em}\lvert D_m F(\mu_{\mathbf X_t}, X^i_t)
- D_m F(m_t, \tilde X^i_t)\rvert\\
&\leq \lvert D_m F(\mu_{\mathbf X_t}, X^i_t)
- D_m F(\mu_{\tilde{\mathbf X}_t},\tilde X^i_t)\rvert
+ \lvert D_m F(\mu_{\tilde{\mathbf X}_t},\tilde X^i_t)
- D_m F(m_t, \tilde X^i_t)\rvert\\
&\leq M^F_{mm} W_1 (\mu_{\mathbf X_t}, \mu_{\tilde{\mathbf X}_t})
+ M^F_{mx} \lvert X^i_t - \tilde X^i_t \rvert
+ M^F_{mm} W_1 (\mu_{\tilde{\mathbf X}_t}, m_t).
\end{align*}
Thus, we have
\begin{multline}
\label{eq:w-poc-growth}
\frac{\dd}{\dd t}
\biggl( \sum_{i=1}^N \lvert X^i_t - \tilde X^i_t \rvert^2 \biggr)
\leq \bigl( 1 + 3 (M^F_{mx})^2 \bigr)
\sum_{i=1}^N \lvert X^i_t - \tilde X^i_t \rvert^2
+ 3 N (M^F_{mm})^2 W_2^2 (\mu_{\mathbf X_t}, \mu_{\tilde{\mathbf X}_t}) \\
+ 3 N (M^F_{mm})^2 W_2^2 (\mu_{\tilde{\mathbf X}_t}, m_t).
\end{multline}
For the second term, we have
\[
\Expect \bigl[ W_2^2 (\mu_{{\mathbf X}_t}, \mu_{\tilde {\mathbf X}_t}) \bigr]
\leq \frac 1N \sum_{i=1}^N \Expect
\bigl[ \lvert X^i_t - \tilde X^i_t\rvert^2\bigr],
\]
and for the last term, we have, by the result of Fournier and Guillin
\cite{FournierGuillin},
\begin{align*}
\Expect\bigl[ W_2^2 (\mu_{\tilde{\mathbf X}_t}, m_t) \bigr]
&\leq C(d) \Expect\bigl[|X_t - \Expect X_t|^6\bigr]^{1/3} \delta_d(N) \\
&= C(d) \Expect\bigl[|X_t - \Expect X_t|^6\bigr]^{1/3} \times
\begin{cases}
N^{-1/2} & \text{if}~d < 4, \\
N^{-1/2} \log (1+N) & \text{if}~d = 4, \\
N^{-2/d} & \text{if}~d > 4.
\end{cases}
\end{align*}
Then,
denoting $\tilde X_t = \tilde X^1_t$,
we only need to control
$\Expect\bigl[|\tilde X_t - \Expect \tilde X_t|^6\bigr]$.
Observe that, by It\=o's formula, we have
\begin{multline*}
\frac{\dd}{\dd t}\Expect [ \lvert \tilde X_t - \Expect \tilde X_t \rvert^6 ] \\
= - 6 \Expect \bigl[ |\tilde X_t - \Expect \tilde X_t|^4
(\tilde X_t - \Expect \tilde X_t) \cdot
\bigl(D_m F(m_t, \tilde X_t) - \Expect[D_m F(m_t, \tilde X_t)]\bigr) \bigr] \\
+ (6d + 24) \Expect \bigl[ |\tilde X_t - \Expect \tilde X_t|^4 \bigr].
\end{multline*}
Then we have the following control of the growth,
by using the elementary inequality
$x^4 \leq \frac 23 x^6 + \frac 13$ for $x \geq 0$:
\[
\frac{\dd}{\dd t}
\Expect\bigl[ \lvert \tilde X_t - \Expect \tilde X_t \rvert^6\bigr]
\leq (6 M^F_{mx} + 4d + 16) \Expect [ \lvert \tilde X_t
- \Expect \tilde X_t \rvert^6 ]
+ (2d + 8).
\]
Thus, by Grönwall's lemma, we have
\begin{multline*}
\Expect\bigl[|\tilde X_t - \Expect \tilde X_t|^6\bigr]
\leq e^{(6M^F_{mx} + 4d + 16)t}
\Expect\bigl[ |\tilde X_0 - \Expect \tilde X_0|^6\bigr] \\
+ \frac{d + 4}{3M^F_{mx} + 2d + 8}
(e^{(6 M^F_{mx} + 4d + 16)t} - 1).
\end{multline*}
We take expectations on both side of the differential inequality
\cref{eq:w-poc-growth} and obtain
\begin{multline*}
\frac{\dd}{\dd t}
\Expect \biggl[ \sum_{i=1}^N \lvert X^i_t - \tilde X^i_t \rvert^2 \biggr]
\leq \bigl( 1 + 3(M^F_{mx})^2 + 3(M^F_{mm})^2 \bigr)
\Expect \biggl[ \sum_{i=1}^N \lvert X^i_t - \tilde X^i_t \rvert^2 \biggr] \\
+ 3N (M^F_{mm})^2 C(d) \delta_d(N)
\Expect \bigl[\lvert \tilde X_t - \Expect \tilde X_t\rvert^6\bigr]^{1/3}.
\end{multline*}
We then use Grönwall's lemma to show \cref{eq:w-poc-finite}.

As for the distance under relative entropy, by Girsanov's theorem we have
\[
H(m^N_t | m_t^{\otimes N})
\leq H(m^N_0 | m_0^{\otimes N})
+ \frac 14 \sum_{i=1}^N \int_0^t \Expect\bigl[
\lvert D_m F(\mu_{{\mathbf X}_s}, X^i_s)
- D_m F(m_s, X^i_s) \rvert^2 \bigr]\dd s,
\]
and we can control the last term by
\begin{align*}
\lvert D_m F(\mu_{{\mathbf X}_s}, X^i_s)
- D_m F(m_s, X^i_s) \rvert
&\leq M^F_{mm} W_2 (\mu_{\mathbf X_s}, m_s) \\
&\leq M^F_{mm} \bigl( W_2(\mu_{\mathbf X_s,}, \mu_{\tilde{\mathbf X}_s})
+ W_2 (\mu_{\tilde{\mathbf X}_s}, m_s) \bigr).
\end{align*}
So we can show \cref{eq:h-poc-finite} by using the same method as before.

\paragraph{Long-time behavior.}
The triangle inequality for the \(L^2\)-Wasserstein distance gives us
\(W_2^2(m_t^N, m_t^{\otimes N})
\leq 2\bigl(W_2^2(m_t^N, m_\infty^{\otimes N})
+ W_2^2(m_t^{\otimes N}, m_\infty^{\otimes N})\bigr)\).
By Talagrand's inequality \cref{eq:t2} for \(m_\infty^{\otimes N}\) we bound
the Wasserstein distances by
\begin{gather*}
\rho W_2^2(m_t^N, m^{\otimes N}_\infty)
\leq H(m_t^N | m_\infty^{\otimes N})
\leq \mathcal F^N(m^N_t) - N\mathcal F(m_\infty), \\
\rho W_2^2(m_t^{\otimes N}, m^{\otimes N}_\infty)
= NW_2^2(m_t, m_\infty)
\leq N H(m_t^N | m_\infty)
\leq N \bigl(\mathcal F(m_t) - \mathcal F(m_\infty)\bigr),
\end{gather*}
where we applied \cref{lem:mf-entropy-sandwich,lem:ps-entropy-sandwich}.
We then apply \cref{thm:mf-entropy-convergence,thm:ps-entropy-convergence} to
obtain \cref{eq:w-poc}.

Now suppose additionally \cref{eq:higher-der} and \(h_0 = m_0 / m_\infty \in
L^{p_0} (m_\infty)\) for \(p_0 > 1\).
The relative entropy satisfies
\begin{align*}
H(m^N_t | m^{\otimes N}_t)
&= \int \log \frac{m^N_t(\mathbf x)}{m^{\otimes N}_t(\mathbf x)}
m^N_t (\mathbf x)\dd\mathbf x \\
&= \int \biggl( \log \frac{m^N_t(\mathbf x)}{m^{\otimes N}_\infty(\mathbf x)}
- \log \frac{m^{\otimes N}_t(\mathbf x)}{m^{\otimes N}_\infty(\mathbf x)}
\biggr) m^N_t(\mathbf x)\dd\mathbf x \\
&= H(m^N_t | m_\infty^{\otimes N}) - \sum_{i=1}^N \int \log
\frac{m_t (x)}{m_\infty(x)} m^{N,i}_t(x)\dd x,
\end{align*}
where \(m^{N,i}_t\) is the \(i\)-th marginal of \(m^N_t\).
We then apply \cref{eq:lp-bounds-ce} in \cref{lem:lp-bounds-ce} to summands in
the second term with \(p = 1\) to obtain
\[
- \int \log \frac {m_t(x)}{m_\infty(x)} m^{N,i}_t (x)\dd x
\leq H(m^{N,i}_t | m_\infty) - \log \lVert h_t \rVert_{-1}.
\]
So we have
\begin{align*}
- \sum_{i=1}^N \int \log \frac{m_t(x)}{m_\infty(x)} m^{N,i}_t(x)\dd x
&\leq - N \log \lVert h_t \rVert_{-1} + \sum_{i=1}^N H(m^{N,i}_t | m_\infty) \\
&\leq - N \log \lVert h_t \rVert_{-1} + H(m^N_t | m_\infty^{\otimes N}),
\end{align*}
where we used the information inequality \cref{eq:info-ineqs} in the last
inequality.
Therefore
\[
H(m^N_t | m^{\otimes N}_t) \leq -N \log \lVert h_t \rVert_{-1} + 2 H(m^N_t |
m^{\otimes N}_\infty).
\]
We conclude by applying the results of
\cref{thm:lp-convergence,thm:ps-entropy-convergence}.
\end{proof}
\vspace{2mm}

\begin{proof}[Proof of \cref{cor:poc}]
In the Wasserstein case,
let $C_4$, $C_5$ be the constants in \cref{thm:poc}.
We take $t_0 = \log N / (d\vee 4) C_4$.
Then, for $t \leq t_0$, by using \cref{eq:w-poc-finite}, we have
\begin{multline}
\label{eq:w-bound-short}
\frac 1N W_2^2 (m^N_t, m_t^{\otimes N})
\leq C_5 (e^{C_4 t} - 1) \bigl( v_6(m_0)^{1/3} + 1 \bigr) \delta_d(N) \\
\leq C_5 (N^{1/(d \vee 4)} - 1)\bigl( v_6(m_0)^{1/3} + 1 \bigr) \delta_d(N),
\end{multline}
where $N^{1/(d\vee 4)} \delta_d (N) \leq N^{-1/(d\vee 4)} \log (1 + N)$
for all $d$.
For $t \geq t_0$, by using \cref{eq:w-poc}, we have
\begin{multline}
\label{eq:w-bound-long}
\frac 1N W_2^2 (m^N_t, m_t^{\otimes N})
\leq 2 \bigl(\mathcal F(m_0) - \mathcal F(m_\infty)\bigr)
N^{- 4\rho/(d\vee 4)C_4} \\
+ \frac{2}{N}\bigl( \mathcal F^N(m_0^{\otimes N}) - N \mathcal F (m_\infty)\bigr)
N^{-(4\rho' - C_1 N^{-1})/(d \vee 4)C_4} \\
+ \frac{2C_2}{4N\rho' - C_1},
\end{multline}
if $N > C_1 / 4\rho'$, where $\rho' \in (0,\rho)$ and $C_1$, $C_2$
are defined in \cref{thm:ps-entropy-convergence}.
By expanding the functional $F$, we also have
\begin{multline*}
F ( \mu_{\mathbf X_0} ) - F(m_0)
= \int \frac{\delta F}{\delta m} (m_0, x) (\mu_{\mathbf X_0} - m_0) (\dd x) \\
+ \int_0^1 \biggl( \frac{\delta F}{\delta m}
\bigl((1-t) \mu_{\mathbf X_0} + tm_0, x\bigr)
- \frac{\delta F}{\delta m}(m_0, x) \biggr) (\mu_{\mathbf X_0} - m_0) (\dd x) dt
\end{multline*}
with
\[
\Expect\biggl[
\int \frac{\delta F}{\delta m} (m_0, x) (\mu_{\mathbf X_0} - m_0) (\dd x)
\biggr] = 0
\]
and
\begin{align*}
&\hspace{-1em}
\Expect\biggl[\int_0^1 \biggl( \frac{\delta F}{\delta m}
\bigl((1-t) \mu_{\mathbf X_0} + tm_0, x\bigr)
- \frac{\delta F}{\delta m}(m_0, x) \biggr)
(\mu_{\mathbf X_0} - m_0)(\dd x)\dd t
\biggr] \\
&\leq \Expect \biggl[ \int_0^1
\bigl\lVert D_m F
\bigl((1-t) \mu_{\mathbf X_0} + tm_0, \cdot\bigr)
- D_m F(m_0, \cdot)\bigr\rVert_\infty W_1(\mu_{\mathbf X_0}, m_0)\dd t \biggr]
\\
&\leq \frac{M^F_{mm}}{2} \Expect\bigl[W_2^2(\mu_{\mathbf X_0}, m_0)\bigr]
\leq M^F_{mm} \Var m_0.
\end{align*}
Thus, we obtain
\begin{multline}
\label{eq:free-energy-initial-bound}
\mathcal F^N(m_0^{\otimes N})
= N \Expect \bigl[ F(\mu_{\mathbf X_0}) \bigr] + H(m_0^{\otimes N})
\leq N F(m_0) + N M^F_{mm} \Var m_0 + N H(m_0) \\
= N \mathcal F(m_0) + N M^F_{mm} \Var m_0.
\end{multline}
Taking $\rho' = \rho/2$,
we obtain the uniform-in-time Wasserstein bound
\cref{eq:w-poc-unif} from \cref{eq:w-bound-short,eq:w-bound-long}.

Similarly, to control the relative entropy,
we take $t'_0 = \tau + \frac{\log N}{(d \vee 4)C_4}$, where
$\tau$ is the constant in \cref{thm:poc}.
So, for $t \leq t'_0$, by \cref{eq:h-poc-finite}, we have
\begin{equation}
\label{eq:h-bound-short}
\frac 1N H(m^N_t | m_t^{\otimes N})
\leq C_5 (e^{C_4 \tau} N^{1/(d\vee 4)} - 1)
\bigl( v_6(m_0)^{1/3} + 1 \bigr) \delta_d(N),
\end{equation}
and, for $t \geq t'_0$, by \cref{eq:h-poc}, we have
\begin{multline}
\label{eq:h-bound-long}
\frac 1N H(m^N_t | m_t^{\otimes N})
\leq C_3 e^{-4\rho'\tau} N^{- 4\rho'/(d \vee 4)} \\
+ \frac 2N \bigl( \mathcal F^N(m_0^{\otimes N}) - N\mathcal F(m_\infty)\bigr)
e^{-(4\rho' - C_1N^{-1}) \tau} N^{- (4\rho' - C_1N^{-1})/ (d\vee 4)C_4} \\
+ \frac{2C_2}{4N\rho' - C_1}.
\end{multline}
So, using again \cref{eq:free-energy-initial-bound},
we can combine \cref{eq:h-bound-short,eq:h-bound-long} to obtain
the uniform-in-time entropic bound \cref{eq:h-poc-unif}.
\end{proof}

\appendix

\section{Proofs of technical results on MFL}
\label{sec:mf-technical}

In the section we provide proofs of technical results on the regularity
properties of the MFL dynamics.

\begin{proof}[Proof of \cref{prop:mf-exist-unique-regular}]
It is classical that under the conditions \cref{eq:lip-in-m,eq:first-der} the
McKean--Vlasov SDE
\[
\dd X_t = - D_m F(m_t, X_t)\dd t + \sqrt 2 \dd W_t, \qquad \Law(X_t) = m_t
\]
has unique global solution defined for \(t \in [0,+\infty)\).
By construction the marginal law \(m_t = \Law(X_t)\) is
in \(C\bigl([0,+\infty); \mathcal P_2(\mathbb R^d)\bigr)\),
proving the existence of solution.
Any solution to the Fokker--Planck equation admits equally this probabilistic
representation,
then the uniqueness in short time follows from Cauchy--Lipschitz bounds.
We extend this uniqueness to the infinity by sewing up the short time intervals,
finishing the proof of the first claim.

Let \(\rho_{t}(x)\) be the density of Gaussian \(\mathcal N(0, 2t)\).
The solution \(m_t\) satisfies Duhamel's formula in the sense of distributions
\begin{align*}
m_t &= \rho_t \star m_0 + \int_0^t \rho_{t-s} \star \nabla
\cdot\bigl(m_s D_m F(m_s, \cdot)\bigr)\dd s \\
&= \rho_t \star m_0
+ \sum_{i=1}^d \int_0^t \nabla_i \rho_{t-s}
\star\bigl(m_s D_m F^i(m_s,\cdot)\bigr)\dd s.
\end{align*}
Note that \(\Vert \nabla \rho_t \Vert_{L^p(\mathbb R^d)} \leq C_{d,p} t^{-\frac
12 + \frac d2(\frac 1p - 1)}\), which is integrable around \(0+\) when \(p <
\frac{d}{d-1}\).
In this case apply Young's convolution inequality to obtain
\[
\Vert m_t \Vert_{L^p(\mathbb R^d)} \leq \Vert \rho_t \Vert_{L^p(\mathbb R^d)}
\Vert m_0 \Vert_{\textnormal{TV}} + \sum_{i=1}^d \int_0^t \Vert \nabla_i
\rho_{t-s} \Vert_{L^p(\mathbb R^d)} \Vert m_s D_m F^i(m_s,\cdot)
\Vert_{\textnormal{TV}}\dd s,
\]
where \(\sup_{s \in [0,t]}\Vert m_s D_m F^i(m_s, \cdot) \Vert_\textnormal{TV}
\leq \sup_{s \in [0,t]} C \int (1 + |x|) m_s(\dd x) < +\infty\).
Hence \(\Vert m_t \Vert_{L^p(\mathbb R^d)} < +\infty\) for all \(t > 0\).
This and the second moment bound
\(\int |x|^2 m_t(\dd x) < +\infty\)
are sufficient for the finiteness of entropy, i.e. the integral \(\int |\log
m_t(x)| m_t(x) \dd x\) is finite, which is our second claim.
Indeed for the lower bound on entropy we use the decomposition in
\cref{eq:entropy-decomposition}, while the upper bounds follows from \(m \log m
\leq \frac{m^p - m}{p-1}\) for all \(p > 1\).

The drift \(D_m F(m_t, x)\) has uniform linear growth in \(x\):
\[
|D_m F(m_s,x)| \leq M^F_{mx} |x| + \sup_{s \in [t_0,t]} |D_m F(m_s,0)|,
\]
where \(M^F_{mx}\) is the constant in \cref{eq:first-der}
and the second term is finite by the compactness of set
\(\{m_s : s \in [t_0,t]\}\) in \(\mathcal P_2\).
As a result,
\[\int_{t_0}^t \int |D_m F(m_s, x)|^2 m_s(\dd x)\dd t < +\infty.\]
We then apply \cite[Theorem 7.4.1]{fpkeq} to obtain the finiteness of
\cref{eq:integrated-fisher}.
Especially,
\(\nabla m \in L^1_\textnormal{loc}\bigl((0,+\infty); L^1(\mathbb R^d)\bigr)\).
Rewrite the Fokker--Planck equations as a continuity equation \(\partial_t m +
\nabla \cdot (m_t v_t) = 0\) where \(v_t(x) = - D_m F(m_t,x) - \nabla \log
m_t(x)\).
We have
\begin{multline*}
\int_{t_0}^t \int |v_s(x)|^2 m_s(\dd x) \dd s \\
\leq 2 \biggl(
\int_{t_0}^t \int |D_m F(m_s,x)|^2 m_s(\dd x)\dd s
+ \int_{t_0}^t \int \frac{|\nabla m_s(x)|^2}{m_s(x)}\dd x \dd s
\biggr) < +\infty.
\end{multline*}
Hence by \cite[Theorem 8.3.1]{gf} the flow \(m_t\) is locally \(AC^2\) in
\((\mathcal P_2, W_2)\).
The vector field \(v_t (x) = - D_m F(m_t,x) - \nabla \log m_t(x)\) solves the
continuity equation
\begin{equation}
\label{eq:ce}
\partial_t m_t + \nabla \cdot (m_t v_t) = 0
\end{equation}
in the sense of distributions and \(v_t\) writes in the gradient form
\(v_t = - \nabla\bigl(\frac{\delta F}{\delta m} (m_t,x) + \log m_t(x)\bigr)
= - \nabla \varphi_t\).

We finally verify \(v_t\) is indeed a tangent vector of \(m_t\) according to
\cite[Definition 8.4.1]{gf},
i.e. \(v_t \in \Tan_{m_t} \mathcal P_2(\mathbb R^d)= \overline{\{ \nabla
\varphi : \varphi \in C^\infty_c (\mathbb R^d) \}}^{L^2(m_t)}\).
Let \(\eta_R : \mathbb R^d \to [0,1]\) be a smooth function supported on
\(B(2R)\), has the constant value \(1\) on \(B(R)\) and satisfies \(|\nabla
\eta(x)| \leq 2/R\) for all \(x\).
We have
\[
\int |\nabla\varphi_t - \nabla(\varphi_t\eta_R)|^2 m_t
\leq 2 \int_{B(2R) \setminus B(R)}
\bigl(|\varphi_t|^2 |\nabla \eta_R|^2
+ |\nabla \varphi_t|^2 |1 - \eta_R|^2\bigr) m_t.
\]
The second term tends to \(0\) when \(R \to \infty\), while the first satisfies
\begin{align*}
&\hspace{-1em}\int_{B(2R)\setminus B(R)} |\varphi_t|^2 |\nabla\eta_R|^2 m_t \\
&\leq \frac{2}{R^2} \int_{B(2R)\setminus B(R)}
\biggl( \biggl| \frac{\delta F}{\delta m}(m_t,x)\biggr|^2
+ \lvert\log m_t(x)\rvert^2 \biggr) m_t \\
&\leq \frac{2C}{R^2} \int_{B(2R)\setminus B(R)} (1 + |x|^4) m_t(\dd x)
+ \frac{2}{R^2} \int_{B(2R) \setminus B(R)} \lvert\log m_t\rvert^2 m_t \\
&\leq \frac{2C}{R^2} \int_{B(2R)\setminus B(R)} (1 + 4R^2|x|^2) m_t(\dd x)
+ \frac{2}{R^2} \int_{B(2R) \setminus B(R)} \lvert\log m_t\rvert^2 m_t.
\end{align*}
Here the first term tends to \(0\) since \(m_t \in \mathcal P_2\),
while the second term tends to \(0\) by the integrability of \(|\log m_t|^2
m_t\),
which follows from the elementary inequality
\[
m \lvert\log m\rvert^2 \leq C_p m^p \mathbf 1_{m \geq 1}
+ 2\Bigl(|x|^2 m + \sup_{t \in [0,1]} t (\log t)^2 e^{-|x|}\Bigr)
\mathbf 1_{m < 1}
\]
for \(p > 1\) and \(x \in \mathbb R^d\).
Hence \(\nabla (\varphi_t \eta_R) \to \nabla \varphi_t\) in \(L^2(m_t)\).
It then suffices to approximate the (essentially) compactly supported function
\(\varphi_t \eta_R\) by \(C^\infty_c\) functions
in the \(L^2(m_t)\)-norm.
We can do this by taking a sequence of compacted supported mollifiers \(\rho_n\)
and applying them to obtain
\(\nabla (\varphi_t \eta_R) \star \rho_n \to \nabla (\varphi_t \eta_R)\)
in \(L^2(m_t)\) when \(n \to \infty\).
\end{proof}

\begin{proof}[Proof of \cref{prop:std-alg-density}]
Let \(h\) be a positive function.
Define the functions \(k_n = \mathbf 1_{B(n)} (h \wedge n) \vee 1/n\) and
\(k_{n,m} = \rho_{m} \star k_n \),
where \((\rho_{m})_{m \in \mathbb N}\) is a sequence of \(C^\infty\) mollifiers.
They satisfy
\[
\forall x \in \mathbb R^d,\qquad
\frac 1n \leq k_n(x),k_{n,m}(x) \leq n
\text{ and }
|\nabla^\ell k_{n,m} (x)| \leq n \Vert \nabla^\ell \rho_m \Vert_\infty <
+\infty.
\]
In particular \(k_{n,m} \in \mathcal A_+\).
We have \(k_n \to h\) in \(L^p(\mu)\) whenever \(h \in L^p(\mu)\) for \(p \geq
1\)
and \(\Vert k_n \Vert_q \to \Vert h \Vert_q\) whenever \(h \in L^q(\mu)\) for
\(q \leq 1\)
by the dominated convergence theorem.
Since for all \(n \in \mathbb N\) the function \(k_n \in L^1(\mathbb R^d)\), we
have
\(k_{n,m} \to k_n\) in \(L^1(\mathbb R^d)\) when \(m \to \infty\).
Hence \(k_{n,m} \to k_n\) a.e. when \(m \to \infty\) along a subsequence.
Then we can apply again the dominated convergence to obtain
\(k_{n,m} \to k_n\) in \(L^p(\mu)\) for all \(p \geq 1\)
and \(\Vert k_{n,m} \Vert_q \to \Vert k_n \Vert_q\) for all \(q < 1\).
We can thus taking a subsequence of \((n,m) \to (+\infty,+\infty)\)
so that \(k_{n,m} \to h\) in the desired ways.
\end{proof}

\begin{proof}[Proof of \cref{prop:std-alg-stability}]
Fix \(T > t_0\).
We denote by \(C\) a positive constant that depends on
\(\max_{k=1,2,3}\sup_{m,x} |\nabla^k D_m F (m, x) |\)
and on the initial condition \(h' \in \mathcal A_+\);
and by \(C_Q\) a positive constant that depends additionally on the quantity
\(Q\).
The constants $C$, $C_Q$ may change from line to line.
Define \(g(t,x) = \nabla \cdot (b_t - b_\infty) + (b_t - b_\infty) \cdot
b_\infty\).
It satisfies \(|g(t,x)| \leq C(1+|x|)\)
for all \((t,x) \in [t_0,T] \times \mathbb R^d\) as
\(\Vert \nabla^k (b_t - b_\infty ) \Vert_\infty \leq C\)
for $k = 0$, $1$ and \(t \in [t_0, T]\).
Fix \(t \in [t_0,T]\).
Let \((X^{t,x}_s)_{s\in [0,t-t_0]}\) be the stochastic process solving
\begin{equation}
\label{eq:std-alg-aux-sde}
\dd X^{t,x}_s = (2b_\infty - b_{t - s}) \dd s + \sqrt 2 \dd W_s
\end{equation}
with \(X^{t,x}_0 = x\)
and define as well its extremal process \(M^{t,x}_s = \sup_{0 \leq u \leq s}
|X_u|\) for \(s \in [0,t - t_0]\).
Since the drift satisfies \((2b_\infty - b_t) \cdot x \leq C_T |x|^2 + C_T\)
for all \((t,x) \in [t_0,T] \times \mathbb R^d\),
we obtain the Gaussian moment bound
\[\Expect \exp\bigl(C_T^{-1} |M^{t,x}_{t-t_0}|^2\bigr)
\leq C_T \exp (C_T|x|^2)\]
by Itō's formula and Doob's maximal inequality.
As a consequence the exponential moments are finite:
\[
\forall \alpha \geq 0,\qquad
\Expect \exp\bigl(\alpha |M^{t,x}_{t-t_0}|\bigr)
\leq C_{T,\alpha} \exp (C_{T,\alpha}|x|).
\]
Set \(h(t_0, \cdot) = h'\).
We construct the solution by the Feynman--Kac formula for \cref{eq:h}
\[
h(t,x) \coloneqq \Expect \biggl[
\exp \biggl( - \int_0^{t-t_0} g(t - s, X^{t,x}_s)\dd s \biggr)
h(t_0, X^{t,x}_{t-t_0}) \biggr].
\]
It is standard that the \(h\) constructed above solves \cref{eq:h} in the sense
of distributions.
We verify \(h_t \in \mathcal A_+\) for all \(t \in [t_0, T]\).
For the upper bound we apply the Cauchy--Schwarz inequality to obtain
\begin{align*}
h(t,x) &\leq \Expect
\biggl[\exp \biggl( - 2\int_0^{t-t_0} g(t - s, X^{t,x}_s)\dd s\biggr)\biggr]^{1/2}
\Expect\bigl[ h(t_0,X^{t,x}_{t-t_0})^2\bigr]^{1/2} \\
&\leq \Expect\bigl[ \exp\bigl(C_T (1 + |M^{t,x}_{t-t_0}|)\bigr)\bigr]^{1/2}
\Expect \bigl[ \exp \bigl(C_T(1 + |X^{t,x}_{t-t_0}|)\bigr)\bigr]^{1/2} \\
&\leq \Expect\bigl[ \exp\bigl(C_T (1 + |M^{t,x}_{t-t_0}| )\bigr)\bigr]
\leq \exp\bigl(C_T(1 + |x|)\bigr).
\end{align*}
We applied the bound on \(g\) and \(h\) in the second inequality and
used the exponential moment bound on \(M_{t-t_0}\) in the last.
For the lower bound we use Cauchy--Schwarz from the other direction:
\begin{align*}
h(t,x)
&\geq \Expect \biggl[ \exp \biggl( \int_0^{t-t_0} g(t-s,X^{t,x}_s) \dd s
\biggr) \biggr]^{-1}
\Expect\bigl[ h(t_0, X^{t,x}_{t-t_0})^{1/2}\bigr]^2 \\
&\geq C^{-1}_T \Expect\bigl[ \exp\bigl( C_T |M^{t,x}_{t-t_0}|\bigr)\bigr]^{-1}
\Expect\bigl[\exp\bigl(- C_T|X^{t,x}_{t-t_0}|\bigr)\bigr]^2 \\
&\geq C^{-1}_T \Expect\bigl[ \exp\bigl( C_T |M^{t,x}_{t-t_0}|\bigr)\bigr]^{-1}
\Expect\bigl[\exp\bigl(C_T|X^{t,x}_{t-t_0}|\bigr)\bigr]^{-2} \\
&\geq C^{-1}_T \Expect\bigl[ \exp\bigl( C_T |M^{t,x}_{t-t_0}|\bigr)\bigr]^{-3}
\geq C^{-1}_T \exp (- C_T|x|).
\end{align*}
Again we applied the bound on \(g\) and \(h\) on the second inequality
and used the exponential moment bound on \(M_{t-t_0}\) on the last line.
So we have proved the bound of both sides
\(\lvert\log h(t,x)\rvert \leq C_T(1+|x|)\),
that is, the ``zeroth-order'' condition of \(\mathcal A_+\).

Now derive the continuity of \(x \mapsto h(t,x)\).
Let the stochastic processes \((X^{t,x}_\cdot)_{x \in \mathbb R^d}\) be coupled
by sharing the same Brownian motion in their defining SDEs
\cref{eq:std-alg-aux-sde}.
The mapping \(x \mapsto X^{t,x}_s\) is continuous almost surely as
its matrix-valued partial derivative
$\partial X^{t,x}_\cdot/\partial x$ solves the SDE
\[
\dd\frac{\partial X^{t,x}_s}{\partial x}
= \nabla\bigl(2b_\infty (X^{t,x}_s) - b_{t-s}(X^{t,x}_s)\bigr)
\frac{\partial X^{t,x}_s}{\partial x} \dd s
\]
whose wellposedness is guaranteed by the bound
\[
|\nabla^2 (2b_\infty - b_{t-s})(x)| \leq 3 \sup_{m \in \mathcal P_2(\mathbb
R^d)} \sup_{x \in \mathbb R} | \nabla^2 D_m F(m,x) | \leq C.
\]
The norm of \(\frac{\partial X^{t,x}_s}{\partial x}\)
satisfies
\[
\forall s \in [0,t-t_0],~\forall x \in \mathbb R^d,
\qquad\biggl\lvert\frac{\partial X^{t,x}_s}{\partial x}\biggr\rvert\leq C_T
\quad\text{a.s.}
\]
by Grönwall's lemma.
Therefore we have
\[
\Expect\Bigl[ \exp\Bigl( C_T^{-1} \sup_{x : |x - x_0| \leq 1}
|M^{t,x}_{t-t_0}|^2\Bigr)\Bigr]
\leq C_T \exp(C_T|x_0|^2)
\]
for all \(x_0 \in \mathbb R^d\).
We obtain \(h(t, x) \to h(t,x_0)\)
when \(x \to x_0\)
by applying the dominated convergence theorem to the Feynman--Kac formula.

We sketch the part for verifying the conditions on derivatives.
Differentiate the evolution equation \cref{eq:h}.
We obtain for $k = 1$, $2$,
\begin{align*}
\partial_t \nabla^k h &= \Delta \nabla^k h + (2 b_\infty - b_t) \cdot \nabla
\nabla^k h + \sum_{i=2}^{k} \binom{k}{i} \nabla^i (2b_\infty - b_t) \cdot
\nabla \nabla^{k-i} h \\
&\quad+ \sum_{i=1}^k \binom{k}{i} \nabla^i g(t,x) \nabla^{k-i} h
+ \bigl(\nabla (2b_\infty - b_t) \cdot \nabla \nabla^{k-1} h
+ g(t,x) \nabla^k h\bigr).
\end{align*}
We then write the Feynman--Kac formula for \(\nabla^k h\), $k = 1$, $2$.
The first two terms on the right hand side of the equation corresponds to the
same stochastic process, to which the Gaussian moment bound applies.
The third and fourth term are lower-order derivatives,
continuous in space and have bound
\(|\nabla^{k-i} h(t,x)| \leq \exp \bigl(C_T(1 + |x|)\bigr)\)
by the induction hypothesis.
The last term corresponds to the exponential in the Feynman--Kac formula, whose
growth in \(x\) remains linear.
So we can argue as before to derive
\(|\nabla^k h(t,x)| \leq \exp\bigl(C_T(1+|x|)\bigr)\)
for all \((t,x) \in [t_0,T]\times\mathbb R^d\).
The continuity of \(x \mapsto \nabla^k h(t,x)\)
for $k=1$, $2$ follows analogously.
Since \(x \mapsto h(t,x)\) are twice-differentiable the generalized derivative
\(\partial_t h\) exists by the evolution equation \cref{eq:h}.
Finally all the constants in the bounds depend only additionally on \(T\),
so \((h_t)_{t \in [t_0,T]} \subset \mathcal A_+\) uniformly.
\end{proof}

\section{Proof of modified Bochner's theorem}
\label{sec:bochner}

\begin{proof}[Proof of \cref{thm:bochner}]
We prove the theorem by showing
$\text{(i)} \Rightarrow \text{(ii)} \Rightarrow \text{(iii)} \Rightarrow
\text{(i)}$.

\paragraph{$\text{(i)} \Rightarrow \text{(ii)}$.}
Suppose (i) holds, i.e., \(m \mapsto F_\textnormal{Int}(m)\) is convex.
Let $\mu$ be a compactly supported signed measure with $\int \dd\mu = 0$.
Then it admits decomposition into positive and negative parts:
$\mu = \mu_+ - \mu_-$.
We define the probability measure
\[
m \coloneqq \frac{\lvert \mu \rvert}{\bigl\lVert \lvert\mu\rvert \bigr\rVert%
_\textnormal{TV}}
= \frac{\mu_+ + \mu_-}{\lVert\mu_+\rVert_\textnormal{TV}
+ \lVert \mu_- \rVert_\textnormal{TV}}.
\]
Then, for all
$t < (\lVert\mu_+\rVert_\textnormal{TV} + \lVert \mu_-
\rVert_\textnormal{TV})^{-1}
\eqqcolon t_0$,
we have $m_t \coloneqq m + t \mu \in \mathcal P(\mathbb R^d)$.
Thus, the mapping
\[
t \mapsto F_\textnormal{Int} (m_t)
= F_\textnormal{Int} (m) + t \iint V(x - y) m(\dd x) \mu(\dd y)
+ \frac {t^2}2 \iint V(x - y) \mu(\dd x) \mu(\dd y)
\]
is convex on the interval $(-t_0,t_0)$, and therefore,
$\iint V(x - y) \mu(\dd x) \mu(\dd y) \geq 0$,
which proves (ii).

\paragraph{$\text{(ii)} \Rightarrow \text{(iii)}$.}
Suppose (ii) holds.
For non-zero $s \in \mathbb R^d$, we define the bounded and continuous function
$W_s(t) \coloneqq 2V(t) - V(t + s) - V(t - s)$.
Then, for every $\boldsymbol\xi \in \mathbb R^N$
and every $x^1,\ldots,x^N \in \mathbb R^d$, we have
\begin{align*}
&\hspace{-1em}
\sum_{i,j=1}^N \xi^i \xi^j W_s(x^i - x^j) \\
&= \sum_{i,j=1}^N \xi^i \xi^j V(x^i - x^j)
+ \sum_{i,j=1}^N \xi^i \xi^j V\bigl( (x^i+s - (x^j+s) \bigr) \\
&\phantom{={}}\quad - \sum_{i,j=1}^N \xi^i \xi^j V\bigl( (x^i+s) - x^j\bigr)
- \sum_{i,j=1}^N \xi^i \xi^j V\bigl(x^i - (x^j + s)\bigr) \\
&= \iint V(x - y) \hat\mu(\dd x) \hat\mu(\dd y) \geq 0,
\qquad\text{for}~
\hat\mu = \sum_{i=1}^N \xi^i \delta_{x^i} - \sum_{i=1}^N \xi^i
\delta_{x^i+s}
\end{align*}
as the measure $\hat\mu$ has zero net mass.
Thus, $W_s$ is a function of positive type, and according to the classical
Bochner's theorem \cite[Theorem IX.9]{reed1975},
its Fourier transform $\hat W_s$ is a positive
and finite measure on $\mathbb R^d$.
On the other hand, denoting by $\hat V$, $\hat W_s$
the Fourier transforms of $V$, $W_s$ respectively,
we have
\[
\hat W_s (k) = 2 \bigl(1 - \cos (k \cdot s)\bigr) \hat V(k)
\]
in the sense of tempered distributions.
For every $k \neq 0$, we can find a non-zero $s \in \mathbb R^d$
such that the mapping $k' \mapsto 1 - \cos(k'\cdot s)$ is lower bounded away
from $0$ in a neighborhood of $k$.
Thus, in this neighborhood,
we have
\[
\hat V(k') = \frac{\hat W_s(k')}{2 \bigl( 1 - \cos(k' \cdot s) \bigr)}.
\]
Therefore, the distribution $\hat V$ restricted on $\mathbb R^d \setminus \{0\}$
is a positive and locally finite measure, which we denote by $\lambda$.
The difference $\hat V - \lambda$, being a Schwartz distribution,
is supported on the singleton $\{0\}$, and by the structure theorem
(see e.g., \cite[Théorème XXXV]{SchwartzDistributions}
and \cite[Theorem 2.3.4]{Hormander1}),
admits decomposition
\[
\hat V - \lambda = \sum_{\lvert n\rvert = 0}^m (-1)^{\lvert n\rvert}
c_n D^n \delta_0,
\]
$n$ being multi-indices, for some $m \in \mathbb N$ and $c_n \in \mathbb C$.
Denote the heat kernel by
\[
\rho^\varepsilon (x) = (2\pi\varepsilon)^{-d/2}
\exp(-\lvert x\rvert^2/2\varepsilon)
\]
and its Fourier transform reads
$\hat\rho^\varepsilon(k)
= (2\pi)^{-d/2} \exp (-2\pi^2\varepsilon\lvert k\rvert^2)$.
Define $V^\varepsilon = V \star \rho^\varepsilon$.
We then have
\begin{align*}
V^\varepsilon(0) = \langle \rho^\varepsilon, V \rangle
= \langle \hat \rho^\varepsilon, \hat V \rangle
&= \Bigl\langle \hat \rho^\varepsilon, \lambda
+ \sum_{\lvert n\rvert = 0}^m (-1)^{\lvert n\rvert}
c_n D^n \delta_0 \Bigr\rangle \\
&= \int_{\mathbb R^d \setminus \{0\}} \hat \rho^\varepsilon \dd \lambda
+ \frac{c_0}{(2\pi)^{d/2}}
+ \sum_{\lvert n\rvert = 1}^m c_n \nabla^n \hat \rho^\varepsilon (0),
\end{align*}
where $\langle\hat\rho^\varepsilon, \hat V\rangle$ is well defined,
since $\hat \rho^\varepsilon \in \mathcal S$ and $\hat V \in \mathcal S'$.
Thanks to the fact that
\[
\int_{\mathbb R^d \setminus \{0\}} \hat \rho^\varepsilon d \lambda
\nearrow \lambda(\mathbb R^d\setminus\{0\}),\qquad
V^\varepsilon(0) \to V(0),\qquad
\nabla^n \hat \rho^\varepsilon(0) \to 0
\]
when $\varepsilon \searrow 0$, for $n$ such that $\lvert n \rvert \geq 1$,
we can take the limit and obtain that the mass
$\lambda(\mathbb R^d\setminus\{0\})$ is finite and $c_0 \in \mathbb R$.
Then the original potential $V$ reads
\[
V(x) = \frac{1}{(2\pi)^{d/2}} \int_{\mathbb R^d \setminus \{0\}}
e^{ik \cdot x} \lambda (\dd k)
+ \frac{c_0}{(2\pi)^{d/2}} + P(x),
\]
where $P$ is an $m$-th-order polynomial with $P(0) = 0$.
The boundedness of $V$ implies that $P$ must be identically zero,
which concludes.

\paragraph{$\text{(iii)} \Rightarrow \text{(i)}$.}
Suppose (iii) holds.
Let $\mu$ be an arbitrary signed measure with $\int d\mu = 0$.
Then its Fourier transform $\hat\mu$ is even, real-valued,
belongs to the class $C_0$
and satisfies $\hat \mu(0) = 0$.
Thus, we have
\begin{align*}
\iint V(x - y) \mu(\dd x) \mu(\dd y)
= \langle V \star \mu, \mu \rangle
&= (2\pi)^{d/2} \langle \hat V \hat \mu, \hat \mu \rangle \\
&= (2\pi)^{d/2}
\int_{\mathbb R^d \setminus \{0\}} \bigl(\hat \mu(k)\bigr)^2 \hat V(\dd k)
\geq 0,
\end{align*}
which proves (ii). Finally, from the computation in the first paragraph,
we see that (i) is a consequence of (ii).
\end{proof}

\paragraph{Acknowledgements.}
We would like to thank an anonymous reviewer
whose comments motivated us to improve the quality of this paper.

\paragraph{Funding.}
The research of Zhenjie Ren is supported by the
\href{https://www.fime-lab.org}{FIME Research Initiative}.

\bibliography{ref}
\bibliographystyle{plain}

\end{document}